\documentclass[12 pt]{amsart}
\title{Murphy's Law for algebraic stacks}
\author{Daniel Bragg and Max Lieblich}
\usepackage{macros}

\begin{document}

\maketitle

\begin{abstract}
We show that various natural algebro-geometric moduli stacks, including the stack of curves, have the property that every Deligne--Mumford gerbe over a field appears as the residual gerbe of one of their points. These gerbes are universal obstructions for objects of the stack to be defined over their fields of moduli, and for the corresponding coarse moduli space to be fine. Thus, our results show that many natural moduli stacks hold objects that are obstructed from being defined over their fields of moduli in every possible way, and have coarse spaces which fail to be fine moduli spaces in every possible way. A basic insight enabling our arguments is that many classical constructions in equivariant projective geometry generalize to the setting of relative geometry over an arbitrary Deligne--Mumford gerbe over a field.
\end{abstract}

\setcounter{tocdepth}{1}
\tableofcontents

\section{Introduction}\label{sec:introduction}

In 1972, Shimura \cite{MR306215} showed that the complex genus 2 curve $C$ given by the smooth projective model of the affine curve
    \begin{equation}\label{eq:shimuras curve}
        y^2=x^6+ax^5+bx^4+x^3-\overline{b}x^2+\overline{a}x-1
    \end{equation}
for a general choice of complex numbers $a,b\in\mathbf{C}$ has a curious property: it is isomorphic to its complex conjugate, but is \emph{not\/} the base-change of a curve defined over the real numbers.
To place this example in a more general context, let us consider a field $K$ and a curve $C$ defined over the separable closure $K^s$ of $K$. If $K^s/L/K$ is an intermediate field extension such that $C$ is isomorphic to the base change of a curve defined over $L$, we say that $L$ is a \emph{field of definition} for $C$. The \emph{field of moduli} of $C$ is the fixed field of the group of automorphisms $\sigma\in\Gal(K^s/K)$ such that $C^{\sigma}\cong C$. The field of moduli is always contained in any field of definition, but Shimura's example shows that a curve need not be defined over its field of moduli.

Various authors have approached the question of how to quantify the difference between the field of moduli and the field of definition. For instance, D\`{e}bes and Emsalem \cite{MR1656571} showed how to attach to a curve a certain collection of cohomology classes, which have the property that the vanishing of at least one of these classes is necessary and sufficient for the curve to be defined over its field of moduli. In this paper, we take a modern approach to these types of questions centered on the \emph{gerbe of definition\/} of the curve $C$, which is the residual gerbe of the corresponding point of the moduli stack of curves. The gerbe of definition is a canonically defined geometric object which encodes the fields of moduli and of definition of the curve: the field of moduli is the residue field of the gerbe, and the fields of definition are the splitting fields of the gerbe. In particular, the statement that Shimura's curve fails to be defined over its field of moduli is equivalent to the statement that the gerbe of definition of $C$ is non-split. Thus, the moduli stack $\ms M_2$ of genus $2$ curves contains a non-split residual gerbe. Shimura's example has the further consequence of showing that the coarse moduli space $M_2$ of genus 2 curves fails to be a fine moduli space: the curve $C$ gives a point of $M_2$ with residue field $\mathbf{R}$ over which there does not exist a universal family. To assist the reader, we have included a more detailed discussion of this translation in Appendix \ref{sec:gerbes and stacks}, as well as the definition of residual gerbes and proofs of their key properties.
We remark that a related discussion of this collection of ideas in the stack theoretic context has recently been given by Bresciani and Vistoli \cite[\S3]{bresciani2022fields}.

In this paper, we solve the corresponding inverse problem: we show that \emph{every possible mode of failure\/}, that is, every Deligne--Mumford gerbe over a field, is realized in the stack of curves, as well as various other natural moduli stacks.

\subsection*{Notational convention}
In this paper, when we call an algebraic stack ``Deligne--Mumford'', we will be assuming that it is quasi-separated with separated diagonal. This was part of Deligne and Mumford's original definition, but the modern definition in the Stacks Project \cite[03YO]{stacks-project} no longer requires quasi-separatedness or a separated diagonal.

\subsection{Statement of the main results}

Let $\ms M$ be a Deligne--Mumford stack over a base field $K_0$. 

\begin{Definition}
We say that $\ms M$ \textit{satisfies Stacky Murphy's Law} if, for every finitely generated field extension $K/K_0$, every Deligne--Mumford gerbe $\ms G$ over $K$ is isomorphic to the residual gerbe of some point of $\ms M$. 
\end{Definition}

A moduli stack satisfying Stacky Murphy's Law parameterizes objects which are obstructed from being defined over their fields of moduli in every possible way. Said another way, the coarse moduli space of such a stack fails to be a fine moduli space in every possible way. This is an incarnation of Murphy's Law as described by Harris and by Vakil \cite{MR2227692}, but here concerning the stacky structure of a moduli problem rather than its singularities.
Our main result is the following.

\begin{Theorem}\label{thm:stacky murphy's law}
	Let $K_0$ be a field. The following moduli stacks satisfy Stacky Murphy's Law.
	\begin{enumerate}
		\item The stack
		\[
		\ms M_{\bullet}=\ms M_{0,3}\sqcup\ms M_{1,1}\sqcup\ms M_2\sqcup\ms M_3\sqcup\dots
		\]
		of smooth, proper, geometrically integral curves over $K_0$.
		\item The stack $\ms A_{\bullet}$ of principally polarized abelian varieties over $K_0$.	
        \item For any $d\geqslant 1$, the stack $\ms M^d_{\omega}$ of smooth projective varieties over $K_0$ of dimension $d$ with ample canonical bundle and reduced automorphism group scheme.
	\end{enumerate}
\end{Theorem}

The proof in the case of the moduli stack of curves consists of an explicit, but somewhat elaborate, construction which starts with a Deligne--Mumford gerbe $\ms G$ over a field $K$ and produces a relative curve $C\to\ms G$ with the property that the induced map $\ms G\to\ms M_{\bullet}$ identifies $\ms G$ with a residual gerbe of $\ms M_{\bullet}$. This construction has two main steps. First, we show the existence of a curve $C$ equipped with a morphism $C\to\ms G$ which is a relative curve. This is carried out in \S\ref{sec:complete intersections over gerbes}, and is a generalization of a classical construction of Serre \cite{MR0098097} of complete intersections equipped with free actions of a finite group. The second, and by far the most involved, part of the construction is to refine an initial choice of such a curve $C$ by taking a carefully constructed finite cover $C'\to C$. This is carried out in \S\ref{sec:curves over gerbes}, to which we refer the reader for further discussion of the strategy. We complete the proof of Theorem \ref{thm:stacky murphy's law} in \S\ref{sec:checking murphys law}. We first use the constructions of \S\ref{sec:curves over gerbes} to prove Theorem \ref{thm:stacky murphy's law} for the stack of curves, and then deduce the result for the stacks of abelian varieties and canonically polarized varieties using geometric constructions starting from curves.

\begin{Philosophy}
    While the geometric constructions in the situation at hand are involved, the primary observation enabling our results is simple to state: there are many classical results in equivariant projective geometry --- for example, constructions of varieties with prescribed automorphism groups --- which may be interpreted as describing geometric constructions over a split gerbe $\B G$; moreover, once rephrased in this way, many of these results seem to generalize to the setting of relative projective geometry of varieties over a (possibly nonsplit) gerbe. 
    This illustrates the conceptual advantage of working directly with gerbes as geometric objects, as opposed to working with their shadows in cohomology.
\end{Philosophy}

\begin{Remark}
     Our constructions for curves in fact show the stronger result that for any fixed integer $g$ the moduli stack $\ms M_{\geqslant g}$ of curves of genus $\geqslant g$ satisfies Stacky Murphy's Law. Similarly, our results show that for any fixed $g$ the stack $\ms A_{\geqslant g}$ of principally polarized abelian varieties of dimension $\geqslant g$ satisfies Stacky Murphy's Law.
\end{Remark}

\begin{Remark}
    Given a Deligne--Mumford gerbe $\ms G$ over $K_0$, we can ask for the minimal $g$ so that $\ms G$ is isomorphic to a residual gerbe of the stack $\ms M_g$ of curves of genus $g$. It is an interesting problem to bound this quantity in terms of invariants of $\ms G$, for instance, the \emph{index} of $\ms G$ (the g.c.d. of the degrees of field extensions splitting $\ms G$), the minimal rank $\rdim(\ms G)$ of a faithful vector bundle over $\ms G$, and the degree of the inertia of $G$. The bound implicit in our construction is at least doubly exponential in $\rdim(\ms G)$.
\end{Remark}

\begin{Remark}
    The moduli stacks in Theorem \ref{thm:stacky murphy's law} are all non quasi-compact, and have irreducible components of arbitrarily large dimension. These properties are necessary. For instance, a Deligne--Mumford stack satisfying Stacky Murphy's Law contains points with residue fields of arbitrarily large transcendence degree over the ground field. In fact, even a stack satisfying Stacky Murphy's Law only for gerbes over the ground field must contain points whose stabilizers are arbitrary finite groups, and thus cannot be quasi-compact.
\end{Remark}

\subsection{Consequences of Murphy's Law}\label{ssec:further results}

We explain some consequences of Theorem \ref{thm:stacky murphy's law}. A stack $\ms M$ satisfying Stacky Murphy's Law in particular contains as a residual gerbe any gerbe of the form $\B G$, where $K/K_0$ is a finitely generated field extension and $G$ is a finite \'{e}tale group scheme over $K$. A point of $\ms M$ whose residual gerbe is isomorphic to $\B G$ gives rise to a $K$-point of $\ms M$ whose automorphism group scheme is isomorphic to $G$ and which is not defined over any proper subfield of $K$. Thus, we obtain the following.

\begin{Theorem}\label{thm:Murphy's law for BG}
    Let $K_0$ be a field and let $\ms M$ be a Deligne--Mumford stack over $K_0$ which satisfies Stacky Murphy's Law (eg. any of the moduli stacks over $K_0$ in the statement of Theorem \ref{thm:stacky murphy's law}). If $K/K_0$ is a finitely generated field extension and $G$ is a finite \'{e}tale group scheme over $K$, then there exists a $K$-point $x_K\in\ms M(K)$ such that $\sAut_K(x_K)\cong G$ and such that $x_K$ is not defined over any intermediate subfield $K/L/K_0$ with $K\neq L$.
\end{Theorem}

In light of this result, one might view the problem of realizing all Deligne--Mumford gerbes as residual gerbes of a moduli stack as a natural generalization of the problem of constructing objects of the moduli stack with prescribed automorphism groups and fields of definition. Theorem \ref{thm:stacky murphy's law} thus gives in particular a solution to this problem for any of the listed moduli stacks. We highlight the following result for curves, which incorporates some slightly stronger properties which follow from our constructions. A similar existence statement holds for the objects parameterized by the other moduli stacks in Theorem \ref{thm:stacky murphy's law}. 

\begin{Corollary}\label{cor:presecribed automorphisms of curves 1}
     Let $K/K_0$ be a finitely generated extension of fields and let $G$ be a finite \'{e}tale group scheme over $K$. There exists a curve $C$ over $K$ equipped with a free $G$-action such that
     \begin{enumerate}
        \item $\sAut_K(C)\cong G$,
        \item $\sAut_K(C/G)=1$, and
        \item neither $C$ nor $C/G$ are defined over any intermediate subfield $K/L/K_0$ with $K\neq L$.
     \end{enumerate}
     Moreover, we may find such curves with arbitrarily large genus.
\end{Corollary}

In particular, this result implies that over an arbitrary field there exist curves whose automorphism group scheme is isomorphic to any prescribed finite \'{e}tale group scheme. This recovers a recent result of the first author \cite{bragg2023automorphism}.

Specializing further to the case of a trivial gerbe, we obtain the following result, which shows that the coarse moduli spaces of any of the moduli stacks in Theorem \ref{thm:stacky murphy's law} contain points with arbitrary residue fields.

\begin{Theorem}\label{thm:residue fields of M_g}
  Let $K_0$ be a field, let $\ms M$ be a Deligne--Mumford stack over $K_0$ with finite diagonal that satisfies Stacky Murphy's Law (eg. any of the moduli stacks over $K_0$ in the statement of Theorem \ref{thm:stacky murphy's law}), and let $M$ be the coarse moduli space of $\ms M$. If $K/K_0$ is a finitely generated field extension, then $K$ is isomorphic to the residue field of a point of $M$.
\end{Theorem}

We note that Theorems \ref{thm:Murphy's law for BG} and \ref{thm:residue fields of M_g} only use the existence of arbitrary \emph{split} residual gerbes in a stack satisfying Stacky Murphy's Law. However, such a stack has the much stronger property of containing arbitrary non-split residual gerbes. We explain a consequence of this for the possible ranks of certain vector bundles. Let $\ms M$ be a Deligne--Mumford stack over a field $K_0$. A vector bundle $\ms V$ on $\ms M$ is \emph{faithful} if the inertial action on $\ms V$ is faithful (Definition \ref{def:faithful vector bundle}). 
We are interested in determining bounds on the minimal rank of a faithful vector bundle on $\ms M$. Let
\[
    \ms M=\bigsqcup_{i\in I}\ms M_i
\]
be the decomposition of $\ms M$ into its connected components. For each $i\in I$, let $r_i$ be the minimal rank of a faithful vector bundle on $\ms M_i$, or $\infty$ if no such vector bundle exists, and set 
\[
    \rdim(\ms M)=\sup_{i\in I}r_i.
\]

\begin{Example}
    If $G$ is a finite group then $\rdim(\B G)$ is equal to the minimal dimension of a faithful $G$-representation over $K_0$. In the literature, this quantity is known as the \textit{representation dimension} of $G$.
\end{Example}

An immediate lower bound for $\rdim(\ms M)$ is provided by the observation that the rank of a faithful vector bundle at a geometric point $x_{\overline{K}}$ of $\ms M$ must be greater than or equal to the representation dimension of the stabilizer group $\Aut_{\ms M(\overline{K})}(x_{\overline{K}})$ of $\ms M$ at $x_{\overline{K}}$. It is known that the representation dimension of finite groups can be arbitrarily large (see eg. \cite{MR2818953}). As a consequence, if one knows that every finite group acts faithfully on some object of $\ms M$ defined over an algebraically closed extension of $K_0$, then one concludes that a faithful vector bundle defined on all of $\ms M$ must have unbounded rank, and hence $\rdim(\ms M)=\infty$. In particular, if $\ms M$ satisfies Stacky Murphy's Law, then $\rdim(\ms M)=\infty$. To obtain an invariant sensitive to the potentially non-split residual gerbes of $\ms M$, we restrict our attention to certain substacks. 
Suppose that $\ms M$ has finite diagonal. Given a finite group $G$, we let $\ms M^G\subset\ms M$ denote the locally closed substack consisting of those geometric points $x$ of $\ms M$ which have automorphism group isomorphic to $G$.

\begin{Theorem}\label{thm:ranks of faithful vector bundles}
    Let $\ms M$ be a Deligne--Mumford stack over $K_0$ with finite diagonal which satisfies Stacky Murphy's Law (eg. any of the moduli stacks in Theorem \ref{thm:stacky murphy's law}). If $G$ is a finite group containing a nontrivial central element of order coprime to the characteristic of $K_0$, then $\rdim(\ms M^G)=\infty$.
\end{Theorem}

\begin{Remark}
    While our formulation of Murphy's Law involves only Deligne--Mumford gerbes, one can also consider more general classes of gerbes. It would be interesting to find examples of geometrically natural moduli stacks which contain all gerbes with affine inertia as residual gerbes, or even all algebraic gerbes. By Theorem \ref{thm:Murphy's law for BG}, this includes as a special case the problem of constructing objects with prescribed automorphism group schemes. A recent result in this direction, due to Mathieu Florence \cite{florence2023realisation}, is that every linear algebraic group over an arbitrary field is the automorphism group of a smooth projective variety. Brion and Schr\"{o}er \cite{brion2022inverse} have also recently shown that every connected algebraic group over an arbitrary field is the connected component of the automorphism group of a smooth projective variety. We refer to the introduction of loc. cit. for a discussion of the history of such construction problems.
\end{Remark}

\subsection{Some explicit examples}\label{ssec:some explicit examples}

We give some explicit examples of residual gerbes of points of the moduli stack of curves of genus $g$ for small $g$. As $\ms M_{0,3}\cong\Spec K_0$ is a single non-stacky point, the first interesting case is genus 1.

\begin{Example}
    Let $\ms M_{1,1}$ be the stack of elliptic curves over a field $K_0$ of characteristic $\neq 2,3$. We will observe that every residual gerbe of $\ms M_{1,1}$ is split. Let $x\in|\ms M_{1,1}|$ be a topological point and let $k(x)$ and $\ms G(x)$ be the residue field and residual gerbe of $\ms M_{1,1}$ at $x$ (see Definitions \ref{def:topological points}, \ref{def:residue field}, and \ref{def:residual gerbe}) . We recall that the $j$-invariant map $j:\ms M_{1,1}\to\mathbf{A}^1\setminus\left\{0,1728\right\}$ induces an isomorphism $M_{1,1}\iso\mathbf{A}^1\setminus\left\{0,1728\right\}$. Thus, we obtain a diagram
\[
    \begin{tikzcd}
        \ms G(x)\arrow[hook]{r}\arrow{d}&\ms M_{1,1}\arrow{d}{j}\\
        \Spec k(x)\arrow[hook]{r}\arrow[dashed]{ur}\arrow[dashed, bend left=35]{u}&\mathbf{A}^1\setminus\left\{0,1728\right\}
    \end{tikzcd}
\]
of solid arrows. The morphism $\Spec k(x)\hookrightarrow\mathbf{A}^1\setminus\left\{0,1728\right\}$, which is isomorphic to the inclusion of the residue field of a point of $\mathbf{A}^1$, corresponds to a scalar say $j_x\in \mathbf{A}^1(k(x))=k(x)$. It is a classical fact that, for any field $K$ of characteristic $\neq 2,3$ and any scalar $j_0\in K$, there exists an elliptic curve over $K$ with $j$-invariant $j_0$. Thus, we may find an elliptic curve say $E$ over $k(x)$ with $j$--invariant $j(E)=j_x$. The curve $E$ gives rise to a dashed diagonal arrow filling in the diagram, which induces a splitting of the residual gerbe $\ms G(x)\to\Spec k(x)$. We therefore obtain an isomorphism $\ms G(x)\cong\B \sAut_{k(x)}(E)$. The problem of classifying the residual gerbes of $\ms M_{1,1}$ therefore reduces to classifying the possible automorphism group schemes of elliptic curves over residue fields of points of $\mathbf{A}^1\setminus\left\{0,1728\right\}$.
\end{Example}

\begin{Example}
    Let $\ms M_2$ be the moduli stack of genus 2 curves over $\mathbf{R}$ and let $x\in |\ms M_2|$ be the point classifying Shimura's curve $C$, given by the equation~\eqref{eq:shimuras curve} for a general choice of complex numbers $a,b\in\mathbf{C}$. The residue field of $x$ is isomorphic to $\mathbf{R}$, and together with the residual gerbe fits into a diagram
    \[
        \begin{tikzcd}
            \ms G(x)\arrow{d}\arrow[hook]{r}&\ms M_2\\
            \Spec\mathbf{R}.&
        \end{tikzcd}
    \]
    Using that $a$ and $b$ are general, one can show that the automorphism group $\Aut_{\mathbf{C}}(C)$ of $C$ over the complex numbers has order 2, generated by the hyperelliptic involution $\iota(x,y)=(x,-y)$. Thus, the base change of the residual gerbe $\ms G(x)$ to $\mathbf{C}$ is isomorphic to $\B\,(\mathbf{Z}/2)$, and so $\ms G(x)$ is classified by the unique nontrivial element of the cohomology group $\H^2(\Spec\mathbf{R},\mathbf{Z}/2)\cong\mathbf{Z}/2$.
    
    We can be more explicit. Consider the canonical short exact sequence
    \[
        0\to\Aut_{\mathbf{C}}(C)\to\Aut_{\mathbf{R}}(C)\xrightarrow{\beta}\Gal(\mathbf{C}/\mathbf{R})\to 0.
    \]
    One can show that the group $\Aut_{\mathbf{R}}(C)$ is isomorphic to $\mathbf{Z}/4$, generated by the $\mathbf{R}$-linear isomorphism $\tau:C\iso C$ which is the composition of the $\mathbf{C}$-linear isomorphism $C\iso\overline{C}$ defined by $(x,y)\mapsto (-x^{-1},ix^{-3}y)$ and the isomorphism $\overline{C}\iso C$ obtained by the base change of the complex conjugation map. One computes that $\tau^2=\iota$, so this generator indeed has order 4, and thus the above sequence is isomorphic to the unique extension
    \[
        0\to\mathbf{Z}/2\to\mathbf{Z}/4\xrightarrow{\beta}\mathbf{Z}/2\to 0.
    \]
    We will show that the residual gerbe $\ms G(x)$ of $\ms M_2$ at $x$ is isomorphic to the quotient stack $[\Spec\mathbf{C}/(\mathbf{Z}/4)]$, where $\mathbf{Z}/4=\Aut_{\mathbf{R}}(C)$ acts on $\Spec\mathbf{C}$ via $\beta$. To see this, consider the diagram
    \[
        \begin{tikzcd}
            C\arrow{d}\arrow{r}&\left[C/(\mathbf{Z}/2)\right]\arrow{r}\arrow{d}&\left[C/(\mathbf{Z}/4)\right]\arrow{d}{\pi}\\
            \Spec\mathbf{C}\arrow[equals]{dr}\arrow{r}&\B\,(\mathbf{Z}/2)\arrow{r}\arrow{d}&\left[\Spec\mathbf{C}/(\mathbf{Z}/4)\right]\arrow{d}\\
            &\Spec\mathbf{C}\arrow{r}&\Spec\mathbf{R}
        \end{tikzcd}
    \]
    in which the squares are Cartesian and we have used the canonical isomorphism $$[\Spec\mathbf{C}/(\mathbf{Z}/2)]=\Spec\mathbf{R}.$$ The map $\pi$ is a relative curve of genus 2, and so induces by descent a map $$[\Spec\mathbf{C}/(\mathbf{Z}/4)]\to\ms M_2,$$ which, by Proposition \ref{prop:gerbe --> stack}, induces an isomorphism $[\Spec\mathbf{C}/(\mathbf{Z}/4)]\iso\ms G(x)$.
\end{Example}

\begin{Example}\label{ex:double cover gerbe}
    Let $\ms M_2$ be the moduli stack of genus 2 curves over a field $K_0$ of characteristic $\neq 2$. We sketch a construction which produces nonsplit residual gerbes in $\ms M_2$ associated to index 2 elements of the Brauer groups of extensions of $K_0$. Let $K/K_0$ be a field extension and let $\alpha\in\Br(K)[2]$ be a 2-torsion Brauer class which has index 2.
    Let $P$ be a Brauer--Severi variety over $K$ of dimension 1 whose cohomology class maps to $\alpha$ under the (injective) boundary map
    \[
        \begin{tikzcd}
            \H^1(\Spec K,\PGL_2)\arrow[hook]{r}{\delta}&\Br(K)
        \end{tikzcd}
    \]
    in nonabelian cohomology. Let $\ms O_P(2)\defeq\omega_P^{\vee}$ be the anticanonical bundle of $P$. Write $\ms O_P(6)=\ms O_P(2)^{\otimes 3}$. We consider the stack $\ms G$ over $\Spec K$ whose fiber over a $K$-scheme $T$ is the groupoid of pairs $(\ms L,\varphi)$, where $\ms L$ is an invertible sheaf on $P_T$ and $\varphi:\ms L^{\otimes 2}\iso\ms O_P(6)_L$ is an isomorphism. The stack $\ms G$ has a canonical structure of $\mu_2$-gerbe over $K$, and one can show furthermore that the cohomology class $[\ms G]\in\H^2(\Spec K,\mu_2)$ maps to $\alpha\in\Br(K)$.
    
    We now describe a morphism $\ms G\to\ms M_2$. Over $P\times\ms G$ there is a universal square root of $\ms O_P(6)$, say $\ms L_u$. Choose a generic section $s\in\H^0(P\times\ms G,\ms L_u^{\otimes 2})=\H^0(P\times\ms G,\ms O_P(6)|_{\ms G})$, and let $C\to P\times\ms G$ be the branched double cover corresponding to the data $(\ms L_u,s)$. We obtain a diagram
    \[
        \begin{tikzcd}
            C\arrow{d}\arrow{r}&P\\
            \ms G.&
        \end{tikzcd}
    \]
    The morphism $C\to\ms G$ has geometric fiber isomorphic to a double cover of $\mathbf{P}^1$ branched over 6 points. Thus, $C\to\ms G$ is a family of genus 2 curves, and so corresponds to a morphism
    \[
        m_{C}:\ms G\to\ms M_2.
    \]
    If, for example, $K=K_0$ and $s$ is chosen generically, one can show that this morphism identifies $\ms G$ with a residual gerbe of a point of $\ms M_2$. We remark that this construction is closely related to Mestre's Brauer group obstruction \cite{MR1106431}.
\end{Example}

\begin{Remark}
    Cardona and Quer show \cite[Theorem 2]{MR2181874} that if $C/L$ is a genus 2 curve over a separably closed field of characteristic $\neq 2$ and the order of the automorphism group of $C$ is greater than $2$, then $C$ is defined over its field of moduli. Thus, if $x\in|\ms M_2|$ is a point whose automorphism group has order greater than 2, then the residual gerbe of $\ms M_2$ at $x$ is split.
\end{Remark}

\begin{Example}\label{ex:cyclic cover gerbe}
    The construction of Example \ref{ex:double cover gerbe} can be generalized by replacing the Brauer--Severi curve $P$ with a complete intersection curve in a Brauer--Severi variety. This produces for any extension $K/K_0$ and Brauer class $\alpha\in\Br(K)$ a map from a corresponding $\mu_n$-gerbe to the stack $\ms M_g$. In particular, this construction produces nonsplit residual gerbes in $\ms M_g$ for arbitrarily large $g$.
\end{Example}

\subsection{Organization of this paper}\label{ssec:organization of this paper}

In \S\ref{ssec:varieties over gerbes}, we define the \emph{inertial action} associated to a scheme over a stack. We give a brief recollection of gerbes over fields, and then discuss in more detail the inertial action on a scheme over a gerbe. In \S\ref{sec:complete intersections over gerbes} we show that over any gerbe there exist smooth and irreducible complete intersections of every dimension. In \S\ref{sec:curves over gerbes} we construct curves over gerbes with prescribed automorphism groups and fields of definition. This construction forms the main part of our proof of Theorem \ref{thm:stacky murphy's law}. In \S\ref{sec:checking murphys law} we give the proofs of the results stated in \S\ref{sec:introduction}. We first use our constructions from \S\ref{sec:curves over gerbes} to prove that the stack of curves satisfies Stacky Murphy's Law. Using various geometric constructions starting from curves, we then deduce the same for the other moduli stacks in Theorem \ref{thm:stacky murphy's law}. Finally, we give the proofs of the results stated in \S\ref{ssec:further results}.

Appendix \ref{sec:gerbes and stacks} consists of some background on algebraic stacks. We define the residual gerbe and residue field of a point of an algebraic stack and prove some results relating them to the field of moduli of a geometric point and to the coarse moduli space. This section contains no new results and is included for lack of a concise reference.

\subsection{Acknowledgments}

The authors thank Martin Olsson and Jakob Stix for interesting discussions, and Bjorn Poonen for some helpful correspondence.

\section{Gerbes over a field and the inertial action}\label{ssec:varieties over gerbes}

The main results of this paper require the construction of morphisms from a gerbe to various moduli stacks of varieties. Giving such a morphism is the same thing as constructing an appropriate variety over a gerbe. With this as motivation, we consider in this section a morphism from a scheme to a gerbe and the resulting \emph{inertial action}. We also give a brief account of some of the key facts about gerbes, and introduce some notation which we will continue throughout this paper.

\subsection{The inertial action}

A scheme equipped with a morphism to an algebraic stack acquires a canonical action by the inertia of the stack. This action will play a central role in the rest of this paper. Before giving the definition we introduce some notation.

\begin{Definition}\label{def:automorphism groups}
    If $\pi:X\to S$ is a morphism of schemes, we write $\Aut_S(X)$ for the group of automorphisms $\alpha$ of $X$ such that $\pi\circ\alpha=\pi$. Given another morphism $f:X\to Y$ of schemes, we let $\Aut_S(X/Y)$ denote the subgroup of $\Aut_S(X)$ consisting of those automorphisms $\alpha$ of $X$ such that $\pi\circ\alpha=\pi$ and $f\circ\alpha=f$. 
    (We note that in this definition the roles played by $S$ and $Y$ differ only psychologically. We will typically be thinking of $S$ as a base object.)
\end{Definition}

Let $S$ be a scheme, let $\ms S$ be an algebraic stack over $S$, and suppose given $S$-schemes $X$ and $Y$ and a diagram
 \[
    \begin{tikzcd}
        X\arrow{r}{f}\arrow{d}[swap]{\pi}&Y\\
        \ms S&
    \end{tikzcd}
 \]
over $S$.

\begin{Definition}\label{def:automorphism group scheaf}
     We define a sheaf of groups $\sAut_{\ms S}(X)$ on $\ms S$ by sending a morphism $T\to\ms S$ to the group $\Aut_T(X_T)$ of automorphisms of $X_T:=X\times_{\ms S}T$ over $T$.
We define a sub sheaf of groups $\sAut_{\ms S}(X/Y)\subset\sAut_{\ms S}(X)$ on $\ms G$ by sending a morphism $T\to\ms S$ to the subgroup $\Aut_T(X_T/Y)\subset\Aut_T(X_T)$ of automorphisms $\alpha$ of $X_T$ such that the diagram
\[
    \begin{tikzcd}
        X_T\arrow[dr,"\alpha","\sim"' sloped]\arrow[bend left=20]{drr}\arrow[bend right=20]{ddr}&&\\
        &X_T\arrow{r}{f\times\id_T}\arrow{d}{\pi_T}&Y\times T\\
        &T&
    \end{tikzcd}
\]
commutes.
\end{Definition}

We note that the sheaves $\sAut_{\ms S}(X)$ and $\sAut_{\ms S}(X/Y)$ are compatible with arbitrary base change on $\ms S$.

\begin{Definition}
    The \emph{inertia} of an algebraic stack $\ms S$ is the sheaf of groups $\ms I_{\ms S}$ on $\ms S$ whose value over a morphism $t:T\to\ms S$ is the group $\Aut_{\ms S(T)}(t)$ of automorphisms of $t$ in the groupoid $\ms S(T)$.

    Given a scheme $T$ and a morphism $t:T\to\ms S$, we write
    \[
        \sAut_{T}(t):=t^{-1}\ms I_{\ms S}
    \]
    for the pullback of $\ms I_{\ms S}$ along $t$. Explicitly, $\sAut_T(t)$ is the sheaf of groups on $T$ whose value on a morphism $u:T'\to T$ is the group $\Aut_{\ms S(T')}(t\circ u)$ of automorphisms of the object $t\circ u$ of the groupoid $\ms S(T')$.
\end{Definition}

\begin{Remark}
    Viewed as a stack over $\ms S$, the inertia $\ms I_{\ms S}$ fits into a canonical 2-Cartesian diagram
    \[
        \begin{tikzcd}
            \ms I_{\ms S}\arrow{d}\arrow{r}&\ms S\arrow{d}{\Delta_{\ms S}}\\
            \ms S\arrow{r}{\Delta_{\ms S}}&\ms S\times\ms S.
        \end{tikzcd}
    \]
    As a consequence, the morphism $\ms I_{\ms S}\to\ms S$ is always representable by algebraic spaces and is locally of finite type \cite[04XS]{stacks-project}.
\end{Remark}

We now define the canonical action of the inertia $\ms I_{\ms S}$ on $X$ over $\ms S$. Let $T$ be an $S$-scheme and let $t:T\to\ms S$ be an $S$-morphism. We form the 2-fiber product $X_T:=X\times_{\ms S}T$. Viewed as a $T$-scheme, $X_T$ has the following functorial description: given a $T$-scheme $u:T'\to T$, the set $X_T(T')$ of dashed arrows in the diagram
\[
    \begin{tikzcd}
        &X_T\arrow{d}{\pi_T}\arrow{r}&X\arrow{d}{\pi}\\
        T'\arrow{r}{u}\arrow[dashed]{ur}&T\arrow{r}{t}&\ms S
    \end{tikzcd}
\]
is identified with the set of pairs $(v,\varphi)$, where $v\in X(T')$ and $\varphi:\pi\circ v\iso t\circ u$ is an isomorphism in the groupoid $\ms S(T')$. We define a group homomorphism
\begin{equation}\label{eq:inertial action on T points}
    \Aut_{\ms S(T)}(t)\to\Aut_T(X_T)
\end{equation}
by sending an automorphism $\beta\in\Aut_{\ms S(T)}(t)$ to the automorphism of $X_T$ over $T$ defined at the level of the functor of points by
\[
    (v,\varphi)\mapsto (v,u^*(\beta)\circ\varphi).
\]

\begin{Definition}\label{def:inertial action}
    The \textit{inertial action map} for $X$ is the homomorphism
    \begin{equation}\label{eq:inertial action}
        \ms I_{\ms S}\to\sAut_{\ms S}(X)
    \end{equation}
    of sheaves of groups on $\ms S$ induced by the maps~\eqref{eq:inertial action on T points}. We note that the maps~\eqref{eq:inertial action on T points} factor through the subgroups $\Aut_T(X_T/Y)\subset\Aut_T(X_T)$, and hence the inertial action map~\eqref{eq:inertial action} factors through a map $\ms I_{\ms S}\to\sAut_{\ms S}(X/Y)$.
\end{Definition}

\subsection{Deligne--Mumford gerbes over a field}

In this section we recall the definition of a gerbe over a field and some fundamental results. General references for the theory are Giraud \cite{MR0344253} and Breen \cite{MR1301844}. Let $K$ be a field.

\begin{Definition}\label{def:gerbes}
    An algebraic stack $\ms G$ over $K$ is a \textit{gerbe} if
    \begin{enumerate}
        \item there exists a finite separable field extension $L/K$ such that $\ms G(L)$ is nonempty, and
        \item for any $K$-scheme $T$ and objects $x,y\in\ms G(T)$, there exists an fppf cover $T'\to T$ such that $x_{T'}\cong y_{T'}$ in the groupoid $\ms G(T)$.
    \end{enumerate}
\end{Definition}

\begin{Notation}
    Given a gerbe $\ms G$ over $K$, we will write $G_{\ms G}:=\ms I_{\ms G}$ for the inertia of $\ms G$. If $\ms G$ is Deligne--Mumford, we will write $|G_{\ms G}|$ for the degree of the morphism $G_{\ms G}\to\ms G$.
\end{Notation}

\begin{Remark}
    It need not be the case that $G_{\ms G}\to\ms G$ is the pullback of a group space defined over $K$, although this is automatic if $G_{\ms G}$ is commutative.
\end{Remark}

By an \emph{algebraic group} over a field $K$ we will mean a group algebraic space locally of finite type over $K$. 

\begin{Example}\label{ex:split gerbe}
    The fundamental example of a gerbe is the classifying stack $\B_K\!G$ associated to an algebraic group $G$ over $K$, which is defined as the quotient stack
    \[
        \B_K\!G:=\left[\Spec K/G\right]
    \]
    where $G$ acts on $\Spec K$ by the trivial action. If we think the ground field $K$ is clear from context, we may write simply $\B G$ for $\B_K\!G$. There is a canonical map $\rho:\B G\to\Spec K$ via which we regard $\B G$ as a stack over $K$. By the definition of the quotient stack, if $T$ is a $K$-scheme the fiber $\B G(T)$ is the groupoid of $G_T$-torsors over $T$. We note that the trivial $G$-torsor over $K$ gives an element of $\B G(K)$, and furthermore any $G_T$-torsor over a scheme $T$ is fppf locally trivial. Thus, $\B G$ is indeed a gerbe over $K$.
       
    The automorphism group of a $G_T$-torsor over a $K$-scheme $T$ is naturally identified with $G(T)$. These identifications give rise to a canonical isomorphism 
    \[
        G|_{\B G}=G\times_{\Spec K}\B G=G_{\B G},
    \]
    identifying the pullback of the group scheme $G$ along $\rho$ with the inertia $G_{\B G}$ of $\B G$. 
\end{Example}

\begin{Definition}\label{def:splitting of a gerbe}
    Let $L/K$ be a field extension. A \textit{splitting} of $\ms G$ over $L$ is an element of $\ms G(L)$, or in other words a morphism $s:\Spec L\to\ms G$ rendering the diagram
    \[
        \begin{tikzcd}[column sep=small]
            &\ms G\arrow{d}\\
            \Spec L\arrow{ur}{s}\arrow{r}&\Spec K
        \end{tikzcd}
    \]
    commutative. Equivalently, this is a section of the $L$-gerbe $\ms G_L:=\ms G\otimes_KL\to\Spec L$. We say that $L/K$ \textit{splits} $\ms G$ or that $\ms G$ is \textit{split} over $L$ if there exists a splitting of $\ms G$ over $L$.
\end{Definition}

If $G$ is an algebraic group over $K$, then the quotient map $s:\Spec K\to\B G$ gives a canonical splitting of $\B G$ over $K$. In fact, every split gerbe is of this form. Indeed, let $\ms G$ be a gerbe over $K$ and let $s:\Spec K\to\ms G$ be a splitting of $\ms G$ over $K$. Write $G:=\sAut_K(s)=s^{-1}G_{\ms G}$ for the pullback of the inertia $G_{\ms G}$ along $s$. The functor $\sIsom_{\ms G}(\Spec K,\_)$ defines an isomorphism
\[
    \ms G\iso\B G
\]
of stacks over $K$ which sends $s$ to the canonical splitting of $\B G$.

\begin{Remark}
    If $\ms G$ is any gerbe over $K$, then by definition $\ms G$ splits over some finite extension $L/K$, and so there exists an isomorphism $\ms G_L\cong\B_L\!G_L$ for some algebraic group $G_L$ over $L$. Thus, a gerbe over $K$ is a $K$-form of the classifying stack $\B_L\!G_L$ for some algebraic group $G_L$.
\end{Remark}

\begin{Remark}\label{rem:classification of gerbes}
    An important fact about gerbes is that they may be classified cohomologically. In the commutative case, this classification is easy to state. Let $G$ be a commutative algebraic group over $K$. A $G$-\textit{gerbe} over $K$ is a gerbe $\ms G$ over $K$ equipped with an isomorphism $G_{\ms G}\cong G|_{\ms G}$ of sheaves of groups over $\ms G$. The  set of isomorphism classes of $G$-gerbes over $K$ is in canonical bijection with the \'{e}tale cohomology group $\H^2(\Spec K,G)$, with the split $G$-gerbe $\B G$ corresponding to zero. The extension of this classification to gerbes with nonabelian inertia is more involved, and is the object of Giraud's theory of nonabelian cohomology \cite{MR0344253}.
\end{Remark}

We now consider Deligne--Mumford gerbes.
\begin{Lemma}\label{lem:split DM gerbe}
    Let $G$ be an algebraic group over $K$. The gerbe $\B G$ is Deligne--Mumford if and only if $G$ is a finite \'{e}tale group scheme.
\end{Lemma}
\begin{proof}
    The quotient morphism $s:\Spec K\to\B G$ is represented by $G$-torsors, and hence is an fppf cover. Similarly, the diagonal morphism $\Delta_{\B G}$ is also an fppf cover. Consider the diagram
    \[
        \begin{tikzcd}
            G\arrow{d}\arrow{r}&G_{\B G}\arrow{d}\arrow{r}&\B G\arrow{d}{\Delta_{\B G}}\\
            \Spec K\arrow{r}{s}&\B G\arrow{r}{\Delta_{\B G}}&\B G\times\B G
        \end{tikzcd}
    \]
    with 2-Cartesian squares. By definition, $\B G$ is Deligne--Mumford if and only if it admits an \'{e}tale cover by a scheme and the diagonal $\Delta_{\B G}$ is quasi-compact and separated, or equivalently \cite[06N3]{stacks-project} if and only if $\Delta_{\B G}$ is unramified, quasi-compact, and separated. These conditions may be checked after an fppf cover, so $\B G$ is Deligne--Mumford if and only if the morphism $G\to\Spec K$ satisfies the same conditions. As the target of this morphism is the spectrum of a field, these conditions hold if and only if $G\to\Spec K$ is finite \'{e}tale.
\end{proof}

\begin{Lemma}\label{lem:DM gerbe over a field}
    An algebraic stack $\ms G$ over $K$ is a Deligne--Mumford gerbe over $K$ if and only if there exists a finite separable field extension $L/K$ and a finite group $G$ such that $\ms G_L\cong\B_L\!G$ as stacks over $L$.
\end{Lemma}
\begin{proof}
     Let $\ms G$ be a Deligne--Mumford gerbe over $K$. By definition, we may find a finite separable extension $L/K$ splitting $\ms G$, and hence by the above discussion an isomorphism $\ms G_L\cong\B_L\!G_L$ over $L$ for some algebraic group $G_L$ over $L$. The stack $\B_L\!G_L$ is Deligne--Mumford, so by Lemma \ref{lem:split DM gerbe} $G_L$ is a finite \'{e}tale group scheme over $L$. Taking a further separable extension, we may arrange so that $G_L$ is the group scheme associated to a finite group. Conversely, the property of being a gerbe over $K$ may be checked after making a finite separable extension \cite[0CPR,0CPS]{stacks-project}, so any such $\ms G$ is a gerbe over $K$. Furthermore, the canonical splitting of the stack $\B_L\!G$ gives rise to a finite \'{e}tale cover $\Spec L\to\ms G$, so any such $\ms G$ is Deligne--Mumford.
\end{proof}

\subsection{The inertial action on a scheme over a gerbe}

We now consider in more detail the inertial action on a scheme equipped with a morphism to a gerbe over a field. We first consider the case of a split gerbe.

\begin{Example}\label{ex:split gerbe 2}
Let $G$ be an algebraic group over $K$ and consider the split gerbe $\B G$ over $K$. We consider the inertial action associated to the quotient map $s:\Spec K\to\B G$. As described in Example \ref{ex:split gerbe}, we have a canonical identification $G|_{\B G}=G_{\B G}$. Thus, the inertial action gives a map
    \begin{equation}\label{eq:canonical map for split gerbe}
        G|_{\B G}=G_{\B G}\to\sAut_{\B G}(\Spec K).
    \end{equation}
The morphism $s$ equipped with this map is the \emph{universal $G$-torsor}, in the following sense. For a $K$-scheme $T$, the definition of the quotient stack gives rise to an identification between the groupoid of morphisms $T\to\B G$ and the groupoid of left $G_T$-torsors $P\to T$ over $T$. This identification has the following description: given a morphism $t:T\to\B G$, we define $P$ as the 2-fiber product in the 2-Cartesian diagram
    \[
        \begin{tikzcd}
            P\arrow{d}\arrow{r}&\Spec K\arrow{d}{s}\\
            T\arrow{r}{t}&\B G.
        \end{tikzcd}
    \]
The pullback of the inertial action map~\eqref{eq:canonical map for split gerbe} to $T$ gives a map $G_T\to\sAut_T(P)$, which defines a left $G_T$-torsor structure on $P\to T$.   
\end{Example}

\begin{Example}\label{ex:non split gerbe}
Now consider a possibly non-split gerbe $\ms G$ over $K$ with inertia $G_{\ms G}:=\ms I_{\ms G}$. Let $X$ be a $K$-scheme equipped with a $K$-morphism $\pi:X\to\ms G$. Suppose given a field extension $L/K$ and a splitting $s\in\ms G(L)$ of $\ms G$ over $L$. Write $G_L=G_{\ms G}\times_{\ms G,s}\Spec L$ for the pullback of $G_{\ms G}$ along $s$. This is an algebraic group over $L$. 
The morphism $s$ gives rise to a splitting $s'\in\ms G_L(L)$ of the gerbe $\ms G_L$ over $L$, which induces an isomorphism $\ms G_L\cong\B G_L$ of gerbes over $L$. We set $X_L=X\otimes_KL$ and $\widetilde{X}_L=X\times_{\ms G,s}\Spec L$, and so obtain a diagram
\begin{equation}\label{eq:non split gerbe diagram}
    \begin{tikzcd}
        \widetilde{X}_L\arrow{d}\arrow{r}&X_L\arrow{r}\arrow{d}&X\arrow{d}{\pi}\\
        \Spec L\arrow{r}{s'}\arrow[equals]{dr}&\B G_L\arrow{r}\arrow{d}&\ms G\arrow{d}\\
        &\Spec L\arrow{r}&\Spec K
    \end{tikzcd}
\end{equation}
with 2-Cartesian squares, where the horizontal composition $\Spec L\to\ms G$ is the map $s$. Pulling back the inertial action map
\begin{equation}\label{eq:canonical map for split gerbe part 3}
        G_{\ms G}\to\sAut_{\ms G}(X)
\end{equation}
along $s$ gives a map
\[
     G_L\to\sAut_{L}(\widetilde{X}_L)
\]
of sheaves of groups over $L$, which corresponds to a left action of $G_L$ on $\widetilde{X}_L$. As described in Example \ref{ex:split gerbe 2}, this action makes $\widetilde{X}_L$ into a $G_L$-torsor over $X_L$. In particular, this action is free,~\eqref{eq:canonical map for split gerbe part 3} is injective, and the morphism $\widetilde{X}_L\to X_L$ induces an isomorphism $[\widetilde{X}_L/ G_L]\iso X_L$.
\end{Example}

\begin{Example}\label{ex:non split gerbe, once again}
We will have occasion to consider the following situation. With the notation of Example \ref{ex:non split gerbe}, suppose given also a $K$-scheme $Y$ and a $K$-morphism $f:X\to Y$, so that we have a diagram
\[
    \begin{tikzcd}
        X\arrow{r}{f}\arrow{d}[swap]{\pi}&Y\\
        \ms G&
    \end{tikzcd}
\]
over $K$. Set $Y_L=Y\otimes_KL$. After pullback along the splitting $s\in\ms G(L)$, we obtain a diagram
\[
    \begin{tikzcd}
        \widetilde{X}_L\arrow{d}\arrow{r}&X_L\arrow{d}\arrow{r}{f_L}&Y_L\\
        \Spec L\arrow{r}{s'}&\B G_L&
    \end{tikzcd}
\]
of $L$-stacks, in which the square is 2-Cartesian and $\widetilde{X}_L\to X_L$ is a $G_L$-torsor. The inertial action map~\eqref{eq:canonical map for split gerbe part 3} factors through $\sAut_{\ms G}(X/Y)$, yielding injections
\begin{equation}\label{eq:inertial actions inclusion}
        G_{\ms G}\hookrightarrow\sAut_{\ms G}(X/Y)\hookrightarrow\sAut_{\ms G}(X).
\end{equation}
Pulling back these maps along $s$, we obtain injections
\[
     G_L\hookrightarrow\sAut_L(\widetilde{X}_L/Y_L)\hookrightarrow\sAut_L(\widetilde{X}_L).
\]
\end{Example}

\begin{Remark}
    Let $X$ be a $K$-scheme. Given a $K$-morphism $X\to\ms G$, the above discussion shows that there is a finite separable field extension $L/K$, a group scheme $ G_L$, and an $L$-scheme $\widetilde{X}_L$ equipped with a $ G_L$-action such that $X_L\cong\widetilde{X}_L/ G_L$. Thus, $X$ is a twisted form of a quotient $\widetilde{X}_L/ G_L$, the twisting being encoded in the structure of the stack $\ms G$. We remark that in general neither the variety $\widetilde{X}_L$ nor the group scheme $ G_L$ need be defined over $K$, and even if they are, the $G_L$-action need not descend to a group action over $K$.
\end{Remark}

\section{Complete intersections over gerbes}\label{sec:complete intersections over gerbes}

Let $K$ be a field and let $\ms G$ be a Deligne--Mumford gerbe over $K$ with inertia $G_{\ms G}:=\ms I_{\ms G}$. In this section we will show that there exist smooth complete intersections over $\ms G$ of every dimension. This generalizes Serre's construction of smooth complete intersections equipped with a free action of a finite group \cite{MR0098097}.

We begin by proving a certain stacky Bertini theorem. Let $\ms X$ be a Deligne--Mumford stack over $K$ with finite diagonal and coarse moduli space $\rho:\ms X\to X$. Let $\pi:\ms X\to\ms G$ be a $K$-morphism, so that we have a diagram
\[
    \begin{tikzcd}
        \ms X\arrow{d}[swap]{\pi}\arrow{r}{\rho}&X\\
        \ms G.&
    \end{tikzcd}
\]

\begin{Lemma}\label{lem:pullback from coarse space}
    If $\ms L$ is an invertible sheaf on $X$, then the pullback map
    \[
        \rho^*:\H^0(X,\ms L)\to\H^0(\ms X,\rho^*\ms L)
    \]
    is an isomorphism.
\end{Lemma}
\begin{proof}
    The pullback map $\rho^{\#}:\ms O_X\iso\rho_*\ms O_{\ms X}$ is an isomorphism (see Remark \ref{rem:coarse moduli space remark}). The composition
    \[
        \begin{tikzcd}
            \ms L\arrow{r}{\sim}[swap]{\id_{\ms L}\otimes\rho^{\#}}&\ms L\otimes\rho_*\ms O_{\ms X}\arrow{r}{\sim}&\rho_*\rho^*\ms L
        \end{tikzcd}
    \]
    is equal to the adjunction map $\ms L\to\rho_*\rho^*\ms L$, where the second map is the isomorphism appearing in the projection formula. The pullback map $\rho^*:\H^0(X,\ms L)\to\H^0(\ms X,\rho^*\ms L)$ is the map on global sections induced by the adjunction map, and therefore is an isomorphism.
\end{proof}

\begin{Definition}
    An invertible sheaf $\ms L$ on $\ms X$ is \textit{(very) ample} if it is the pullback of a (very) ample invertible sheaf on $X$. 
\end{Definition}

\begin{Theorem}\label{thm:Bertini for stacks}
    Suppose that the morphism $\pi:\ms X\to\ms G$ is smooth and representable and has connected geometric fibers of dimension $d\geqslant 2$. Let $\ms L$ be a very ample invertible sheaf on $\ms X$. If $s\in\H^0(\ms X,\ms L)$ is a general section, then the closed substack $Z=V(s)\subset\ms X$ has the property that the map $Z\to\ms G$ is smooth and representable and has connected geometric fibers of dimension $d-1$.
\end{Theorem}
\begin{proof}
    Say $\ms L=\rho^*\ms L_0$ where $\ms L_0$ is a very ample invertible sheaf on $X$. Choose a finite separable extension $L/K$ and a splitting $s\in\ms G(L)$ of $\ms G$. As in Example \ref{ex:non split gerbe}, we write $G_L=s^{-1}G_{\ms G}$ for the pullback of $G_{\ms G}$ along $s$ and set $X_L=X\otimes_KL$, $\ms X_L=\ms X\otimes_KL$, and $\widetilde{X}_L=\ms X\times_{\ms G,s}\Spec L$. We obtain a diagram
    \[
        \begin{tikzcd}
            \widetilde{\ms X}_L\arrow{r}{q_L}\arrow[bend left=25]{rr}{\widetilde{\rho}_L}\arrow{d}&\ms X_L\arrow{r}{\rho_L}\arrow{d}&X_L\\
            \Spec L\arrow{r}&\B G_L&
        \end{tikzcd}
    \]
    of $L$-stacks with 2-Cartesian square. It follows from our assumptions on $\pi$ that $\widetilde{\ms X}_L$ is a scheme and is smooth and geometrically connected over $L$ of dimension $d$. Furthermore, $\widetilde{X}_L$ is equipped with an action of $G_L$, and the morphism $q_L:\widetilde{\ms X}_L\to\ms X_L$ induces an isomorphism $[\widetilde{\ms X}_L/G_L]\iso\ms X_L$, and in particular is finite \'{e}tale of degree $|G_{\ms G}|$. Finally, $\ms X_L$ is a Deligne--Mumford stack and $\rho_L:\ms X_L\to X_L$ is a coarse moduli space morphism. Let $\ms L_L$ be the base change of $\ms L$ to $\ms X_L$ and let $\ms L_{0,L}$ be the pullback of $\ms L_0$ to $X_L$. Let $V\subset\H^0(\widetilde{\ms X}_L,q_L^*\ms L_L)$ be the image of the pullback map
    \[
        \widetilde{\rho}_L^*:\H^0(X_L,\ms L_{0,L})\to\H^0(\widetilde{\ms X}_L,q_L^*\ms L_L).
    \]
    As $\ms L_{0,L}$ is very ample, the linear system $V$ is base point free and is separable. We may therefore  apply Spreafico's Bertini theorem \cite{MR1612610} to conclude that the subscheme of $\widetilde{\ms X}_L$ cut out by a general element of $V$ is a smooth and geometrically connected $L$-scheme of dimension $d-1$. On the other hand, the pullback map $\widetilde{\rho}_L^*$ factors as the composition
    \[
    \begin{tikzcd}[column sep=small]
        \H^0(X_L,\ms L_{0,L})\arrow{r}{\sim}[swap]{\rho_L^*}&\H^0(\ms X_L,\ms L_L)\arrow[hook]{r}[swap]{q_L^*}&\H^0(\widetilde{\ms X}_L,q_L^*\ms L_L).
    \end{tikzcd}
    \]
    Here, the first map is an isomorphism by Lemma \ref{lem:pullback from coarse space} and the second map is injective because $q_L$ is flat. Thus, pullback along $q_L$ gives an isomorphism
    \[
        \begin{tikzcd}[column sep=small]
            \H^0(\ms X,\ms L)\otimes_KL=\H^0(\ms X_L,\ms L_{L})\arrow{r}{\sim}[swap]{q^*_L}& V.
        \end{tikzcd}
    \]
    It follows that a general element of $\H^0(\ms X,\ms L)$ gives rise to a general element of $V$, which proves the result.
\end{proof}

\begin{Definition}\label{def:faithful vector bundle}
    We say that a vector bundle $\ms V$ on an algebraic stack $\ms S$ is \emph{faithful} if the inertial action map $\ms I_{\ms S}\to\sAut_{\ms S}(\ms V)$ is injective. A locally free coherent sheaf on $\ms S$ is \emph{faithful} if the associated vector bundle is faithful.
\end{Definition}

\begin{Lemma}\label{lem:there is a faithful locally free sheaf}
    There exists a faithful locally free coherent sheaf on $\ms G$.
\end{Lemma}
\begin{proof}
    Choose a finite field extension $L/K$ which splits $\ms G$ and let $s\in\ms G(L)$ be a section, so that we have a diagram
    \[
        \begin{tikzcd}[column sep=small]
            &\ms G\arrow{d}\\
            \Spec L\arrow{ur}{s}\arrow{r}&\Spec K.
        \end{tikzcd}
    \]
    The pushforward $\ms E:=s_*\ms O_{\Spec L}$ is a faithful locally free sheaf on $\ms G$. We remark that $s$ is finite \'{e}tale of degree $|G_{\ms G}|[L:K]$ (see eg. the lower two rows of the diagram~\eqref{eq:non split gerbe diagram}), so $\ms E$ has rank $|G_{\ms G}|[L:K]$.
\end{proof}

Let $\ms E$ be a faithful locally free coherent sheaf on $\ms G$. Let $\pi:\ms P:=\mathbf{P}_{\ms G}(\ms E)\to\ms G$ be the associated projective bundle. Then $\ms P$ is itself a separated Deligne--Mumford stack over $K$, and hence admits a coarse moduli space $\rho:\ms P\to P$. We obtain a diagram
\[
    \begin{tikzcd}
        \ms P\arrow{d}[swap]{\pi}\arrow{r}{\rho}&P\\
        \ms G.&
    \end{tikzcd}
\]
As a projective bundle, $\ms P$ comes equipped with the tautological invertible sheaf $\ms O_{\ms P}(1)$.
    
\begin{Remark}\label{rem:quotient of projective space}
    To orient the reader, we record what our constructions so far amount to when $\ms G=\B G$ is the classifying stack associated to a finite group $G$. In this case, pulling back $\ms E$ along the canonical section yields a finite dimensional $K$-vector space, say $V$, equipped with a faithful $G$-action. We have $\ms P=\mathbf{P}_{\B G}(\ms E)\cong [\mathbf{P}(V)/G]$, and the coarse moduli space map $\rho:\ms P\to P$ is identified with the map $[\mathbf{P}(V)/G]\to\mathbf{P}(V)/G$. Thus, we have a diagram
    \[
        \begin{tikzcd}
            \mathbf{P}(V)\arrow{d}\arrow{r}&\left[\mathbf{P}(V)/G\right]\arrow{d}{\pi}\arrow{r}{\rho}&\mathbf{P}(V)/G\\
            \Spec K\arrow{r}&\B G&
        \end{tikzcd}
    \]
    where the square is 2-Cartesian. In the general case, where $\ms G$ is not necessarily split, we obtain this picture after making a base change along a finite separable extension $L/K$ splitting $\ms G$ and the inertia.
\end{Remark}

\begin{Lemma}\label{lem:Pic is torsion}
    The Picard group $\Pic(\ms G)$ is torsion.
\end{Lemma}
\begin{proof}
    We will use the fact that, if $f:\ms X\to\ms Y$ is a finite locally free morphism of algebraic stacks of degree $e$, then there exists a canonical \emph{norm map} for $f$, which is a group homomorphism
    \[
        N_f:\Pic(\ms X)\to\Pic(\ms Y)
    \]
    with the property that $N_f([f^*\ms L])=[\ms L^{\otimes e}]$ for any invertible sheaf $\ms L$ on $\ms Y$. The existence of such a map in the case when $\ms X$ and $\ms Y$ are schemes is given in the Stacks Project \cite[0BD2,0BCY]{stacks-project}, and the existence in general follows from this case by descent.
    
    Choose a splitting $s\in\ms G(L)$ over a finite separable extension $L/K$. Then the morphism $s:\Spec L\to\ms G$ is finite \'{e}tale of degree $e:=|G_{\ms G}|[L:K]$. For any invertible sheaf $\ms L$ on $\ms G$ we have that $N_s([s^*\ms L])=[\ms L^{\otimes e}]$. As the Picard group of $\Spec L$ is trivial, this implies that $\ms L^{\otimes e}$ is trivial, so $\Pic(\ms G)$ is $e$-torsion.
\end{proof}

\begin{Remark}
    If $G$ is an algebraic group over $K$ then the Picard group $\Pic(\B G)$ may be identified with the group of characters of $G$.
\end{Remark}

\begin{Lemma}\label{lem:the tautological sheaf descends}
    For all sufficiently divisible positive integers $m$, the invertible sheaf $\ms O_{\ms P}(m)$ on $\ms P$ is very ample.
\end{Lemma}
\begin{proof}
    The coarse moduli space $P$ is a $K$-form of the quotient space $\mathbf{P}(V)/G$ for some finite group $G$ acting on a finite dimensional vector space $V$. In particular, $P$ is projective, so we may find a very ample invertible sheaf say $\ms L_0$ on $P$. The usual formula for the Picard group of a projective bundle gives
    \[
        \Pic(\ms P)=\mathbf{Z}[\ms O_{\ms P}(1)]\oplus\Pic(\ms G).
    \]
    By Lemma \ref{lem:Pic is torsion}, the Picard group $\Pic(\ms G)$ is torsion. Thus, a power of $\rho^*\ms L_0$ is isomorphic to $\ms O_{\ms P}(m)$ for some integer $m$. 
\end{proof}

\begin{Definition}
    Let $\ms X$ be a Deligne--Mumford stack over $K$ with finite diagonal. Then the inertia stack $\ms I_{\ms X}$ of $\ms X$ is finite and unramified over $\ms X$. It follows that the identity section $e:\ms X\to\ms I_{\ms X}$ is open and closed, and therefore the map $|\ms I_{\ms X}\setminus e(\ms X)|\to|\ms X|$ has closed image. The \textit{stacky locus} of $\ms X$ is the unique reduced closed substack $\Gamma_{\ms X}\subset\ms X$ such that $|\Gamma_{\ms X}|$ is equal to the image of the map $|\ms I_{\ms X}\setminus e(\ms X)|\to|\ms X|$.
\end{Definition}

We will show that we can choose the faithful locally free sheaf $\ms E$ on $\ms G$ so that the codimension of $\Gamma_{\ms P}$ in $\ms P$ is as large as we like. We make the following definition.

\begin{Definition}
    Let $\ms E$ be a locally free coherent sheaf on an algebraic stack $\ms X$. A locally free sheaf $\ms F$ on $\ms X$ is said to be a \textit{polynomial tensorial construction in }$\ms E$ if
    \[
        \ms F\cong\bigoplus_{i=1}^N\ms E^{\otimes l_i}
    \]
    for some nonnegative integers $l_1,\dots,l_N$.
\end{Definition}

\begin{Proposition}\label{prop:tensorial construction}
    Let $\ms E$ be a faithful locally free coherent sheaf on $\ms G$. For any integer $d$, there exists a faithful locally free sheaf $\ms F$ on $\ms G$ which is a polynomial tensorial construction in $\ms E$ such that
    \begin{equation}\label{eq:inequality 2}
        \dim\mathbf{P}_{\ms G}(\ms F)-\dim\Gamma_{\mathbf{P}_{\ms G}(\ms F)}\geqslant d.
    \end{equation}
\end{Proposition}
\begin{proof}
    It will suffice to prove the result when $K$ is algebraically closed and $\ms G$ is the split gerbe $\B G$ associated to a finite group $G$. As described in Remark \ref{rem:quotient of projective space}, in this case $\ms E$ gives rise to a finite dimension $K$-vector space say $V$ equipped with a faithful action of $G$, corresponding to an injective group homomorphism $\varphi:G\hookrightarrow\GL(V)$. For an element $g\in G$, let $\Gamma_{\varphi(g)}\subset\mathbf{P}(V)$ be the subvariety of points $x\in\mathbf{P}(V)$ such that $g\cdot x=x$. Explicitly, $\Gamma_{\varphi(g)}$ is the union of the projectivizations of the eigenspaces of the linear operator $\varphi(g)$. Let $\Gamma_{\varphi}\subset\mathbf{P}(V)$ be the union of the $\Gamma_{\varphi(g)}$ as $g$ ranges over the nontrivial elements of $G$. The claimed inequality is then equivalent to the inequality
    \begin{equation}\label{eq:goal inequality}
        \dim\mathbf{P}(V)-\dim\Gamma_{\varphi}\geqslant d.
    \end{equation}
    We will first show that this can be achieved for $d=1$ by taking a suitable polynomial tensorial construction. Consider the intersection of the subgroup $\varphi(G)\subset\GL(V)$ with the subgroup $\ms Z(\GL(V))=\mathbf{G}_m\cdot\mathrm{I}$ of scalar matrices. As $G$ is finite, this intersection is contained in the subgroup $\mu_m\cdot\mathrm{I}$ for some integer $m$. The polynomial tensorial construction $K\oplus V\oplus V^{\otimes 2}\oplus\dots\oplus V^{\otimes m-1}$ equipped with the induced $G$-action then has the property that no element of $g$ acts by a scalar matrix, and so the above inequality~\eqref{eq:goal inequality} holds with $d=1$.

    We now assume that the inequality~\eqref{eq:goal inequality} holds for $V$ with $d=1$, and show that by taking a further polynomial tensorial construction it can also be made to hold for any given $d\geqslant 1$. For an element $g\in G$, we have that the dimension of $\Gamma_{\varphi(g)}$ at any point is less than or equal to $r_{\varphi(g)}-1$, where $r_{\varphi(g)}$ is the largest dimension of an eigenspace of $\varphi(g)$. Consider the polynomial tensorial construction $V^{\oplus d}$ equipped with the direct sum representation
    \[
        \psi:G\hookrightarrow\GL(V^{\oplus d}).
    \]
    We note that $r_{\psi(g)}=dr_{\varphi(g)}$. Taking the union over the finitely many nontrivial elements of $G$, it follows that 
    \[
        \dim\Gamma_{\psi}+1=d(\dim\Gamma_{\varphi}+1).
    \]
    A computation then shows that
    \[
        \dim\mathbf{P}(V^{\oplus d})-\dim\Gamma_{\psi}=d(\dim\mathbf{P}(V)-\dim\Gamma_\varphi).
    \]
    We assume that $\dim\mathbf{P}(V)-\dim\Gamma_\varphi\geqslant 1$, so the right hand side is $\geqslant d$.
\end{proof}

\begin{Corollary}\label{cor:there exist nice faithful vector bundles}
    For any integer $d\geqslant 1$, there exists a faithful locally free coherent sheaf $\ms E$ on $\ms G$ such that the stacky locus $\Gamma_{\ms P}\subset\ms P$ of the projective bundle $\ms P=\mathbf{P}_{\ms G}(\ms E)$ has codimension $\geqslant d$.
\end{Corollary}
\begin{proof}
    By Lemma \ref{lem:there is a faithful locally free sheaf} we may find a faithful locally free coherent sheaf on $\ms G$. By Proposition \ref{prop:tensorial construction} a suitable polynomial tensorial construction gives the desired faithful locally free sheaf.
\end{proof}

Let $\ms E$ be a faithful locally free sheaf on $\ms G$ of rank $n+1$. Let $d$ be a nonnegative integer. We will consider complete intersections of dimension $d$ in the projective bundle $\ms P:=\mathbf{P}_{\ms G}(\ms E)$.

\begin{Definition}
    A \emph{multidegree of dimension $d$} is a sequence $\underline{m}=(m_1,\dots,m_{n-d})$ of positive integers. A \textit{complete intersection in} $\ms P$ \textit{of dimension} $d$ \textit{and multidegree} $\underline{m}=(m_1,\dots,m_{n-d})$ is a closed substack $\ms Z\subset\ms P$ of pure dimension $d$ of the form $V(f_1)\cap\dots\cap V(f_{n-d})$ for some sections $f_i\in\H^0(\ms P,\ms O_{\ms P}(m_i))$.
\end{Definition}

By Lemma \ref{lem:the tautological sheaf descends}, we may choose a positive integer $m$ such that $\ms O_{\ms P}(m)$ descends to a very ample invertible sheaf, say $\ms O_P(1)$, on $P$. Let $d\geqslant 1$ be a positive integer and let $\underline{m}=(m_1,\dots,m_{n-d})$ be a multidegree of dimension $d$ such that each $m_i$ is divisible by $m$. By Lemma \ref{lem:pullback from coarse space}, pullback along $\rho$ induces isomorphisms
\[
    \H^0(P,\ms O_P(m_i/m)\iso\H^0(\ms P,\ms O_{\ms P}(m_i)).
\]
In particular, this shows that the space of complete intersections in $\ms P$ of dimension $d$ and multidegree $\underline{m}$ is positive dimensional.

\begin{Theorem}\label{thm:general complete intersections over a gerbe}
    If the stacky locus $\Gamma_{\ms P}\subset\ms P$ has codimension $\geqslant d+1$, then a general complete intersection $Z\subset\ms P$ of dimension $d$ and multidegree $\underline{m}$ has the following properties.
    \begin{enumerate}
        \item\label{item:item 1 CI} $Z$ is a scheme and $Z$ is smooth and geometrically connected over $K$.
        \item\label{item:item 2 CI} The morphism $Z\to\ms G$ is smooth and geometrically connected.
        \item\label{item:item 3 CI} The inertial action of $G_{\ms G}$ on $Z$ is free.
    \end{enumerate}
\end{Theorem}
\begin{proof}
    We assume that the stacky locus $\Gamma_{\ms P}$ has codimension $\geqslant d+1$ in $\ms P$, so a general complete intersection $Z\subset\ms P$ of dimension $d$ and multidegree $\underline{m}$ does not intersect the stacky locus, and hence is a scheme. This implies that the inertial action of $G_{\ms G}$ on $Z$ is free, so~\eqref{item:item 3 CI} holds. Claim~\eqref{item:item 2 CI} follows from our stacky Bertini's theorem (Theorem \ref{thm:Bertini for stacks}) applied to the morphism $\pi:\ms P\to\ms G$, and the remainder of Claim~\eqref{item:item 1 CI} follows from the stacky Bertini's theorem applied to the morphism $\ms P\to\Spec K$.
\end{proof}

\begin{Corollary}\label{cor:general complete intersections over a gerbe}
    Suppose that $K$ is infinite. For any integer $d\geqslant 1$, there exists a scheme $Z$ of dimension $d$ and a morphism $Z\to\ms G$ which is a complete intersection in a projective bundle over $\ms G$ and which satisfies conditions~\eqref{item:item 1 CI},~\eqref{item:item 2 CI}, and~\eqref{item:item 3 CI} of Theorem \ref{thm:general complete intersections over a gerbe}.
\end{Corollary}
\begin{proof}
    Combine Corollary \ref{cor:there exist nice faithful vector bundles} and Theorem \ref{thm:general complete intersections over a gerbe}.
\end{proof}

\subsection{Complete intersections over a gerbe, over a finite field}

We use Poonen's Bertini-type results over finite fields \cite{MR2144974} to prove that the existence statement of Corollary \ref{cor:general complete intersections over a gerbe} holds also over a finite field.

\begin{Theorem}\label{thm:there exist complete intersections, finite field case}
    Corollary \ref{cor:general complete intersections over a gerbe} holds also when $K$ is finite.
\end{Theorem}
\begin{proof}
    Let $d$ be a positive integer. Let $\ms E$ be a faithful locally free sheaf on $\ms G$ such that the stacky locus $\Gamma$ of the projective bundle $\ms P:=\mathbf{P}_{\ms G}(\ms E)$ has codimension $\geqslant d+1$. Let $\rho:\ms P\to P$ be the coarse space of $\ms P$ and let $m$ be a positive integer so that $\ms O_{\ms P}(m)$ descends to a very ample invertible sheaf say $\ms O_P(1)$ on $P$. Consider the projective embedding $P\subset\mathbf{P}^N$ induced by $\ms O_P(1)$. Applying Poonen's Bertini theorem \cite[Theorem 1.1]{MR2144974} iteratively to a smooth stratification of the subvarieties $P$ and $\rho(\Gamma)$ of $\mathbf{P}^N$, we find a sequence $f_1,\dots,f_{n-d}$ of forms $f_i\in\H^0(\mathbf{P}^N,\ms O_{\mathbf{P}^N}(d_i))$ such that the intersection $X:=P\cap V(f_1)\cap\dots\cap V(f_{n-d})$ is smooth over $K$ of dimension $d$ and does not intersect $\rho(\Gamma)$. Let $Z\subset\ms P$ be the preimage of $X$ under $\rho$. The map $\rho:\ms P\to P$ restricts to an isomorphism on the complement of $\Gamma$, so $\rho$ induces an isomorphism $Z\iso X$. We claim that $Z$ satisfies properties~\eqref{item:item 1 CI},~\eqref{item:item 2 CI}, and~\eqref{item:item 3 CI} of Theorem \ref{thm:general complete intersections over a gerbe}.
    
    For~\eqref{item:item 1 CI}, we note that by construction $Z$ is smooth over $K$. Moreover, $Z$ is isomorphic to $X$, which is a positive dimensional intersection of $P$ with some hypersurfaces in a projective space, and therefore is automatically geometrically connected. To verify the remaining conditions, choose a splitting $s\in\ms G(L)$ of $\ms G$ over some finite extension $L/K$ (in fact, this field extension can be shown to be unnecessary: see Remark \ref{rem:gerbes over finite fields are split}). Following the notation of Example \ref{ex:non split gerbe}, we set $G_L=s^{-1}G_{\ms G}$, $\widetilde{Z}_L=Z\times_{\ms G,s}\Spec L$, and $Z_L=Z\otimes_KL$. As described in Remark \ref{rem:quotient of projective space}, the pullback of $\ms P$ along the splitting $s$ is isomorphic to a projective space $\mathbf{P}_L$ over $L$ which is equipped with a faithful action of the group scheme $G_L$, and the base change of $\ms P$ to $L$ is isomorphic to the quotient stack $[\mathbf{P}_L/ G_L]$. We have a diagram
    \[
        \begin{tikzcd}
            \widetilde{Z}_L\arrow[hook]{d}\arrow{r}{q_Z}&Z_L\arrow[hook]{d}\\
            \mathbf{P}_L\arrow{r}{q}\arrow{d}&\left[\mathbf{P}_L/ G_L\right]\arrow{d}\\
            \Spec L\arrow{r}&\B G_L
        \end{tikzcd}
    \]
    with 2-Cartesian squares. As $Z$ has trivial intersection with $\Gamma$, the $G_L$-action on $\widetilde{Z}_L$ is free, and therefore~\eqref{item:item 3 CI} holds. To verify~\eqref{item:item 2 CI}, we note that the map $q$ restricts to an \'{e}tale morphism over the complement of the stacky locus $\Gamma_L\subset [\mathbf{P}_L/ G_L]$. Thus the map $q_Z:\widetilde{Z}_L\to Z_L$ is \'{e}tale. As $Z$ was smooth over $K$, $Z_L$ is smooth over $L$, so also $\widetilde{Z}_L$ is smooth over $L$, and therefore $Z\to\ms G$ is smooth. Finally, the subvariety $\widetilde{Z}_L\subset\mathbf{P}_L$ is a complete intersection of positive dimension, and hence is geometrically connected. This implies that $Z\to\ms G$ has geometrically connected fibers, and completes the verification of~\eqref{item:item 2 CI}.
\end{proof}

\section{Curves over gerbes}\label{sec:curves over gerbes}

Let $K/K_0$ be a finitely generated field extension. Let $\ms G$ be a Deligne--Mumford gerbe over $K$ with inertia $G_{\ms G}:=\ms I_{\ms G}$. The goal of this section is to prove the following result.

\begin{Theorem}\label{thm:main theorem for curves over gerbes}
    For any integer $N$, there exists a curve $C$ over $K$ of genus $\geqslant N$ and a morphism $C\to\ms G$ which is a relative curve such that
    \begin{enumerate}
        \item\label{item:main curve thm 1} the inertial action induces an isomorphism $G_{\ms G}\iso\sAut_{\ms G}(C)$,
        \item\label{item:main curve thm 2} $\sAut_K(C)=1$, and
        \item\label{item:main curve thm 3} $C$ is not defined over any intermediate field extension $K/L/K_0$ with $K\neq L$.
    \end{enumerate}
\end{Theorem}

Here, by a \emph{curve} over $K$ we mean a smooth proper geometrically integral $K$-scheme of dimension $1$, and by a \emph{relative curve} we mean a smooth proper morphism all of whose geometric fibers are curves. The strategy of the proof is to first select a curve $C$ and a morphism $C\to\ms G$ which is a relative curve, using results of section \S\ref{sec:complete intersections over gerbes}. The inertial action of $G$ on $C$ is faithful, and therefore gives rise to an inclusion $G_{\ms G}\subset\sAut_{\ms G}(C)$. We then refine this initial choice of curve by constructing a curve $C'$ and a finite morphism $C'\to C$ such that $C'\to\ms G$ has the desired properties. The refinement $C'$ is obtained by taking fiber products of certain carefully chosen finite morphisms of curves. We summarize the output of this construction in the following existence result.

\begin{Theorem}\label{thm:finite cover of curve over a gerbe}
    Let $C$ be a curve over $K$ and let $C\to\ms G$ be a morphism which is a relative curve. For any integer $N$, there exists a curve $C'$ over $K$ of genus $\geqslant N$ and a finite separable morphism $C'\to C$ such that $C'\to\ms G$ is a relative curve and $C'$ and $C'\to\ms G$ satisfy conditions~\eqref{item:main curve thm 1},~\eqref{item:main curve thm 2}, and~\eqref{item:main curve thm 3} of Theorem \ref{thm:main theorem for curves over gerbes}.
\end{Theorem}

The existence of the initial curve $C\to\ms G$ is given by Corollary \ref{cor:general complete intersections over a gerbe} (when $K$ is infinite) and Theorem \ref{thm:there exist complete intersections, finite field case} (when $K$ is finite). Thus, from Theorem \ref{thm:finite cover of curve over a gerbe} we deduce Theorem \ref{thm:main theorem for curves over gerbes} as an immediate consequence.

The proof of Theorem \ref{thm:finite cover of curve over a gerbe} occupies the remainder of this section, whose organization we outline. In section \S\ref{ssec:pencils on curves over gerbes} we show that every curve over a gerbe admits pencils with prescribed automorphism groups. Sections \S \ref{ssec:results on curves 1}, \S\ref{ssec:results on curves 2}, and \S \ref{ssec:incompressible morphisms} consist of some preparatory results concerning the descent of isomorphisms of curves along finite morphisms of curves with certain properties. In \S \ref{ssec:curves over gerbes} we prove some results on automorphism groups of curves over gerbes. In \S\ref{ssec:curves over gerbes over an infinite field}, we restrict to the case when $K$ is infinite, and prove Theorem \ref{thm:finite cover of curve over a gerbe, infinite field case}, which is the partial result that there exists a refinement $C'\to C$ of arbitrarily large genus satisfying conditions~\eqref{item:main curve thm 1} and~\eqref{item:main curve thm 2} of Theorem \ref{thm:main theorem for curves over gerbes}. We then consider fields of definition. In \S\ref{ssec:descent} we recall some definitions related to descent data, and show that a descent datum in curves with respect to an arbitrary field extension is effective. In \S\ref{ssec:fields of definition} we give some results which will allow us to control the fields of definition of curves related by certain finite morphisms. In \S\ref{ssec:curves over gerbes over an infinite field + field of definition} we complete the proof of Theorem \ref{thm:finite cover of curve over a gerbe} in the case when $K$ is infinite (Theorem \ref{thm:finite cover of curve over a gerbe, field of definition, infinite field case}). Finally, in \S\ref{ssec:curves over gerbes with prescribed field of definition, finite field case}, we consider the case when $K$ is finite. After some preparations, we prove the finite field case of Theorem \ref{thm:finite cover of curve over a gerbe} (Theorem \ref{thm:finite cover of curve over a gerbe, field of definition, finite field case}). We refer to \S\ref{ssec:curves over gerbes with prescribed field of definition, finite field case} for a discussion of the differences between the infinite and finite field cases.

We will use the following terminology with respect to curves and their morphisms.
\begin{Notation}
    A \emph{curve} over a field $K$ is a smooth proper geometrically integral $K$-scheme of dimension $1$. A \emph{relative curve} or a \emph{family of curves} over an algebraic stack $\ms S$ is a smooth proper morphism $\ms C\to\ms S$ all of whose geometric fibers are curves. Let $f:C\to D$ be a finite separable morphism of curves over $K$. The \textit{ramification locus} of $f$ is the set of closed points $P\in C$ such that $f$ is not smooth at $P$. The \textit{branch locus} of $f$ is the image in $D$ of the ramification locus. We say that $f$ is \emph{separable} (resp. \emph{Galois}) if the corresponding field extension $k(C)/k(D)$ is separable (resp. Galois).
A \textit{Galois closure} of $f$ is a finite morphism $\widetilde{C}\to C$ of curves over $K$ such that the composition $\widetilde{C}\to C\to D$ is Galois. A \textit{minimal Galois closure} of $f$ is a Galois closure of minimal degree.
\end{Notation}

\subsection{Pencils on curves over gerbes}\label{ssec:pencils on curves over gerbes}

In this section we consider pencils on curves over gerbes. We begin with some results concerning pencils on curves in a projective space. Let $C$ be a curve over $K$. Let $\ms L$ be a very ample invertible sheaf on $C$.
Let $\mathbf{P}:=\mathbf{P}(\H^0(C,\ms L))$ be the projective space classifying codimension one subspaces of $\H^0(C,\ms L)$. The canonical surjection $\H^0(C,\ms L)\otimes\ms O_C\twoheadrightarrow\ms L$ gives rise to a projective embedding $C\hookrightarrow\mathbf{P}$ which is \emph{nondegenerate} (that is, does not factor through a proper linear subvariety of $\mathbf{P}$). We write $|\ms L|$ for the complete linear system consisting of those divisors in $C$ which are the vanishing loci of global sections of $\ms L$. The complete linear system $|\ms L|$ is isomorphic to the dual projective space $\mathbf{P}^{\vee}=\mathbf{P}(\H^0(C,\ms L)^{\vee})$ classifying hyperplanes in $\mathbf{P}$.

\begin{Proposition}\label{prop:disjoint divisors}
    Fix a nontrivial automorphism $\alpha\in\Aut_K(C)$. If $Z$ is a general divisor in $|\ms L|$, then the divisors $Z$ and $\alpha\cdot Z$ are disjoint.
\end{Proposition}
\begin{proof}
    Consider the \emph{incidence correspondence}
    \[
        \begin{tikzcd}
            \Sigma\subset C\times\mathbf{P}^{\vee}
        \end{tikzcd}
    \]
    where $\Sigma$ is the scheme classifying pairs $(P,H)$, where $H$ is a hyperplane in $\mathbf{P}$ and $P$ is a point in $C$ which is contained in $H$. Define $\Sigma^{\alpha}\subset C\times\mathbf{P}^{\vee}$ to be the pullback of $\Sigma$ under the automorphism $\alpha\times\id:C\times\mathbf{P}^{\vee}\iso C\times\mathbf{P}^{\vee}$. Thus $\Sigma^{\alpha}$ consists of pairs $(P,H)$ such that $\alpha(P)\in H$. We claim that the intersection $\Sigma\cap\Sigma^{\alpha}$ does not dominate $\mathbf{P}^{\vee}$. To see this, we note that the morphism $\Sigma\to C$ is a projective bundle, so $\Sigma$ and $\Sigma^{\alpha}$ are irreducible. The morphisms $\Sigma\to\mathbf{P}^{\vee}$ and $\Sigma^{\alpha}\to\mathbf{P}^{\vee}$ are finite, so if the intersection $\Sigma\cap \Sigma^{\alpha}$ were to dominate $\mathbf{P}^{\vee}$, then we would have $\Sigma=\Sigma^{\alpha}$ as subvarieties of $C\times\mathbf{P}^{\vee}$. This would imply that for every hyperplane $H$ the intersection $C\cap H$ is preserved by $\alpha$. But if $P$ is any point of $C$, we can find two hyperplanes whose common intersection with $C$ consists of exactly $P$. This implies that every point of $C$ is fixed by $\alpha$, so $\alpha$ is the identity, contrary to our assumption.
\end{proof}

\begin{Definition}
    A \emph{pencil} in the complete linear system $|\ms L|$ is a one dimensional linear subvariety $\mathbf{L}\subset|\ms L|$. Equivalently, a pencil in $|\ms L|$ may be determined by giving a one dimensional linear subvariety of $\mathbf{P}^{\vee}$, a one-dimensional linear family of hyperplanes in $\mathbf{P}$, or a two dimensional subspace of the vector space $\H^0(C,\ms L)$. 
\end{Definition}

As the embedding $C\hookrightarrow\mathbf{P}$ is nondegenerate, any pencil $\mathbf{L}\subset|\ms L|$ gives rise to a nonconstant rational map $C\dashrightarrow\mathbf{L}$ defined away from the base locus of $\mathbf{L}$. The fiber of this map over a point $x\in\mathbf{L}\subset|\ms L|=\mathbf{P}^{\vee}$ which is not in the base locus is the divisor $C\cap H_x$, where $H_x$ is the hyperplane in $\mathbf{P}$ classified by $x$. If $\mathbf{L}$ is basepoint free (a generic condition), then we obtain a finite morphism $C\to\mathbf{L}$. We say that $\mathbf{L}$ is \emph{separable} if this morphism 
is separable. By Bertini's theorem, a general hyperplane in $\mathbf{P}$ intersects $C$ transversely, so a generic pencil is separable.
\begin{Notation}
    Given a basepoint free pencil $\mathbf{L}\subset |\ms L|$ and a choice of an isomorphism $\mathbf{L}\cong\mathbf{P}^1$, we will write 
    \[
        f_{\mathbf{L}}:C\to\mathbf{P}^1
    \]
    for the resulting finite morphism to the projective line.
\end{Notation}

\begin{Proposition}\label{prop:general pencil has no aut}
    Assume that $C$ has genus $g_C\geqslant 2$. If $\mathbf{L}$ is a general pencil in the complete linear system $|\ms L|$ and $\mathbf{L}\cong\mathbf{P}^1$ is an isomorphism, then the resulting morphism $f_{\mathbf{L}}:C\to\mathbf{P}^1$ is separable and satisfies $\Aut_K(C/\mathbf{P}^1)=1$.
\end{Proposition}
\begin{proof}
    If $\alpha$ is an automorphism of $C$, then by Proposition \ref{prop:disjoint divisors} there is a nonempty open subset of $\mathbf{P}^{\vee}$ corresponding to divisors of $C$ which are not preserved by $\alpha$. The union of these open subsets over the finitely many elements of $\Aut_K(C)$ is again nonempty. Thus, a general pencil $\mathbf{L}\subset\mathbf{P}^{\vee}$ is basepoint free, separable, and intersects this open subset, and so gives rise to a finite morphism $f_{\mathbf{L}}:C\to\mathbf{P}^1$ with the property that $\Aut_K(C/\mathbf{P}^1)=1$.
\end{proof}

We now consider the following situation. Let $G$ be a finite group. Let $\widetilde{C}$ be a curve over $K$ equipped with a faithful $G$-action, set $C=\widetilde{C}/G$, and let
\[
    q:\widetilde{C}\to C
\]
be the quotient map. As before, we let $\ms L$ be a very ample invertible sheaf on $C$, write $\mathbf{P}=\mathbf{P}(\H^0(C,\ms L))$, and let $C\hookrightarrow\mathbf{P}$ be the resulting projective embedding. By replacing $\ms L$ with a suitable tensor power, we may assume that there exists a $G$-equivariant invertible sheaf $\widetilde{\ms L}$ on $\widetilde{C}$ which is very ample such that the pullback $q^*\ms L$ is isomorphic to $\widetilde{\ms L}^{\otimes |G|}$.

\begin{Proposition}\label{prop:Gamma=G}
    Fix an automorphism $\alpha\in\Aut_K(\widetilde{C})$. If $Z\subset C$ is a general divisor in the linear system $|\ms L|$, then $\alpha$ preserves $q^{-1}(Z)$ setwise if and only if $\alpha\in G$.
\end{Proposition}
\begin{proof}
    Let $W\subset\widetilde{C}$ be a general divisor in the linear system $|\widetilde{\ms L}|$ and set
    \[
        Z:=q\left(\sum_{g\in G}g\cdot W\right).
    \]
    Then $Z\subset C$ is a divisor in the linear system $|\ms L|$. We claim that $Z$ has the desired property, that is, that an automorphism $\alpha\in\Aut_K(\widetilde{C})$ preserves $q^{-1}(Z)$ setwise if and only if $\alpha\in G$. For this, we note that we have
    \[
        q^{-1}(Z)=\sum_{g\in G}g\cdot W.
    \]
    Thus, if $\alpha\in G$, then $\alpha$ certainly preserves $q^{-1}(Z)$ setwise. Conversely, if an automorphism $\alpha\in\Aut_K(\widetilde{C})$ preserves $q^{-1}(Z)$, then $\alpha\cdot W$ intersects $g\cdot W$ nontrivially for some $g\in G$. This implies that $(g^{-1}\alpha)\cdot W$ intersects $W$ nontrivially. But $W$ was chosen to be general, so by Proposition \ref{prop:disjoint divisors} we have that $g^{-1}\alpha$ is the identity. We conclude that $\alpha=g^{-1}$ and hence $\alpha\in G$.
\end{proof}

We now consider finite morphisms $\widetilde{C}\to\mathbf{P}^1$ induced by pencils of divisors on $C$. Given a basepoint free pencil $\mathbf{L}\subset |\ms L|$ and an isomorphism $\mathbf{L}\cong\mathbf{P}^1$, we write $\widetilde{f}_{\mathbf{L}}=f_{\mathbf{L}}\circ q:\widetilde{C}\to\mathbf{P}^1$. By construction, the morphism $\widetilde{f}_{\mathbf{L}}$ is $G$-equivariant, so we have an inclusion $G\subset\Aut_K(\widetilde{C}/\mathbf{P}^1)$.

\begin{Theorem}\label{thm:generic G equivariant pencils on a curve}
    Assume that $C$ has genus $g_{C}\geqslant 2$. If $\mathbf{L}$ is a general pencil of hyperplanes in $\mathbf{P}$ and $\mathbf{L}\cong\mathbf{P}^1$ is an isomorphism, then the resulting morphisms $f_{\mathbf{L}}:C\to\mathbf{P}^1$ and $\widetilde{f}_{\mathbf{L}}:\widetilde{C}\to\mathbf{P}^1$ are separable, and we have
    \begin{enumerate}
        \item $\Aut_K(\widetilde{C}/\mathbf{P}^1)=G$ and
        \item $\Aut_K(C/\mathbf{P}^1)=1$.
    \end{enumerate}
\end{Theorem}
\begin{proof}
    By Proposition \ref{prop:general pencil has no aut}, if $\mathbf{L}$ is a general pencil in $|\ms L|$ and $\mathbf{L}\cong\mathbf{P}^1$ is an isomorphism then the resulting morphism $f_{\mathbf{L}}:C\to\mathbf{P}^1$ is separable and satisfies $\Aut_K(C/\mathbf{P}^1)=1$. The map $\widetilde{C}\to C$ is generically \'{e}tale, so the morphism $\widetilde{f}_{\mathbf{L}}:\widetilde{C}\to\mathbf{P}^1$ is also separable. It remains to show that $\Aut_K(\widetilde{C}/\mathbf{P}^1)=G$. Fix an automorphism $\alpha\in\Aut_K(\widetilde{C})$. By Proposition \ref{prop:Gamma=G}, there is a nonempty open subset of $\mathbf{P}^{\vee}$ corresponding to divisors $Z\subset C$ in the linear system $|\ms L|$ which have the property that $\alpha$ preserves $q^{-1}(Z)$ setwise if and only if $\alpha\in G$. Thus, if $\mathbf{L}$ is a general pencil of hyperplanes in $\mathbf{P}$ which intersects this open subset, then we have $\alpha\in\Aut_K(\widetilde{C}/\mathbf{P}^1)$ if and only if $\alpha\in G$. The union of these open subsets taken over the finitely many elements of $\Aut_K(\widetilde{C})$ is again a nonempty open subset, and any basepoint free pencil $\mathbf{L}\subset\mathbf{P}^{\vee}$ intersecting this subset gives rise to a finite morphism $\widetilde{f}_{\mathbf{L}}:\widetilde{C}\to\mathbf{P}^1$ with the property that $\Aut_K(\widetilde{C}/\mathbf{P}^1)=G$.
\end{proof}

We now work over a gerbe. Let $\ms G$ be a Deligne--Mumford gerbe over $K$ with inertia $G_{\ms G}:=\ms I_{\ms G}$. Let $C$ be a curve over $K$ and let $\pi:C\to\ms G$ be a morphism which is a relative curve. As before, we let $\ms L$ be a very ample invertible sheaf on $C$, write $\mathbf{P}=\mathbf{P}(\H^0(C,\ms L))$, and let $C\hookrightarrow\mathbf{P}$ be the resulting projective embedding.

\begin{Theorem}\label{thm:there exists a curve over a gerbe with a pencil}
    Assume that $C$ has genus $g_C\geqslant 2$. After possibly replacing $\ms L$ by a positive multiple, we have that if $\mathbf{L}$ is a general pencil of hyperplanes in $\mathbf{P}$ and $\mathbf{L}\cong\mathbf{P}^1$ is an isomorphism, then the resulting morphism $f_{\mathbf{L}}:C\to\mathbf{P}^1$ is finite and separable, and furthermore
    \begin{enumerate}
        \item\label{item:CI curve 1} the inertial action induces an isomorphism $G_{\ms G}\iso\sAut_{\ms G}(C/\mathbf{P}^1)$, and
        \item\label{item:CI curve 2} $\sAut_K(C/\mathbf{P}^1)=1$.
    \end{enumerate}
\end{Theorem}
\begin{proof}
    Let $L$ be a separable closure of $K$ and choose a section $s\in\ms G(L)$. The pullback of $G_{\ms G}$ along $s$ is then the group scheme associated to a finite group, say $G$. Set $\widetilde{C}_L=C\times_{\ms G,s}\Spec L$. We obtain a diagram
    \[
        \begin{tikzcd}
            \widetilde{C}_L\arrow{d}\arrow{r}{q_L}&C_L\arrow{d}\arrow[hook]{r}&\mathbf{P}_L\\
            \Spec L\arrow{r}&\B_L\!G&
        \end{tikzcd}
    \]
    where the square is 2-Cartesian, the morphism $q_L$ is a $G$-torsor, and the morphism $C_L\hookrightarrow\mathbf{P}_L$ is the base change of the projective embedding $C\hookrightarrow\mathbf{P}$ to $L$. Let $\ms L_L$ be the base change of $\ms L$ to $C_L$. After replacing $\ms L$ by a suitable tensor power, we may arrange so that we are in the situation of Theorem \ref{thm:generic G equivariant pencils on a curve}; namely, so that there exists a $G$-equivariant very ample invertible sheaf $\widetilde{\ms L}_L$ on $\widetilde{C}_L$ such that $q_L^*\ms L_L$ is isomorphic to $\widetilde{\ms L}_L^{\otimes |G|}$. Consider a morphism $f:C\to\mathbf{P}^1$ induced by a general pencil of hyperplanes in $\mathbf{P}$. Let $f_L:C_L\to\mathbf{P}^1_L$ be the base change of $f$ to $L$ and set $\widetilde{f}_L=f_L\circ q_L:\widetilde{C}_L\to\mathbf{P}^1_L$. The morphism $f_L$ is induced by a general pencil of hyperplanes in $\mathbf{P}_L$, so by Theorem \ref{thm:generic G equivariant pencils on a curve} we have that $f_L$ and $\widetilde{f}_L$ are separable, hence $f$ is separable, and furthermore we have $\Aut_L(C_L/\mathbf{P}^1_L)=1$ and $G=\Aut_L(\widetilde{C}_L/\mathbf{P}^1_L)$. The group $\Aut_L(C_L/\mathbf{P}^1_L)$ is the group of $L$-points of $\sAut_K(C/\mathbf{P}^1)$, so $\sAut_K(C/\mathbf{P}^1)=1$. The natural map $G\to\Aut_L(\widetilde{C}_L/\mathbf{P}^1_L)$ is the pullback along $s$ of the inertial action map $G_{\ms G}\to\sAut_{\ms G}(C/\mathbf{P}^1)$, which is therefore an isomorphism.
\end{proof}

\begin{Corollary}\label{cor:there exists a curve over a gerbe with a pencil}
    If $K$ is infinite, then there exists a finite separable morphism $f:C\to\mathbf{P}^1$ satisfying conditions~\eqref{item:CI curve 1} and~\eqref{item:CI curve 2} of Theorem \ref{thm:there exists a curve over a gerbe with a pencil}.
\end{Corollary}

\begin{Remark}
    We will show later that Corollary \ref{cor:there exists a curve over a gerbe with a pencil} also holds when $K$ is finite (combine Proposition \ref{prop:pencils with auts over finite fields} and Lemma \ref{lem:pencils over finite fields}).
\end{Remark}

As a further input in our constructions, we will need a curve over $K$ with trivial automorphism group. The existence of such a curve when $K$ is infinite, or finite of sufficiently large cardinality, has been shown by Katz--Sarnak \cite[Theorem 10.6.14, Remark 10.6.24]{MR1659828}, over an arbitrary finite field by Rzedowski--Calder\'{o}n and Villa--Salvador \cite[Theorem 2]{MR1120718}, and in general by Poonen \cite{MR1748288}. We will use Poonen's construction, which has the useful additional property of exhibiting curves which admit a degree 3 morphism to $\mathbf{P}^1$ which is totally ramified over $\infty$. The specific existence result we will need is the following.


\begin{Theorem}[{Poonen \cite[Theorem 1]{MR1748288}}]\label{thm:curves with no auts}
    If $K$ is any field, then there exists a curve $X$ over $K$ such that $\sAut_K(X)=1$. Furthermore, we may choose such an $X$ so that there exists a finite separable morphism $X\to\mathbf{P}^1$ which is totally ramified over $\infty$.
\end{Theorem}

\subsection{Finite morphisms of curves 1: condition $(\ast)$}\label{ssec:results on curves 1}

Let $K$ be a field. Inspired by Madden--Valentini \cite{MR705883}, we consider the following condition on a finite morphism of curves.

\begin{Definition}\label{def:condition star}
    We say that a finite morphism $f:C\to D$ of curves over $K$ \emph{satisfies condition} $(\ast)$ if given any curve $V$ over $K$ and a factorization $C\to V\to D$ of $f$ such that $V\to D$ is not an isomorphism, we have
    \[
        g_V\geqslant\deg(V/D)^2+2(g_D-1)\deg(V/D)+2
    \]
    where $g_D$ is the genus of $D$ and $g_V$ is the genus of $V$.
\end{Definition}

We will show that isomorphisms of curves descend along morphisms which satisfy condition $(\ast)$. We will use the following result, known as the inequality of Castelnuovo--Severi. A reference for this result is Stichtenoth \cite[Satz 1]{MR733931}. For another proof and a historical discussion, we refer the reader to Kani \cite{MR758693} and the associated references.

\begin{Theorem}
    Let $f_1:C\to D_1$ and $f_2:C\to D_2$ be finite morphisms of curves over $K$ of degrees $d_1$ and $d_2$. If the induced morphism $f_1\times f_2:C\to D_1\times D_2$ is birational onto its image, then
    \[
        g_C\leqslant (d_1-1)(d_2-1)+d_1g_{D_1}+d_2g_{D_2}.
    \]
\end{Theorem}

The following consequence of the Castelnuovo--Severi inequality is due to Madden--Valentini \cite[Lemma 1]{MR705883}. For the sake of completeness, we give the proof.

\begin{Proposition}\label{prop:CSI consequence}
    Let $f_1:C_1\to D_1$ and $f_2:C_2\to D_2$ be finite morphisms of curves over $K$. Suppose that $f_1$ and $f_2$ have the same degree, that $D_1$ and $D_2$ have the same genus, and that $f_1$ satisfies condition $(\ast)$.
    If $\alpha:C_1\iso C_2$ is an isomorphism, then there exists a unique isomorphism $\beta:D_1\iso D_2$ such that the diagram
    \[
        \begin{tikzcd}
            C_1\arrow{r}{\alpha}[swap]{\sim}\arrow{d}[swap]{f_1}&C_2\arrow{d}{f_2}\\
            D_1\arrow{r}{\beta}[swap]{\sim}&D_2
        \end{tikzcd}
    \]
    commutes.
\end{Proposition}
\begin{proof}
    Let $\alpha:C_1\iso C_2$ be an isomorphism. Let $V$ be the normalization of the image of the morphism 
    \[
        f_1\times (f_2\circ\alpha):C_1\to D_1\times D_2.
    \]
    Thus $V$ is a regular curve over $K$. The morphism $C_1\to V$ is flat, and $C_1$ is smooth over $K$, so in fact $V$ is smooth over $K$ as well (we remark that the distinction between a regular and smooth curve is only relevant when $K$ is imperfect). For $i=1,2$ let $p_i:V\to D_i$ be the composition of the normalization map and the projection $\pi_i:D_1\times D_2\to D_i$. We have a diagram
    \[
        \begin{tikzcd}
            C_1\arrow{r}\arrow{dr}&V\arrow{r}\arrow{d}{p_i}&D_1\times D_2\arrow{dl}{\pi_i}\\
            &D_i.&
        \end{tikzcd}
    \]
    We apply the inequality of Castelnuovo--Severi to the two morphisms $p_1,p_2$. Let $d$ be the degree of $p_1$. We assume that the degrees of $f_1$ and $f_2$ are equal, so the degree of $p_2$ is also $d$. Let $g$ be the genus of $D_1$, which by assumption is equal to the genus of $D_2$. We obtain
    \[
        g_V\leqslant (d-1)^2+2dg=d^2+2(g-1)+1.
    \]
    This contradicts the inequality in condition $(\ast)$ applied to the factorization $C_1\to V\xrightarrow{p_1}D_1$ of $f_1$ unless $p_1$ is an isomorphism. This must therefore be the case. It follows that $V$ is isomorphic to its image in $D_1\times D_2$, and this image is the graph of an isomorphism $\beta:D_1\iso D_2$. By construction, $\beta\circ f_1=f_2\circ\alpha$. The uniqueness of $\beta$ follows from the fact that $f_1$ is faithfully flat.
\end{proof}

The following extends Proposition \ref{prop:CSI consequence} to families of curves.

\begin{Proposition}\label{prop:CSI consequence in families}
    Let $S$ be a $K$-scheme. Let $C_1,C_2,D_1,D_2$ be relative curves over $S$. Let $f_1:C_1\to D_1$ and $f_2:C_2\to D_2$ be finite separable morphisms over $S$. Suppose that for any field $L$ and morphism $s:\Spec L\to S$ the pullbacks of $C_i,D_i,$ and $f_i$ along $s$ satisfy the assumptions of Proposition \ref{prop:CSI consequence}. If $\alpha:C_1\iso C_2$ is an isomorphism over $S$, then there exists a unique isomorphism $\beta:D_1\iso D_2$ over $S$ such that the diagram
    \[
        \begin{tikzcd}
            C_1\arrow{r}{\alpha}[swap]{\sim}\arrow{d}[swap]{f_1}&C_2\arrow{d}{f_2}\\
            D_1\arrow{r}{\beta}[swap]{\sim}&D_2
        \end{tikzcd}
    \]
    commutes.
\end{Proposition}
\begin{proof}
    The proof is a reduction to Proposition \ref{prop:CSI consequence}. We consider the morphisms
    \begin{equation}\label{eq:b factors through a}
        \begin{tikzcd}
            &\sIsom_S(D_1,D_2)\arrow{d}{a}\\
            \sIsom_S(C_1,C_2)\arrow{r}{b}\arrow[dashed]{ur}&\sHom_S(C_1,D_2)
        \end{tikzcd}
    \end{equation}
    of Hom and Isom schemes over $S$, where $a=\_\circ f_1$ is precomposition with $f_1$ and $b=f_2\circ\_$ is postcomposition with $f_2$. We claim that $b$ factors through $a$, that is, that there exists a dashed arrow filling in the diagram. We will first show that $a$ is an immersion. Consider the group $\sAut_S(D_2)$, which acts freely and transitively on $\sIsom_S(D_1,D_2)$. By Proposition \ref{prop:CSI consequence} the morphism $\sIsom_S(D_1,D_2)\to S$ is surjective, and in particular the scheme $\sIsom_S(D_1,D_2)$ is nonempty. Thus, the quotient $\sIsom_S(D_1,D_2)/\sAut_S(D_2)$ is isomorphic to $S$. The group $\sAut_S(D_2)$ also acts freely on $\sHom_S(C_1,D_2)$, and the morphism $a$ is compatible with these actions. We therefore have a Cartesian diagram
    \[
        \begin{tikzcd}
            \sIsom_S(D_1,D_2)\arrow{r}\arrow{d}{a}&S\arrow{d}{\overline{a}}\\
            \sHom_S(C_1,D_2)\arrow{r}&\sHom_S(C_1,D_2)/\sAut_S(D_2).
        \end{tikzcd}
    \]
    As $C_1$ is projective, the Hom-scheme $\sHom_S(C_1,D_2)$ is quasi-projective \cite[Theorem 5.23]{MR2223407}. The quotient $\sHom_S(C_1,D_2)/\sAut_S(D_2)$ is therefore a scheme \cite[Theorem 5.25]{MR2223407}, \cite{MR0555258}. The morphism $\overline{a}$ is thus a section of a morphism of schemes, and hence is an immersion \cite[01KT]{stacks-project}. It follows that $a$ is an immersion, as claimed.
    
    We now show that $b$ factors through $a$. It will suffice to prove this in the universal case, when $S$ is (a smooth cover of) the Hurwitz stack parametrizing tuples $(C_1,C_2,D_1,D_2,f_1,f_2)$ where $C_1,C_2,D_1,D_2$ are curves and $f_1:C_1\to D_1$ and $f_2:C_2\to D_2$ are finite separable morphisms. As this stack is smooth over $K$, we may therefore assume that $S$ is smooth over $K$, and in particular that $S$ is reduced. As the Hom and Isom schemes involved are smooth over $S$, they are also then smooth over $K$, and hence also are reduced. We have shown that $a$ is an immersion, so to show that $b$ factors through $a$ it will suffice to show that this is true at the level of geometric points, that is, that the image of every geometric point of $\sIsom_S(C_1,C_2)$ in $\sHom_S(C_1,D_2)$ is the image of a geometric point of $\sIsom_S(D_1,D_2)$. We may therefore assume that $S$ is the spectrum of a field, and check that we have the desired factorization at the level of Hom and Isom sets. This is exactly the statement of Proposition \ref{prop:CSI consequence}.
    
    We have shown that $b$ factors through $a$, so we have a morphism
    \[
        \sIsom_S(C_1,C_2)\to\sIsom_S(D_1,D_2)
    \]
    of Isom schemes over $S$ completing the diagram~\eqref{eq:b factors through a}. The image of the global section 
    \[
        \alpha\in\sIsom_S(C_1,C_2)(S)=\Isom_S(C_1,C_2)
    \]
    under this map gives the desired isomorphism $\beta$. The uniqueness of $\beta$ follows from the fact that $f_1$ is faithfully flat.
\end{proof}

\subsection{Finite morphisms of curves 2: condition $(\ast\ast)$}\label{ssec:results on curves 2}

In this section we will consider another condition on a finite morphism of curves which will again imply that certain isomorphisms of curves descend. Consider a commuting square
\begin{equation}\label{eq:balanced square}
    \begin{tikzcd}
        F\arrow{d}[swap]{g'}\arrow{r}{f'}&E\arrow{d}{g}\\
        C\arrow{r}{f}&D
    \end{tikzcd}
\end{equation}
where $C,D,E,$ and $F$ are curves over $K$ and the arrows are finite separable $K$--morphisms.
\begin{Definition}\label{def:balanced}
    We say that the square~\eqref{eq:balanced square} \textit{satisfies condition} $(\ast\ast')$ if (1) it is Cartesian, and (2) for every curve $\widetilde{C}$ over $K$ and every Galois closure $\widetilde{C}\to C$ of $C\to D$ such that $\deg(\widetilde{C}/D)\leqslant (\deg(C/D)^2)!$, the fiber product $\widetilde{C}\times_D E$ is connected.

    We say that~\eqref{eq:balanced square} \textit{satisfies condition} $(\ast\ast)$ if its base change to a separable closure of $K$ satisfies condition $(\ast\ast')$.
\end{Definition}

We will produce squares satisfying the above conditions by taking suitable fiber products of curves. Let $C,D,$ and $E$ be curves over $K$. Let $f:C\to D$ and $g:E\to D$ be finite separable morphisms of curves, and set $F=C\times_DE$.

\begin{Lemma}\label{lem:smoothness of fiber products}
    The fiber product $F=C\times_DE$ is smooth if and only if the branch loci of $f$ and $g$ are disjoint (as subschemes of $D$).
\end{Lemma}
\begin{proof}
    We may assume that $K$ is algebraically closed. Consider closed points $P\in C$ and $Q\in E$ whose images in $D$ are equal, and write $R=f(P)=g(Q)$. The scheme $F$ is smooth at the point $(P,Q)$ if and only if the tangent space $T_{F,(P,Q)}$ of $F$ at $(P,Q)$ is one-dimensional. It follows from the functorial description of the tangent space as pointed maps from the scheme $\Spec K[\varepsilon]$ that this tangent space is given by the fiber product
    \[
        T_{F,(P,Q)}=T_{C,P}\times_{T_{D,R}}T_{E,Q},
    \]
    and therefore may be computed by the exact sequence
    \[
        0\to T_{F,(P,Q)}\to T_{C,P}\oplus T_{E,Q}\xrightarrow{\mathrm{d}f-\mathrm{d}g}T_{D,R}.
    \]
    As $C,D,$ and $E$ are smooth, $T_{C,P},T_{D,R},$ and $T_{E,Q}$ are all one-dimensional $K$-vector spaces. By the rank nullity theorem, $T_{F,(P,Q)}$ is one dimensional if and only if at least one of $\mathrm{d}f$ or $\mathrm{d}g$ is surjective. We conclude that $F$ is smooth at $(P,Q)$ if and only if either $f$ is unramified at $P$ or $g$ is unramified at $Q$, or both.
\end{proof}

\begin{Lemma}\label{lem:compositum omnibus lemma}
    Suppose that the branch loci of $f$ and $g$ are disjoint (as subschemes of $D$) and that at least one of the following two conditions holds.
    \begin{enumerate}
        \item\label{item:fiber product condition 1} $g$ is totally ramified over some closed point of $D$.
        \item\label{item:fiber product condition 2} $\deg(C/D)<\deg(E/D)$ and $\deg(E/D)$ is prime.
    \end{enumerate}
     The fiber product $F=C\times_DE$ is a smooth proper geometrically integral curve, and the commuting square~\eqref{eq:balanced square} satisfies condition $(\ast\ast)$.
\end{Lemma}
\begin{proof}
    By Lemma \ref{lem:smoothness of fiber products} the fiber product $F=C\times_DE$ is smooth. All of the above conditions are stable under field extension, so it will suffice to consider the case when $K$ is separably closed. Let $\widetilde{C}\to C$ be a Galois closure of $C\to D$ such that $\deg(\widetilde{C}/D)\leqslant (\deg(C/D)^2)!$. Set $\widetilde{F}=\widetilde{C}\times_DE$, so that we have a diagram
    \[
        \begin{tikzcd}
            \widetilde{F}\arrow{d}\arrow{r}&F\arrow{d}[swap]{g'}\arrow{r}{f'}&E\arrow{d}{g}\\
            \widetilde{C}\arrow{r}&C\arrow{r}{f}&D
        \end{tikzcd}
    \]
    with Cartesian squares. As $F$ is smooth, to show that it is integral we only need to verify that it is connected. Thus, to complete the proof, it will suffice to verify that under either of the given conditions $\widetilde{F}$ is connected. Suppose that~\eqref{item:fiber product condition 1} holds, and so $g$ is totally ramified over some closed point say $P\in D$. Let $Q$ be a closed point of $\widetilde{C}$ mapping to $P$. Then $\widetilde{F}\to\widetilde{C}$ is also totally ramified over $Q$, and hence its fiber over $Q$ consists of a single closed point. Every irreducible component of $\widetilde{F}$ is proper over $\widetilde{C}$, and hence contains a point mapping to $Q$. Therefore the irreducible components of $\widetilde{F}$ all intersect at a common point, and so $\widetilde{F}$ is connected. Suppose that~\eqref{item:fiber product condition 2} holds. It follows from our assumptions that the morphisms $\widetilde{C}\to D$ and $E\to D$ have coprime degrees, which implies that $\widetilde{F}$ is connected.
\end{proof}

We record the following observation.

\begin{Lemma}\label{lem:Galois base change}
    Consider a commuting square of the form~\eqref{eq:balanced square}, where $C,D,E,$ and $F$ are curves over $K$ and $f,f',g,$ and $g'$ are finite separable $K$-morphisms. Suppose that this square is Cartesian. If $f$ is Galois, then $f'$ is also Galois, and pullback along $g$ induces an isomorphism
    \begin{equation}\label{eq:pullback map galois}
        \Aut_K(C/D)\iso\Aut_K(F/E).
    \end{equation}
\end{Lemma}
\begin{proof}
    As $g$ is faithfully flat, the pullback map~\eqref{eq:pullback map galois} is injective. Using our assumption that $f$ is Galois, we have the inequalities
    \[
        \deg(C/D)=|\Aut_K(C/D)|\leqslant |\Aut_K(F/E)|\leqslant \deg(F/E).
    \]
    But the square is Cartesian, so $\deg(C/D)=\deg(F/E)$. It follows that $|\Aut_K(C/D)|=|\Aut_K(F/E)|$, so the pullback map~\eqref{eq:pullback map galois} is bijective, and that $|\Aut_K(F/E)|=\deg(F/E)$, so $f'$ is Galois.    
\end{proof}

We now give some consequences of condition $(\ast\ast)$.
\begin{Lemma}\label{lem:there is a curve (Galois)}
    Consider a commuting square of the form~\eqref{eq:balanced square}, where $C,D,E,$ and $F$ are curves over $K$ and $f,f',g,$ and $g'$ are finite separable $K$-morphisms. Assume that this square satisfies condition $(\ast\ast)$. Suppose given a curve $V$ over $K$ and a factorization $F\to V\to E$ of $f'$. There exists a curve $U$ over $K$ and a commutative diagram
    \[
        \begin{tikzcd}
            F\arrow{r}\arrow[bend left=25]{rr}{f'}\arrow{d}[swap]{g'}&V\arrow{d}\arrow{r}&E\arrow{d}{g}\\
            C\arrow{r}\arrow[bend right=25]{rr}[swap]{f}&U\arrow{r}&D
        \end{tikzcd}
    \]
    in which both squares are Cartesian. Moreover, the curve $U$ and the above diagram are uniquely determined up to a unique isomorphism.
\end{Lemma}
\begin{proof}
    By Galois descent, it will suffice to consider the case when $K$ is separably closed. Let $\widetilde{C}\to C$ be a minimal Galois closure of $C\to D$. Thus, we have that $\deg(\widetilde{C}/D)\leqslant\deg(C/D)!$. We claim that the fiber product $\widetilde{F}=\widetilde{C}\times_DE$ is again a smooth proper geometrically integral curve. It will suffice to show that $\widetilde{F}$ is both smooth and connected. The latter follows from our assumption that the square satisfies condition $(\ast\ast)$. For the former, we note that as $\widetilde{C}\to C$ is a minimal Galois closure, the branch locus of $\widetilde{C}\to D$ is equal to the branch locus of $C\to D$. As $F$ is smooth, Lemma \ref{lem:smoothness of fiber products} implies that the branch loci of $C\to D$ and $E\to D$ are disjoint. Thus, the branch loci of $\widetilde{C}\to D$ and $E\to D$ are disjoint, so Lemma \ref{lem:smoothness of fiber products} gives that $\widetilde{F}$ is smooth, as claimed.
    
    We now consider the commutative diagram
    \[
        \begin{tikzcd}
            \widetilde{F}\arrow{d}\arrow{r}&F\arrow{d}[swap]{g'}\arrow{r}{f'}&E\arrow{d}{g}\\
            \widetilde{C}\arrow{r}&C\arrow{r}{f}&D
        \end{tikzcd}
    \]
    of curves over $K$ with Cartesian squares. The morphism $\widetilde{C}\to D$ is Galois, so by Lemma \ref{lem:Galois base change} pulling back along $g$ induces an isomorphism
    \begin{equation}\label{eq:pullback map galois closure}
        \Aut(\widetilde{C}/D)\iso\Aut(\widetilde{F}/E)
    \end{equation}
    and furthermore the morphism $\widetilde{F}\to E$ is Galois. Under the Galois correspondence, the intermediate curve $\widetilde{F}\to V\to E$ corresponds to a subgroup of $\Aut(\widetilde{F}/E)$. Let $\Gamma\subset\Aut(\widetilde{C}/D)$ denote the preimage of this subgroup under the isomorphism~\eqref{eq:pullback map galois closure} and set $U=\widetilde{C}/\Gamma$. Equivalently, $U$ is the regular curve whose function field is the intersection of the function fields of $\widetilde{C}$ and of $V$ in the function field of $\widetilde{F}$. We have a commutative diagram
    \[
        \begin{tikzcd}
            \widetilde{F}\arrow{d}\arrow{r}&F\arrow{r}\arrow{d}&V\arrow{d}\arrow{r}&E\arrow{d}{g}\\
            \widetilde{C}\arrow{r}&C\arrow{r}&U\arrow{r}&D.
        \end{tikzcd}
    \]
    By construction, we have $\deg(\widetilde{C}/U)=\deg(\widetilde{F}/V)$, from which it follows that each of the squares in this diagram is Cartesian.
\end{proof}

\begin{Proposition}\label{prop:cube theorem with ramification condition}
    Suppose given a commutative diagram
    \begin{equation}\label{eq:cube of curves}
        \begin{tikzcd}[sep = .5 cm]
	         	                & F  \arrow [rr] \arrow[dd]      & 	            & E\arrow[dd,"g"]		 \\
	         F_0\arrow [ur,"\sim"' sloped] \arrow[rr,crossing over] \arrow[dd]     & 	         & E_0\arrow[ur,"\sim"' sloped]   &			\\
	                           	& C \arrow[rr,near end,"f"]      & 	            & D \\
            C_0\arrow[rr]\arrow[ur,dashed]	     	 &           & D_0\arrow[ur,"\sim"' sloped] \arrow[from=uu,crossing over]    &  \\
    \end{tikzcd}
    \end{equation}
    of solid arrows, where the objects are all curves over $K$, the arrows are finite separable morphisms over $K$, the diagonal arrows are all isomorphisms, and the faces of the cube involving only solid arrows satisfy condition $(\ast\ast)$. There exists a unique isomorphism $C_0\iso C$ filling in the dashed arrow in the diagram.
\end{Proposition}
\begin{proof}
    By Galois descent we may assume that $K$ is separably closed. Using our assumption that the solid faces are Cartesian, and that $F$ and $F_0$ are smooth, Lemma \ref{lem:smoothness of fiber products} implies that the branch loci of $f:C\to D$ and $g:E\to D$ are disjoint in $D$, and that the branch loci of $C_0\to D_0$ and $E_0\to D_0$ are disjoint in $D_0$. By Lemma \ref{lem:smoothness of fiber products} the fiber product $V:=C_0\times_DC$ is smooth. Choose a connected component $U\subset V$ and let $X$ be a minimal Galois closure of $U\to D$. We obtain a commutative diagram
    \[
        \begin{tikzcd}
            X\arrow{r}&U\arrow{rr}\arrow{d}&&C\arrow{d}\\
            &C_0\arrow{r}&D_0\arrow{r}&D
        \end{tikzcd}
    \]
    of curves. As we assume that the solid faces of the cube~\eqref{eq:cube of curves} are Cartesian, we have that $\deg(C/D)=\deg(C_0/D_0)$, and therefore $\deg(U/D)\leqslant\deg(C/D)\deg(C_0/D_0)=\deg(C/D)^2$. This implies that $\deg(X/D)\leqslant (\deg(C/D)^2)!$. Furthermore, as $X$ is a minimal Galois closure of $U\to D$, the branch locus of $X\to D$ is equal to that of $U\to D$, which itself is equal to the union of the branch loci of $C_0\to D$ and $C\to D$. Thus, the branch loci of $X\to D$ and of $E\to D$ are disjoint. Lemma \ref{lem:smoothness of fiber products} implies that the fiber product $Y:=X\times_DE$ is smooth. Our assumption that the solid faces of the cube satisfy condition $(\ast\ast)$ implies that $Y$ is also connected, hence integral. By Lemma \ref{lem:Galois base change}, pulling back along $g$ induces an isomorphism
    \begin{equation}\label{eq:pullback map}
        \Aut(X/D)\iso\Aut(Y/E),
    \end{equation}
    and also $Y\to E$ is Galois. It follows that the automorphism $\alpha$ extends to an automorphism, say $\beta$, of $Y$. As~\eqref{eq:pullback map} is surjective, $\beta$ descends to an automorphism, say $\gamma$, of $X$, and we obtain a diagram
    \begin{equation}\label{eq:cube of curves 2}
        \begin{tikzcd}[sep = .5 cm]
	      &   Y\arrow[rr]\arrow[dd]	&                & F  \arrow [rr]\arrow[dd,crossing over]         & 	            & E\arrow[dd,"g"]		 \\
	 Y\arrow[dd]\arrow[rr,crossing over]\arrow[ur,"\beta","\sim"' sloped]&      &  F_0\arrow [ur,"\alpha","\sim"' sloped] \arrow[rr,crossing over]     & 	         & E_0\arrow[ur,"\sim"' sloped]   &			\\
	&X\arrow[rr]  &  &                       	 C \arrow[rr,near end,"f"]        & 	            & D \\
      X\arrow[ur,"\gamma","\sim"' sloped]\arrow[rr] &   &  C_0\arrow[rr]\arrow[ur,dashed]	\arrow[from=uu,crossing over]     	 &           & D_0\arrow[ur,"\sim"' sloped]  \arrow[from=uu,crossing over]   &  \\
    \end{tikzcd}
    \end{equation}
    of solid arrows. We will show that $\gamma$ descends to an isomorphism $C_0\to C$ filling in the dashed arrow in the diagram. The set of isomorphisms $C_0\iso C$ over $D$ may be described in terms of the Galois group $\Aut(X/D)$ as follows. Let $N$ denote the set of isomorphisms $\gamma:X\iso X$ over $D$ whose action by conjugation on $\Aut(X/D)$ takes the subgroup $\Aut(X/C_0)$ to the subgroup $\Aut(X/C)$. The set of isomorphism $C_0\to C$ over $D$ is then identified with the set of orbits of $N$ under the left conjugation action of $\Aut(X/C_0)$ and right conjugation action of $\Aut(X/C)$. As shown above, pulling back along $g$ maps $\Aut(X/D)$ isomorphically to $\Aut(Y/E)$, and furthermore maps $\Aut(X/C_0)$ to $\Aut(Y/F_0)$ and $\Aut(X/C)$ to $\Aut(Y/E)$. From the above description of the set of isomorphisms $C_0\iso C$ over $D$ in terms of the Galois group, we see that pulling back by $g$ induces a bijection between this set and the set of isomorphisms $F_0\to F$ over $E$.
\end{proof}

The following extends Proposition \ref{prop:cube theorem with ramification condition} to families of curves.

\begin{Proposition}\label{prop:cube theorem in families}
    Let $S$ be a $K$-scheme. Suppose given a commutative diagram of the form~\eqref{eq:cube of curves}, where the objects are curves over $S$, the arrows are finite separable morphisms over $S$, the diagonal arrows are all isomorphisms, and for every geometric point $s\in S$, the pullbacks to $s$ of each of the faces of the cube involving only solid arrows satisfy condition $(\ast\ast)$. There exists a unique dotted arrow $C_0\iso C$ over $S$ filling in the diagram.
\end{Proposition}
\begin{proof}
    The proof is similar to that of Proposition \ref{prop:CSI consequence in families}, consisting of a reduction to the pointwise statement of Proposition \ref{prop:cube theorem with ramification condition}. We let $\sIsom_S(C_0,C/D)$ denote the Isom scheme parametrizing isomorphisms $C_0\iso C$ which commute with the maps to $D$, and similarly let $\sIsom_S(F_0,F/E)$ denote the Isom scheme parametrizing isomorphisms $F_0\iso F$ commuting with the maps to $E$. Pullback along $g$ induces a morphism
    \[
        g^*:\sIsom_S(C_0,C/D)\to\sIsom_S(F_0,F/E).
    \]
    We will show that this map is an isomorphism. Because $g$ is faithfully flat, this map is a monomorphism. To complete the proof, we argue as in the proof of Proposition \ref{prop:CSI consequence in families} to show that it suffices to consider the case when $S$ is smooth over $K$. The above Isom schemes are smooth over $S$, so to show that $g^*$ is an isomorphism it suffices to show that $g^*$ is surjective on geometric points, which follows from Proposition \ref{prop:cube theorem with ramification condition}.
\end{proof}

\subsection{Finite morphisms of curves 3: incompressible morphisms}\label{ssec:incompressible morphisms}

In this section we consider the following condition on a finite morphism of curves.

\begin{Definition}
    Let $\lambda$ be an integer. A finite morphism $f:C\to D$ of curves over $K$ is $\lambda$-\textit{incompressible} if for every curve $V$ and factorization $C\to V\to D$ of $f$ such that $V\to D$ is not an isomorphism, we have $g_V\geqslant\lambda$.
    We say that $f$ is \textit{geometrically} $\lambda$-\textit{incompressible} if the base change of $f$ to an algebraic closure of $K$ is $\lambda$-incompressible. Equivalently, $f$ is geometrically $\lambda$-incompressible if for every finite extension $L/K$ the base change $f_L$ is $\lambda$-incompressible.
\end{Definition}

Following ideas of Stichtenoth \cite[Lemma 2]{MR753434}, we show in the following that a pencil on a curve can be refined by a suitable base change to a pencil which is $\lambda$-incompressible for arbitrarily large $\lambda$.

\begin{Lemma}\label{lem:incompressible lemma}
    Let $C$ be a curve over $K$ and let $f:C\to\mathbf{P}^1$ be a finite separable $K$-morphism. If $K$ is finite, we assume also that $f$ is unramified over at least two $K$-points of $\mathbf{P}^1$. For any integer $\lambda$, there exists a finite separable morphism $\varphi:\mathbf{P}^1\to\mathbf{P}^1$ which is totally ramified over some $K$-point of the target, and which has the following property: let $\widetilde{C}\to C$ be any finite \'{e}tale morphism of curves, and define $D$ and $\widetilde{D}$ as the fiber products in the diagram 
    \begin{equation}\label{eq:diagram for incompressible lemma}
        \begin{tikzcd}
            \widetilde{D}\arrow{r}\arrow{d}\arrow[bend left=30]{rr}{\widetilde{f}'}&D\arrow{d}\arrow{r}{f'}&\mathbf{P}^1\arrow{d}{\varphi}\\
            \widetilde{C}\arrow{r}\arrow[bend right=25]{rr}[swap]{\widetilde{f}}&C\arrow{r}{f}&\mathbf{P}^1
        \end{tikzcd}
    \end{equation}
    with Cartesian squares. Then $D$ and $\widetilde{D}$ are smooth proper geometrically integral curves over $K$, and the morphism $\widetilde{f}':\widetilde{D}\to\mathbf{P}^1$ is geometrically $\lambda$-incompressible.
\end{Lemma}
\begin{proof}
    By composing with a suitable automorphism of $\mathbf{P}^1$, we may assume that $f$ is unramified above $0$ and $\infty$. Let $m$ be an integer which is coprime to the characteristic of $K$ (if positive) and such that $m\geqslant\lambda+1$. Let $\varphi:\mathbf{P}^1\to\mathbf{P}^1$ denote the $m$th power map $[x:y]\mapsto [x^m:y^m]$. We claim that $\varphi$ has the desired property. To verify this, let $\widetilde{C}\to C$ be a finite \'{e}tale morphism of curves, and adopt the notation in the diagram~\eqref{eq:diagram for incompressible lemma}. The branch locus of $\varphi$ is disjoint from that of $f$ and $\widetilde{f}$, so Lemma \ref{lem:compositum omnibus lemma} implies that $D$ and $\widetilde{D}$ are smooth proper geometrically integral curves and the two squares in the diagram~\eqref{eq:diagram for incompressible lemma} both satisfy condition $(\ast\ast)$. We will show that $\widetilde{f}'$ is geometrically $\lambda$-incompressible. For this, we may assume that $K$ is algebraically closed, in which case we need to show that $\widetilde{f}'$ is $\lambda$-incompressible. Suppose that $V$ is a curve over $K$ and 
    \[
        \widetilde{D}\to V\to\mathbf{P}^1
    \]
    is a factorization of $\widetilde{f}'$ such that $V\to\mathbf{P}^1$ is not an isomorphism. We will show that $g_V\geqslant\lambda$. We apply Lemma \ref{lem:there is a curve (Galois)} to obtain a curve $U$ and a commutative diagram
    \begin{equation}\label{eq:RH squares}
        \begin{tikzcd}
            \widetilde{D}\arrow{d}\arrow{r}&V\arrow{r}\arrow{d}{\psi}&\mathbf{P}^1\arrow{d}{\varphi}\\
            \widetilde{C}\arrow{r}&U\arrow{r}&\mathbf{P}^1
        \end{tikzcd}
    \end{equation}
    in which both squares are Cartesian. The Riemann--Hurwitz formula \cite[0C1F]{stacks-project} for $\psi$ is
    \[
        2g_V-2=m(2g_{U}-2)+\sum_{P\in V}(e_P-1)
    \]
    where $e_P$ is the ramification index of $\ms O_{V,P}$ over $\ms O_{U,\psi(P)}$. The map $\varphi$ is totally ramified to order $m$ over $0$ and $\infty$. It follows from our assumptions that the branch loci of $U\to\mathbf{P}^1$ and $\varphi$ are disjoint, so there exist at least four closed points of $U$ lying over $0$ or $\infty$. If $Q$ is such a point, then $\psi$ is totally ramified over $Q$ to order $m$, so the fiber $\psi^{-1}(Q)$ consists of a single closed point at which $\psi$ ramifies to order $m$. Combining this observation with the Riemann-Hurwitz formula we obtain the inequalities
     \begin{align*}
        g_V&=mg_{U}-m+1+\frac{1}{2}\sum_{P\in V}(e_P-1)\\
        &\geqslant mg_{U}-m+1+\frac{1}{2}4(m-1)\\
        &\geqslant m-1,
    \end{align*}
    where the final inequality is due to the trivial bound $g_U\geqslant 0$. We chose $m$ so that $m\geqslant\lambda+1$, and thus deduce that $g_V\geqslant\lambda$. We conclude that $\widetilde{f}'$ is $\lambda$-incompressible.
\end{proof}

\subsection{Automorphisms of curves over gerbes}\label{ssec:curves over gerbes}

We now consider relative curves over a gerbe. Let $\ms G$ be a Deligne--Mumford gerbe over $K$ with inertia $G_{\ms G}:=\ms I_{\ms G}$ and structural morphism $\rho:\ms G\to\Spec K$. Consider a diagram
\[
    \begin{tikzcd}
        C\arrow{d}[swap]{\pi}\arrow{r}{f}&D\\
        \ms G&
    \end{tikzcd}
\]
where $C$ and $D$ are curves over $K$, $f$ is a finite separable morphism over $K$, and $\pi$ is a $K$-morphism which is a relative curve.

\begin{Proposition}\label{prop:SES of automorphism groups}
    Suppose that for every field extension $L/K$ and every splitting $s\in\ms G(L)$ the fiber $f_s:\widetilde{C}_L\to D_L$ satisfies condition $(\ast)$ (here, $\widetilde{C}_L=C\times_{\ms G,s}\Spec L$ is the pullback of $C$ along $s$, as in Example \ref{ex:non split gerbe}).
    There is a natural short exact sequence
    \[
        1\to \sAut_{\ms G}(C/D)\to\sAut_{\ms G}(C)\to\rho^{-1}\sAut_{K}(D)
    \]
    of group schemes on $\ms G$.
\end{Proposition}
\begin{proof}
    Write $D_{\ms G}=D\times_{\Spec K}\ms G$. We have a natural isomorphism $\rho^{-1}\sAut_K(D)=\sAut_{\ms G}(D_{\ms G})$. Consider the morphisms
    \[
        \begin{tikzcd}
            &\sAut_{\ms G}(D_{\ms G})\arrow{d}{a}\\
            \sAut_{\ms G}(C)\arrow{r}{b}&\sHom_{\ms G}(C,D_{\ms G})
        \end{tikzcd}
    \]
    of sheaves on $\ms G$, where $a=\_\circ f$ is precomposition with $f$ and $b=f\circ\_$ is postcomposition with $f$. We claim that $a$ is an immersion and that $b$ factors through $a$. Given this, we will obtain a morphism $\sAut_{\ms G}(C)\to\sAut_{\ms G}(D_{\ms G})$ whose kernel is exactly the subgroup $\sAut_{\ms G}(C/D)\subset\sAut_{\ms G}(C)$. Choose an extension $L/K$ and a splitting $s\in\ms G(L)$. The corresponding morphism $s:\Spec L\to\ms G$ is a faithfully flat cover, so it will suffice to verify both these claims after pulling back along $s$. The result then follows from Proposition \ref{prop:CSI consequence}.
\end{proof}

Let $g:E\to D$ be a finite separable morphism of curves over $K$ and set $F=C\times_DE$. We obtain a diagram
    \begin{equation}\label{eq:compositum diagram 2}
        \begin{tikzcd}
            F\arrow{r}{f'}\arrow{d}[swap]{g'}&E\arrow{d}{g}\\
            C\arrow{r}{f}\arrow{d}[swap]{\pi}&D\\
            \ms G.&
        \end{tikzcd}
    \end{equation}

\begin{Lemma}\label{lem:relative curve lemma 1}
    Suppose that $f$ and $g$ have disjoint branch loci, and that at least one of the following two conditions holds.
    \begin{enumerate}
        \item\label{item:fiber product condition 1 gerbe} $g$ is totally ramified over some $K$-point of $D$.
        \item\label{item:fiber product condition 2 gerbe} $|G_{\ms G}|\deg(C/D)<\deg(E/D)$ and $\deg(E/D)$ is a prime.
    \end{enumerate}
     The fiber product $F=C\times_DE$ is a smooth proper geometrically integral curve over $K$, the composition $\pi\circ g':F\to\ms G$ is a relative curve, and the square in the diagram~\eqref{eq:compositum diagram 2} satisfies condition $(\ast\ast)$.
\end{Lemma}
\begin{proof}
    Lemma \ref{lem:compositum omnibus lemma} shows that under either of the given conditions $F$ is a smooth proper geometrically integral curve over $K$ and the square in the diagram~\eqref{eq:compositum diagram 2} satisfies condition $(\ast\ast)$. Let $L/K$ be a finite separable extension splitting $\ms G$ and fix a section $s\in\ms G(L)$. Following the notation of Example \ref{ex:non split gerbe}, we let $\widetilde{C}_L$ and $\widetilde{F}_L$ denote the respective pullbacks of $C$ and $F$ by $s$. As described in Example \ref{ex:non split gerbe, once again}, we obtain a diagram
    \[
        \begin{tikzcd}
            \widetilde{F}_L\arrow{d}\arrow{r}{|G_{\ms G}|}&F_L\arrow{d}[swap]{g'_L}\arrow{r}{f'_L}&E_L\arrow{d}{g_L}\\
            \widetilde{C}_L\arrow{r}{|G_{\ms G}|}&C_L\arrow{r}{f_L}&D_L
        \end{tikzcd}
    \]
    of $L$-schemes with Cartesian squares. The left two horizontal arrows are torsors under the $L$-group scheme $G_L:=s^{-1}G_{\ms G}$ and in particular are finite \'{e}tale of degree $|G_{\ms G}|$. Applying Lemma \ref{lem:compositum omnibus lemma} to the outer rectangle, we deduce that if either of conditions~\eqref{item:fiber product condition 1 gerbe} and~\eqref{item:fiber product condition 2 gerbe} hold then $\widetilde{F}_L$ is a smooth proper geometrically integral curve over $L$. The map $\widetilde{F}_L\to\Spec L$ is the pullback of $F\to\ms G$ along $s$, so in either case $F\to\ms G$ is a relative curve.
\end{proof}

\begin{Proposition}\label{prop:key galois proposition}
    Assume that $F\to\ms G$ is a relative curve and that the square in the diagram~\eqref{eq:compositum diagram 2} satisfies condition $(\ast\ast)$. Pullback along $g$ induces an isomorphism 
    \[
        \sAut_{\ms G}(C/D)\iso\sAut_{\ms G}(F/E).
    \]
\end{Proposition}
\begin{proof}
    Let $L/K$ be a finite separable extension splitting $\ms G$ and fix a section $s\in\ms G(L)$. We adopt the notation of the proof of Lemma \ref{lem:relative curve lemma 1}. The pullback of the map
    \[
        \sAut_{\ms G}(C/D)\to\sAut_{\ms G}(F/E)
    \]
    by $s$ is the map
    \[
        \sAut_L(\widetilde{C}_L/D_L)\to\sAut_L(\widetilde{F}_L/E_L)
    \]
    induced by pullback along $g_L$ (see Example \ref{ex:non split gerbe, once again}). As $g$ is faithfully flat, this map is injective. It remains to show that it is also surjective. As these groups are smooth, it will suffice to assume that $L$ is separably closed and show surjectivity on $L$-points. This follows from Proposition \ref{prop:cube theorem with ramification condition} by taking $C_0=C=\widetilde{C}_L$, $D_0=D=D_L$, $E_0=E=E_L$, and $F_0=F=\widetilde{F}_L$.
\end{proof}

\subsection{Curves over gerbes with prescribed automorphism group, over an infinite field}\label{ssec:curves over gerbes over an infinite field}

Let $\ms G$ be a Deligne--Mumford gerbe over a field $K$ with inertia $G_{\ms G}:=\ms I_{\ms G}$. In this section we prove the partial result that, if $K$ is infinite, then there exist curves over $\ms G$ satisfying conditions~\eqref{item:main curve thm 1} and~\eqref{item:main curve thm 2} of Theorem \ref{thm:main theorem for curves over gerbes}.

 Let $C$ be a curve over $K$ and let $\pi:C\to\ms G$ be a $K$-morphism which is a relative curve.

\begin{Theorem}\label{thm:finite cover of curve over a gerbe, infinite field case}
    Suppose that $K$ is infinite. For any integer $N$, there exists a curve $C'$ over $K$ of genus $\geqslant N$ and a finite separable morphism $C'\to C$ such that the composition $C'\to C\xrightarrow{\pi}\ms G$ is a relative curve, the inertial action on $C'$ induces an isomorphism $G_{\ms G}\iso\sAut_{\ms G}(C')$, and we have $\sAut_K(C')=1$.
\end{Theorem}
\begin{proof}
    By taking a suitable fiber product of a pencil on $C$ and the $m$th power map on $\mathbf{P}^1$, we may assume that $C$ has genus $\geqslant 2$ (this follows from Riemann--Hurwitz, as spelled out below). By Theorem \ref{thm:there exists a curve over a gerbe with a pencil} we may find a finite morphism $f:C\to\mathbf{P}^1$ such that the inertial action induces an isomorphism $G_{\ms G}\iso\sAut_{\ms G}(C/\mathbf{P}^1)$ and $\sAut_K(C/\mathbf{P}^1)=1$. Let $\lambda$ be a positive integer. Let $\varphi:\mathbf{P}^1\to\mathbf{P}^1$ be a morphism which satisfies the conclusions of Lemma \ref{lem:incompressible lemma} applied to $f:C\to\mathbf{P}^1$ and $\lambda$. Let $X$ be a curve with $\sAut_K(X)=1$ equipped with a finite separable morphism $g:X\to\mathbf{P}^1$ which is totally ramified over $\infty$ (Theorem \ref{thm:curves with no auts}). Composing with an automorphism of $\mathbf{P}^1$, we may assume that the branch locus of the composition $\varphi\circ g:X\to\mathbf{P}^1$ is disjoint from the branch locus of $f$. We define $D$ and $E$ to be the fiber products in the diagram
    \begin{equation}\label{eq:a rectangle of curves and a gerbe}
        \begin{tikzcd}
            E\arrow{d}\arrow{r}{f''}&X\arrow{d}{g}\\
            D\arrow{d}\arrow{r}{f'}&\mathbf{P}^1\arrow{d}{\varphi}\\
            C\arrow{r}{f}\arrow{d}[swap]{\pi}&\mathbf{P}^1\\
            \ms G&
        \end{tikzcd}
    \end{equation}
    with Cartesian squares. Condition~\eqref{item:fiber product condition 1 gerbe} of Lemma \ref{lem:relative curve lemma 1} holds for $\varphi$ and for $g$, and therefore $D$ and $E$ are smooth proper geometrically integral curves, the morphisms $D\to\ms G$ and $E\to\ms G$ are relative curves, and the two squares each satisfy condition $(\ast\ast)$. Furthermore, by construction the morphism $f':D\to\mathbf{P}^1$ is $\lambda$-incompressible. In particular, we have that $g_E\geqslant g_D\geqslant\lambda$, so we may assume that the genus of $E$ is arbitrarily large. We claim that if $\lambda$ is sufficiently large then the morphism $E\to C$ has the desired properties.
    
    We consider the inertial action map $G_{\ms G}\to\sAut_{\ms G}(E)$ associated to the morphism $E\to\ms G$. We will show that if $\lambda$ is sufficiently large then this map is an isomorphism. This map factors as the composition 
    \begin{equation}\label{eq:two maps}
        G_{\ms G}\iso\sAut_{\ms G}(C/\mathbf{P}^1)\xrightarrow{(\varphi\circ g)^*}\sAut_{\ms G}(E/X)\hookrightarrow\sAut_{\ms G}(E).
    \end{equation}

    \begin{Claim}\label{claim:claim 1}
        The map $(\varphi\circ g)^*$ is an isomorphism.
    \end{Claim}
    \begin{proof}
        Apply Proposition \ref{prop:key galois proposition} to each of the two squares of the diagram~\eqref{eq:a rectangle of curves and a gerbe}.
    \end{proof}

    \begin{Claim}\label{claim:claim 2}
        If $\lambda$ is sufficiently large then the inclusion $\sAut_{\ms G}(E/X)\hookrightarrow\sAut_{\ms G}(E)$ is an isomorphism.
    \end{Claim}
    \begin{proof}
    We will show that if $\lambda$ satisfies the inequality
    \begin{equation}\label{eq:bound on lambda}
        \lambda\geqslant|G_{\ms G}|^2\deg(C/\mathbf{P}^1)^2+2(g_X-1)|G_{\ms G}|\deg(C/\mathbf{P}^1)+2
    \end{equation}
    then the conditions of Proposition \ref{prop:SES of automorphism groups} hold for the morphisms $f'':E\to X$ and $E\to\ms G$. This will give the result, because then the conclusion of Proposition \ref{prop:SES of automorphism groups} will yield an exact sequence
    \[
        1\to\sAut_{\ms G}(E/X)\to\sAut_{\ms G}(E)\to\rho^{-1}\sAut_{K}(X),
    \]
    where $\rho:\ms G\to\Spec K$ is the structure map. By our choice of $X$, we have $\sAut_K(X)=1$, and so it will follow that the map $\sAut_{\ms G}(E/X)\to\sAut_{\ms G}(E)$ is an isomorphism, as claimed.

    We proceed with the verification of the assumptions of Proposition \ref{prop:SES of automorphism groups}. Let $L/K$ be a field extension splitting $\ms G$ and fix a section $s\in\ms G(L)$. We adopt the notation introduced in Example \ref{ex:non split gerbe}, so that $\widetilde{C}_L,\widetilde{D}_L,$ and $\widetilde{E}_L$ are the respective pullbacks of $C,D,$ and $E$ by $s$. As in Example \ref{ex:non split gerbe, once again}, the diagram~\eqref{eq:a rectangle of curves and a gerbe} yields a commutative diagram
    \begin{equation}\label{eq:big square of curves over L}
        \begin{tikzcd}
            \widetilde{E}_L\arrow{r}{|G_{L}|}\arrow{d}&E_L\arrow{d}\arrow{r}{f''_L}&X_L\arrow{d}{g_L}\\
            \widetilde{D}_L\arrow{r}{|G_{L}|}\arrow{d}&D_L\arrow{r}{f'_L}\arrow{d}&\mathbf{P}^1_L\arrow{d}{\varphi_L}\\
            \widetilde{C}_L\arrow{r}{|G_{L}|}&C_L\arrow{r}{f_L}&\mathbf{P}^1_L
        \end{tikzcd}
    \end{equation}
    of curves over $L$ with Cartesian squares. The three left hand horizontal arrows are torsors under the $L$-group scheme $G_L:=s^{-1}G$, and in particular are finite \'{e}tale of degree $|G_L|=|G_{\ms G}|$. Let $\widetilde{f}''_L:\widetilde{E}_L\to X_L$ (resp. $\widetilde{f}'_L:\widetilde{D}_L\to\mathbf{P}^1_L$) denote the pullback of $f''$ (resp. $f'$) along $s$. Equivalently, these are the horizontal compositions in the top and middle row of the diagram~\eqref{eq:big square of curves over L}. We will show that if $\lambda$ satisfies the inequality~\eqref{eq:bound on lambda} then $\widetilde{f}''_L$ satisfies condition $(\ast)$. To prove this, we consider the Cartesian square
    \[
        \begin{tikzcd}
            \widetilde{E}_L\arrow{r}{\widetilde{f}''_L}\arrow{d}&X_L\arrow{d}{g_L}\\
            \widetilde{D}_L\arrow{r}{\widetilde{f}'_L}&\mathbf{P}^1_L
        \end{tikzcd}
    \]
    appearing in the top two rows of~\eqref{eq:big square of curves over L}. By Lemma \ref{lem:relative curve lemma 1} this square satisfies condition $(\ast\ast)$.
    Suppose given a curve $V_L$ over $L$ and a factorization
    \[
        \widetilde{E}_L\to V_L\to X_L
    \]
    of $\widetilde{f}''_L$ such that $V_L\to X_L$ is not an isomorphism. We apply Lemma \ref{lem:there is a curve (Galois)} to obtain a curve $U_L$ over $L$ and a diagram
    \[
        \begin{tikzcd}
            \widetilde{E}_L\arrow{d}\arrow{r}&V_L\arrow{d}\arrow{r}&X_L\arrow{d}{g_L}\\
            \widetilde{D}_L\arrow{r}&U_L\arrow{r}&\mathbf{P}^1_L
        \end{tikzcd}
    \]
    where both squares are Cartesian. As the map $V_L\to X_L$ is not an isomorphism, the map $U_L\to\mathbf{P}^1_L$ is also not an isomorphism. By our choice of $\varphi$ the map $\widetilde{f}'_L:\widetilde{D}_L\to \mathbf{P}^1_L$ is $\lambda$-incompressible, so we have $g_{U_L}\geqslant\lambda$. Furthermore, we have the equality
    \[
        |G_{\ms G}|\deg(C/\mathbf{P}^1)=|G_L|\deg(E_L/X_L)=\deg(\widetilde{E}_L/X_L).
    \]
    Combining this with the inequality~\eqref{eq:bound on lambda}, we obtain the inequalities
    \begin{align*}
        g_{V_L}\geqslant g_{U_L}\geqslant \lambda&\geqslant|G_{\ms G}|^2\deg(C/\mathbf{P}^1)^2+2(g_X-1)|G_{\ms G}|\deg(C/\mathbf{P}^1)+2\\
        &=\deg(\widetilde{E}_L/X_L)^2+2(g_X-1)\deg(\widetilde{E}_L/X_L)+2\\
        &\geqslant\deg(V_L/X_L)^2+2(g_X-1)\deg(V_L/X_L)+2.
    \end{align*}
    Here, in the last step we have used that $g_X\geqslant 2$, which is a consequence of the fact that $\sAut_K(X)=1$. This shows that $\widetilde{f}''_L$ satisfies condition $(\ast)$. It follows that the assumptions of Proposition \ref{prop:SES of automorphism groups} hold, which completes the proof of the claim.
\end{proof}

Combining Claims \ref{claim:claim 1} and \ref{claim:claim 2}, we deduce that if $\lambda$ is sufficiently large then the maps~\eqref{eq:two maps} are isomorphisms, and therefore the inertial action induces an isomorphism $G_{\ms G}\iso\sAut_{\ms G}(E)$. It remains only to show the following.
    
\begin{Claim}\label{claim:claim 3}
    If $\lambda$ is sufficiently large then we have $\sAut_K(E)=1$.
\end{Claim}
\begin{proof}
    The proof of this is similar to the proof that $G_{\ms G}\cong\sAut_{\ms G}(E)$, but with the gerbe $\ms G$ replaced with $\Spec K$. We consider the composition
    \[
        \sAut_K(C/\mathbf{P}^1)\xrightarrow{(\varphi\circ g)^*}\sAut_K(E/X)\hookrightarrow\sAut_K(E).
    \]
    By assumption, we have $\sAut_K(C/\mathbf{P}^1)=1$, and by Proposition \ref{prop:key galois proposition} applied with $\ms G=\Spec K$ being the trivial gerbe, we have that $(\varphi\circ g)^*$ is an isomorphism. Suppose that $\lambda$ satisfies the inequality
    \[
        \lambda\geqslant\deg(C/\mathbf{P}^1)^2+2(g_X-1)\deg(C/\mathbf{P}^1)+2.
    \]
    Arguing as in the proof of Claim \ref{claim:claim 2}, we deduce that for any field extension $L/K$ the base change $f''_{L}$ satisfies condition $(\ast)$. Applying Proposition \ref{prop:SES of automorphism groups} with $\ms G$ being the trivial $K$-gerbe $\Spec K$, we get an exact sequence
    \[
        1\to\sAut_K(E/X)\to\sAut_K(E)\to\sAut_K(X).
    \]
    We have already shown that $\sAut_K(E/X)=1$, and by construction $\sAut_K(X)=1$, so $\sAut_K(E)=1$, as claimed.
\end{proof}
\end{proof}

\subsection{Descent of curves}
\label{ssec:descent}

Our next task in pursuit of the proof of Theorem \ref{thm:finite cover of curve over a gerbe} is to produce curves with a prescribed field of definition, while simultaneously preserving the properties we have already ensured. As preparation, in this section we will prove that arbitrary extensions of fields are of effective descent for curves. This is presumably well known, but we include the proof here for lack of a suitable reference. We will also take this opportunity to establish some notation.

Let $K$ be a field. Let $L\subset K$ be a subfield and write $\pi_i:\Spec K\otimes_LK\to\Spec K$ ($i=1,2$) for the two projections.

\begin{Definition}
    Let $C$ be a curve over $K$. A \emph{descent datum} for $C$ with respect to the extension $K/L$ is an isomorphism $\varphi_C:\pi_1^{-1}C\iso\pi_2^{-1}C$ of schemes over $\Spec K\otimes_LK$ which satisfies the cocycle condition \cite[023V]{stacks-project}. A \emph{descent datum in curves} is a pair $(C,\varphi_C)$, where $C$ is a curve over $K$ and $\varphi_C$ is a descent datum for $C$. A \emph{morphism} between two such pairs $(C,\varphi_C)$ and $(D,\varphi_D)$ is a morphism $f:C\to D$ over $K$ such that $\varphi_D\circ(\pi_1^{-1}f)=(\pi_2^{-1}f)\circ\varphi_C$.
\end{Definition}

Given a curve $C_L$ over $L$, the pullback $C_K:=C_L\otimes_LK$ is equipped with a canonical descent datum with respect to the extension $K/L$ \cite[023Z]{stacks-project}. Furthermore, given a morphism $f_L:C_L\to D_L$ over $L$, the pullback $f_K:=f_L\otimes_LK:C_K\to D_K$ yields a morphism of the canonical descent data. We obtain a functor from the category of curves over $L$ and $L$-morphisms to the category of descent data in curves with respect to $K/L$. A descent datum $(C,\varphi_C)$ is said to be \emph{effective} if it is in the essential image of this functor.

\begin{Theorem}[Effective descent for curves]
\label{thm:effective descent!}
    The above functor defines an equivalence between the category of curves over $L$ and $L$-morphisms and the category of descent data in curves with respect to $K/L$. In particular, every descent datum in curves is effective.
\end{Theorem}
\begin{proof}
    The morphism $\Spec K\to\Spec L$ is an fpqc cover. Thus, as a special case of fpqc descent for morphisms of schemes \cite[02W0]{stacks-project}, our pullback functor is fully faithful. It remains to show that the functor is also essentially surjective, that is, that every descent datum in curves is effective. For this, let $C$ be a curve over $K$ and let $\varphi_C$ be a descent datum for $C$ with respect to $K/L$. Following an argument of Bosch--L\"{u}tkebohmert--Raynaud \cite[pg. 160]{MR1045822}, we may reduce to the case when $K/L$ is finite. We then use that every finite set of points of $C$ is contained in an open affine to find an open affine cover of $C$ which is $\varphi_C$-invariant. By fpqc descent for affine schemes \cite[0247]{stacks-project} we may descend this cover to $L$. Using fpqc descent for morphisms of schemes, the morphisms gluing together this cover also descend to $L$. Thus, the descent datum $(C,\varphi_C)$ is effective.
\end{proof}

We will use the following terminology.

\begin{Definition}\label{def:field of definition of some things}
    Let $C$ be a curve over $K$. A \emph{descent} of $C$ to $L$ is a pair $(C_L,\alpha)$, where $C_L$ is a curve over $L$ and $\alpha:C\iso C_L\otimes_LK$ is an isomorphism over $K$. We say that $C$ is \textit{defined over} $L$ if there exists a descent of $C$ to $L$. An \emph{isomorphism} $(C_L,\alpha)\to (C'_L,\alpha')$ between two such pairs consists of an isomorphism $\theta:C_L\iso C'_L$ over $L$ such that $\alpha'=(\theta\otimes_LK)\circ\alpha$.

    Let $f:C\to D$ be a morphism of curves over $K$. A \textit{descent} of $f$ to $L$ is a tuple $(f_L,C_L,\alpha_C,D_L,\alpha_D)$, where $(C_L,\alpha_C)$ is a descent of $C$ to $L$, $(D_L,\alpha_D)$ is a descent of $D$ to $L$, and $f_L:C_L\to D_L$ is a morphism over $L$ such that $\alpha_D\circ f=(f_L\otimes_LK)\circ\alpha_C$. We say that $f$ is \textit{defined over} $L$ if there exists a descent of $f$ to $L$.
\end{Definition}

Given a subfield $L\subset K$, a curve $C$ over $K$, and a descent $(C_L,\alpha_C)$ of $C$ to $L$, the canonical descent datum on $C_L\otimes_LK$ gives rise via the isomorphism $\alpha_C$ to a descent datum on $C$. This association is functorial, and so describes a functor from the category of descents of $C$ to $L$ to the (discrete) category of descent data for $C$ with respect to $K/L$. The following is an immediate consequence of Theorem \ref{thm:effective descent!}

\begin{Corollary}
\label{cor:effective descent for morphisms}
    Let $C$ and $D$ be curves over $K$. Let $(C_L,\alpha_C)$ and $(D_L,\alpha_D)$ be descents of $C$ and $D$ to $L$ and let $\varphi_C$ and $\varphi_D$ be the induced descent data on $C$ and $D$.
    \begin{enumerate}
        \item\label{item:descent for curves 1} Pullback defines an equivalence between the category of descents of $C$ to $L$ and the discrete category of descent data for $C$ with respect to $K/L$.
        \item\label{item:descent for curves 2} Pullback defines a bijection between the set of $L$-morphisms $f_L:C_L\to D_L$ such that $\alpha_D\circ f=(f_L\otimes_LK)\circ\alpha_C$ and the set of $K$-morphisms $f:C\to D$ such that $\varphi_D\circ(\pi_1^{-1}f)=(\pi_2^{-1}f)\circ\varphi_C$.
    \end{enumerate}
\end{Corollary}

\subsection{Finite morphisms of curves 4: fields of definition}\label{ssec:fields of definition}

In this section we will prove some results which under certain conditions relate the possible fields of definition of curves connected by finite morphisms. Fix a finitely generated field extension $K/K_0$. Our first result shows that a curve over $K$ with trivial automorphism group scheme has a unique minimal field of definition.

\begin{Lemma}\label{lem:unique descent}
    Let $C$ be a curve over $K$ such that $\sAut_K(C)=1$. 
    \begin{enumerate}
        \item\label{item:no auts 1} If $L\subset K$ is a subfield and $(C_L,\alpha)$ and $(C_L',\alpha')$ are two descents of $C$ to $L$, then there exists a unique isomorphism $\theta:C_L\iso C_L'$ over $L$ such that $\alpha'=(\theta\otimes_LK)\circ\alpha$.
        \item\label{item:no auts 2} There exists a unique subfield $M_C\subset K$ containing $K_0$ over which $C$ is defined, and which is contained in every other subfield of $K$ containing $K_0$ over which $C$ is defined.
    \end{enumerate}
\end{Lemma}
\begin{proof}
    Our assumption that $\sAut_K(C)=1$ implies that
    \[
        \sAut_{\Spec K\otimes_LK}(\pi_1^{-1}C)=\pi_1^{-1}\sAut_K(C)=1.
    \]
    Thus, there is a \emph{unique} isomorphism $\varphi_C:\pi_1^{-1}C\iso\pi_2^{-1}C$ over $\Spec K\otimes_LK$. It follows that the descent data on $C$ induced by any two descents of $C$ to $L$ are equal, and~\eqref{item:no auts 1} follows from Corollary \ref{cor:effective descent for morphisms} ~\eqref{item:descent for curves 1}.
    
    We now show~\eqref{item:no auts 2}. Our assumption that $\sAut_K(C)=1$ implies that $C$ has genus $g\geqslant 2$. Thus the moduli stack $\ms M_g$ of genus $g$ curves over $K_0$ is Deligne--Mumford. We define $M_C$ to be the residue field of the topological point $x_C\in|\ms M_g|$ represented by $C$. We have inclusions $K_0\subset M_C\subset K$. As $C$ has trivial automorphism group scheme, $\Spec M_C$ is the residual gerbe of $\ms M_g$ at the point $x_C$, so $C$ is defined over $M_C$. By the the universal property of the residual gerbe $M_C$ is contained in every other field of definition for $C$. 
\end{proof}

We now prove some consequences for fields of definition of condition $(\ast)$ (Definition \ref{def:condition star}).

\begin{Proposition}\label{prop:step 1}
    Let $f:C\to D$ be a finite separable morphism of curves over $K$. Suppose that for every field extension $K'/K$, the base change $f_{K'}=f\otimes_KK'$ satisfies condition $(\ast)$. If $C$ is defined over a subfield $L\subset K$, then $f$ is also defined over $L$ in a compatible way. More precisely, if $(C_L,\alpha_C)$ is a descent of $C$ to $L$, then there exists a descent $(D_L,\alpha_D)$ of $D$ to $L$ and a finite morphism $f_L:C_L\to D_L$ over $L$ such that $\alpha_D\circ f=(f_L\otimes_LK)\circ\alpha_C$.
\end{Proposition}
\begin{proof}
    This is a consequence of Proposition \ref{prop:CSI consequence in families} and descent. Let $\pi_i:\Spec K\otimes_LK\to\Spec K$ $(i=1,2)$ be the two projections. Let $(C_L,\alpha_C)$ be a descent of $C$ to $L$ and let $\varphi_C:\pi_1^{-1}C\iso\pi_2^{-1}C$ be the resulting descent datum on $C$. Applying Proposition \ref{prop:CSI consequence in families} with $S=\Spec K\otimes_LK$, we deduce that there exists a unique isomorphism $\varphi_D:\pi^{-1}_1D\iso\pi^{-1}_2D$ over $S$ such that the diagram
    \[
        \begin{tikzcd}
            \pi_1^{-1}C\arrow{d}[swap]{\pi_1^{-1}f}\arrow{r}{\varphi_C}[swap]{\sim}&\pi_2^{-1}C\arrow{d}{\pi_2^{-1}f}\\
            \pi_1^{-1}D\arrow{r}{\varphi_D}[swap]{\sim}&\pi_2^{-1}D
        \end{tikzcd}
    \]
    commutes. Applying the uniqueness statement of Proposition \ref{prop:CSI consequence in families} over $\Spec K\otimes_LK\otimes_LK$ and using that $\varphi_C$ satisfies the cocycle condition, we see that $\varphi_D$ must also satisfy the cocycle condition. Applying parts~\eqref{item:descent for curves 1} and~\eqref{item:descent for curves 2} of Corollary \ref{cor:effective descent for morphisms}, we obtain a descent $(D_L,\alpha_D)$ of $D$ to $L$ and a morphism $f_L:C_L\to D_L$ descending $f$.
\end{proof}

Consider a commuting square
\begin{equation}\label{eq:balanced square 2}
    \begin{tikzcd}
        F\arrow{d}[swap]{g'}\arrow{r}{f'}&E\arrow{d}{g}\\
        C\arrow{r}{f}&D
    \end{tikzcd}
\end{equation}
where $C,D,E,$ and $F$ are curves over $K$ and the arrows are all finite separable $K$-morphisms. We prove the following consequence for fields of definition of condition $(\ast\ast)$ (Definition \ref{def:balanced}).

\begin{Proposition}\label{prop:step 2}
    Assume that the square~\eqref{eq:balanced square 2} satisfies condition $(\ast\ast)$ and that $\sAut_K(E)=1$. If $L\subset K$ is a subfield over which both $f'$ and $g$ are defined, then the entire commutative square~\eqref{eq:balanced square 2} is also defined over $L$. Moreover, we may choose a descent of~\eqref{eq:balanced square 2} to include any given descents of $f'$ and of $g$ up to isomorphism. More precisely, suppose $(f'_L,F_L,\alpha_F,E_L,\alpha_E)$ is a descent of $f'$ to $L$ and $(g_L,\widetilde{E}_L,\widetilde{\alpha}_E,D_L,\alpha_D)$ is a descent of $g$ to $L$. By Lemma \ref{lem:unique descent}~\eqref{item:no auts 1}, the descent of $E$ to $L$ is unique up to unique isomorphism, so after replacing $g_L$ with its composition with a uniquely determined isomorphism, we may assume that $\widetilde{E}_L=E_L$ and $\widetilde{\alpha}_E=\alpha_E$. Then, under the given assumptions, there exists a curve $C_L$ over $L$, morphisms $f_L:C_L\to D_L$ and $g'_L:F_L\to C_L$ over $L$ such that the diagram
    \[
        \begin{tikzcd}
            F_L\arrow{r}{f'_L}\arrow{d}[swap]{g'_L}&E_L\arrow{d}{g_L}\\
            C_L\arrow{r}{f_L}&D_L
        \end{tikzcd}
    \]
    commutes, and an isomorphism $\alpha_C:C\iso C_L\otimes_LK$ such that the diagram
    \[
        \begin{tikzcd}
            F_L\otimes_LK\arrow{rrr}{f'_L\otimes_LK}\arrow{ddd}[swap]{g'_L\otimes_LK}&&&E_L\otimes_LK\arrow{ddd}{g_L\otimes_K}\\
            &F\arrow{r}{f'}\arrow{d}[swap]{g'}\arrow[ul,"\sim"{sloped,auto,swap},"\alpha_F"{auto,swap}]&E\arrow{d}{g}\arrow[ur,"\sim"{sloped,auto,swap},"\alpha_E"{auto}]&\\
            &C\arrow{r}{f}\arrow[dl,"\sim"{sloped,auto,swap},"\alpha_C"{auto,swap}]&D\arrow[dr,"\sim"{sloped,auto,swap},"\alpha_D"{auto}]&\\
            C_L\otimes_LK\arrow{rrr}{f_L\otimes_LK}&&&D_L\otimes_LK
        \end{tikzcd}
    \]
    commutes.
\end{Proposition}
\begin{proof}
    Let $(f'_L,F_L,\alpha_F,E_L,\alpha_E)$ be a descent of $f'$ to $L$ and let $(g_L,E_L,\alpha_E,D_L,\alpha_D)$ be a compatible descent of $g$ to $L$. Let $\pi_i:\Spec K\otimes_LK\to\Spec K$ ($i=1,2$) be the two projections. Let $\varphi_F,\varphi_E,$ and $\varphi_D$ be the induced descent data for the curves $F,E,$ and $D$ with respect to $K/L$. We obtain a diagram
     \[
        \begin{tikzcd}[sep = .5 cm]
	         	                & \pi_2^{-1}F  \arrow [rr]\arrow[dd,crossing over]         & 	            & \pi_2^{-1}E\arrow[dd,"\pi_2^{-1}g"]		 \\
	         \pi_1^{-1}F\arrow [ur,"\varphi_F"] \arrow[rr,crossing over] \arrow[dd]     & 	         & \pi_1^{-1}E\arrow[ur,"\varphi_E"]   &			\\
	                           	& \pi_2^{-1}C \arrow[rr,near end,"\pi_2^{-1}f"]        & 	            & \pi_2^{-1}D \\
            \pi_1^{-1}C\arrow[rr,"\pi_1^{-1}f"]\arrow[ur,dotted,"\varphi_C"]	     	 &           & \pi_1^{-1}D\arrow[ur,"\varphi_D"] \arrow[from=uu,crossing over,near start,swap,"\pi_1^{-1}g"]    &  \\
        \end{tikzcd}
    \]
    of solid arrows over the base scheme $S=\Spec K\otimes_LK$. This diagram satisfies the assumptions of Proposition \ref{prop:cube theorem in families}, so there exists a unique isomorphism $\varphi_C:\pi_1^{-1}C\iso\pi_2^{-1}C$ filling in the diagram. The solid diagonal arrows all satisfy the cocycle condition, so $\varphi_C$ does as well. Applying Corollary \ref{cor:effective descent for morphisms} we find the desired descents of $C$, $f$, and $g'$ to $L$.
\end{proof}

The following result shows that a curve with trivial automorphism group scheme over an infinite field admits pencils with prescribed fields of definition.

\begin{Theorem}\label{thm:fields of definition of pencils}
    Let $K$ be an infinite field and let $K/K_0$ be a finitely generated extension. If $C$ is a curve over $K$ such that $\sAut_K(C)=1$, then there exists a finite separable $K$-morphism $f:C\to\mathbf{P}^1$ such that, if $L\subset K$ is any proper subfield containing $K_0$, then there does not exist a morphism of curves over $L$ with target $\mathbf{P}^1_L$ descending $f$.
\end{Theorem}
\begin{proof}
    By Lemma \ref{lem:unique descent} there exists a unique minimal field of definition for $C$ containing $K_0$. By a change of notation we may assume that this field is $K_0$. Fix a curve $C_0$ over $K_0$ descending $C$. To simplify the notation, we will assume that $C=C_0\otimes_{K_0}K$. Let $\ms L_0$ be a very ample invertible sheaf on $C_0$ and write $\ms L$ for the base change of $\ms L_0$ to $C$. We will show that if the degree of $\ms L_0$ is sufficiently large then there exists a basepoint free pencil $\mathbf{L}\subset|\ms L|$ such that for any choice of isomorphism $\mathbf{L}\cong\mathbf{P}^1$ the resulting morphism $f_{\mathbf{L}}:C\to\mathbf{P}^1$ does not descend over any proper subfield of $K$ containing $K_0$ to a morphism to the projective line.
    
    Let $\Gr_{K_0}$ denote the Grassmannian classifying lines in the projective space $|\ms L_0|$. The base change $\Gr_{K}=\Gr_{K_0}\otimes_{K_0}K$ is then identified with the Grassmannian classifying lines in $|\ms L|$.
    
    \begin{Claim}\label{claim:pencils 1}
        Consider field extensions $K/L/K_0$. Let $\mathbf{L}\subset|\ms L|$ be a basepoint free pencil, let $\mathbf{L}\cong\mathbf{P}^1$ be an isomorphism, and write $f=f_{\mathbf{L}}:C\to\mathbf{P}^1$ for the resulting finite $K$--morphism. The morphism $f$ descends to a morphism over $L$ with target $\mathbf{P}^1_L$ if and only if the $K$-point $[\mathbf{L}]\in\Gr_{K_0}(K)$ classifying $\mathbf{L}$ is in the image of the map $\Gr_{K_0}(L)\to\Gr_{K_0}(K)$.
    \end{Claim}
    \begin{proof}
        We write $C_L:=C_0\otimes_{K_0}L$ and let $\ms L_L$ denote the base change of $\ms L_0$ to $C_L$.
        
        Suppose that the $K$-point $[\mathbf{L}]\in\Gr_{K_0}(K)$ is in the image of the map $\Gr_{K_0}(L)\to\Gr_{K_0}(K)$. Then $\mathbf{L}$ descends to a pencil say $\mathbf{L}_L\subset|\ms L_L|$. As $\mathbf{L}$ was assumed basepoint free, $\mathbf{L}_L$ is also basepoint free. Given any choice of isomorphism $\mathbf{L}_L\cong\mathbf{P}^1_{L}$, the resulting morphism $C_L\to\mathbf{P}^1_L$ gives a descent of $f$ to $L$ with target $\mathbf{P}^1_L$.
        
        Conversely, suppose that $f$ descends to a morphism over $L$ with target $\mathbf{P}^1_L$. By Lemma \ref{lem:unique descent}, a descent of $C$ to a subfield is unique up to unique isomorphism, so $f$ descends to a morphism $f_L:C_L\to\mathbf{P}^1_L$. As $C$ has trivial automorphism group, the condition that $f_L$ is a descent of $f$ implies that we have a $K$-automorphism $\theta$ of $\mathbf{P}^1$ and a commutative diagram
        \begin{equation}\label{eq:descent triangle}
            \begin{tikzcd}
                C\arrow{d}[swap]{f_L\otimes_LK}\arrow{dr}{f}&\\
                \mathbf{P}^1\arrow{r}{\sim}[swap]{\theta}&\mathbf{P}^1
            \end{tikzcd}
        \end{equation}
        over $K$. We claim that $\ms L_L\cong f_L^*\ms O_{\mathbf{P}^1_L}(1)$. Indeed, consider the sheaf $\ms F$ over $L$ whose sections over an $L$-scheme $T$ are isomorphisms $\pi^*_T\ms L_L\cong f_T^*\ms O_{\mathbf{P}^1_T}(1)$ of invertible sheaves over $C_T=C_L\times_{\Spec L}T$, where $\pi_T:C_T\to C_L$ is the projection and $f_T=f_L\times_{\Spec L}\id_T:C_T\to\mathbf{P}^1_T$ is the base change of $f_L$. Then $\ms F$ carries a natural simply transitive action of $\mathbf{G}_m$. The above diagram shows that $\ms L\cong f^*\ms O_{\mathbf{P}^1}(1)\cong (f_L\otimes_LK)^*\ms O_{\mathbf{P}^1}(1)$, so $\ms F(K)$ is nonempty. Thus, $\ms F$ is a $\mathbf{G}_m$-torsor over $L$. Hilbert's theorem 90 implies that $\ms F(L)$ is nonempty, which proves the claim. 

        Fix an isomorphism $\ms L_L\cong f_L^*\ms O_{\mathbf{P}^1_L}(1)$. Under this isomorphism, the morphism $f_L$ corresponds to a pencil say $\mathbf{L}_L$ in the projective space $|\ms L_L|$ together with an isomorphism $\mathbf{L}_L\cong\mathbf{P}^1_L$ over $L$. The commutativity of the above diagram~\eqref{eq:descent triangle} shows that the base change of $\mathbf{L}_L$ to $K$ is exactly the line $\mathbf{L}\subset|\ms L|$.
    \end{proof}

    Set $r=\dim_{K_0}\H^0(C_0,\ms L_0)$ and $N=2(r-2)$. The Grassmannian $\Gr_{K_0}$ is a rational $K_0$-variety of dimension $N$.
    
    \begin{Claim}\label{claim:pencils 2}
        Suppose that $K$ may be generated as a field extension of $K_0$ by $\leqslant N$ elements. If $U_{K_0}\subset\Gr_{K_0}$ is a dense open subvariety, then there exists a $K$-point of $U_{K_0}$ which is not in the image of the inclusion
        \[
            \left(\bigcup_{K\supsetneq L\supset K_0}\Gr_{K_0}(L)\right)\subset\Gr_{K_0}(K).
        \]
    \end{Claim}
    \begin{proof}
        The open subvariety $U_{K_0}$ is rational of dimension $N$. By shrinking $U_{K_0}$ we may assume there exists an open immersion $U_{K_0}\subset\mathbf{A}_{K_0}^{N}$. If $\lambda_1,\ldots,\lambda_{N}\in K$ are scalars that generate $K$ as a field extension of $K_0$, then the $K$-point $(\lambda_1,\ldots,\lambda_{N})\in\mathbf{A}^N_K(K)$ is not in the image of the map 
        \[
            \mathbf{A}^N_{K_0}(L)\to\mathbf{A}^{N}_{K_0}(K)=\mathbf{A}^N_K(K)
        \]
        for any proper subfield $L\subset K$ containing $K_0$. The translates $(\lambda_1+m_1,\ldots,\lambda_N+m_N)$ for integers $m_1,\ldots,m_N\in\mathbf{Z}$ have the same property. Using our assumption that $K$ is infinite, a straightforward induction shows that the set of such points is Zariski dense in $\mathbf{A}^N_K$, so we may find one which is contained in the open subvariety $U_K:=U_{K_0}\otimes_{K_0}K$. Under the canonical bijection $U_K(K)=U_{K_0}(K)$ this point gives a $K$-point of $U_{K_0}$ with the desired property.
    \end{proof}

    We now complete the proof of the theorem. Consider the dense open subvariety $U_{K_0}\subset\Gr_{K_0}$ parameterizing pencils which are basepoint free and separable. By Claim \ref{claim:pencils 2}, there exists a $K$-point of $U_{K_0}$ which is not the image of an $L$-point of $\Gr_{K_0}$ for any proper subfield $L\subset K$ containing $K_0$. By Claim \ref{claim:pencils 1}, if $\mathbf{L}\subset|\ms L|$ is the corresponding pencil over $K$, then for any choice of isomorphism $\mathbf{L}\cong\mathbf{P}^1$ the resulting $K$-morphism $f:C\to\mathbf{P}^1$ has the desired property.
\end{proof}

\subsection{Curves over gerbes with prescribed field of definition, over an infinite field}\label{ssec:curves over gerbes over an infinite field + field of definition}

In this section we will prove Theorem \ref{thm:finite cover of curve over a gerbe} in the case of an infinite field. This extends Theorem \ref{thm:finite cover of curve over a gerbe, infinite field case} to include consideration of the field of definition. In particular, this result will complete the proof of Theorem \ref{thm:main theorem for curves over gerbes} in the case of an infinite field.

We recall the notation: $K/K_0$ is a finitely generated field extension, $\ms G$ is a Deligne--Mumford gerbe over $K$ with inertia $G_{\ms G}:=\ms I_{\ms G}$, $C$ is a curve over $K$, and $\pi:C\to\ms G$ is a morphism which is a relative curve.

\begin{Theorem}\label{thm:finite cover of curve over a gerbe, field of definition, infinite field case}
    Theorem \ref{thm:finite cover of curve over a gerbe} holds if $K$ is infinite.
\end{Theorem}
\begin{proof}
    By Theorem \ref{thm:finite cover of curve over a gerbe, infinite field case}, we may assume by replacing $C$ with a finite cover that the inertial action induces an isomorphism $G_{\ms G}\iso\sAut_{\ms G}(C)$ and that $\sAut_K(C)=1$. We will produce the desired curve by taking a further finite cover, using essentially the same construction as in the proof of Theorem \ref{thm:finite cover of curve over a gerbe, infinite field case}. We follow the proof closely and indicate where modifications are needed. Using Theorem \ref{thm:fields of definition of pencils}, we choose a finite separable $K$-morphism $f:C\to\mathbf{P}^1$ which does not descend over any proper intermediate extension $K/L/K_0$ to a morphism to the projective line. By our initial choice of $C$, the composition
    \[
        G_{\ms G}\hookrightarrow\sAut_{\ms G}(C/\mathbf{P}^1)\hookrightarrow\sAut_{\ms G}(C)
    \]
    is an isomorphism, so we automatically have $G_{\ms G}\cong\sAut_{\ms G}(C/\mathbf{P}^1)$, and because $\sAut_K(C)=1$ we also have $\sAut_K(C/\mathbf{P}^1)=1$. As before, we choose a morphism $\varphi:\mathbf{P}^1\to\mathbf{P}^1$, a curve $X$ over $K$ such that $\sAut_K(X)=1$, and a finite morphism $g:X\to\mathbf{P}^1$ which is totally ramified over $\infty$. We may assume that $X$, as well as the morphisms $g$ and $\varphi$, are defined over $K_0$ (or even over the prime field of $K$), and that the branch loci of $\varphi\circ g$ and $f$ are disjoint.
    We adopt the notation of the diagram~\eqref{eq:a rectangle of curves and a gerbe}. In particular, we consider the commutative square
    \begin{equation}\label{eq:square field of def}
        \begin{tikzcd}
            E\arrow{d}\arrow{r}{f''}&X\arrow{d}{\varphi\circ g}\\
            C\arrow{r}{f}&\mathbf{P}^1
        \end{tikzcd}
    \end{equation}
    appearing as the outer rectangle. As before, we may ensure that $E$ has arbitrarily large genus. We claim that the finite cover $E\to C$ of $C$ has the desired properties. We have only added properties to the construction, so the morphism $E\to\ms G$ still satisfies the conclusions of Theorem \ref{thm:finite cover of curve over a gerbe, infinite field case}. It therefore remains only to show that $E$ is not defined over any proper subfield of $K$ containing $K_0$.
    
    For this, we note that, as before, the squares of~\eqref{eq:a rectangle of curves and a gerbe} satisfy condition $(\ast\ast)$, and therefore the above commutative square~\eqref{eq:square field of def} also satisfies condition $(\ast\ast)$. Furthermore, as shown in the proof of Claim \ref{claim:claim 3}, the morphism $f''$ satisfies condition $(\ast)$, as does its base change to any field containing $K$. Now, suppose that $E$ is defined over an intermediate extension $K/L/K_0$. By Proposition \ref{prop:step 1}, $f''$ is also defined over $L$ in a compatible way. The composition $\varphi\circ g$ descends to a morphism over $K_0$ with target $\mathbf{P}^1_{K_0}$, and $\sAut_K(X)=1$, so by Proposition \ref{prop:step 2} applied to the square~\eqref{eq:square field of def}, $f$ is compatibly defined over $L$, and hence descends to a morphism over $L$ with target $\mathbf{P}^1_L$. By our initial choice of $f$, we conclude that $K=L$.
\end{proof}

\subsection{Curves over gerbes, over a finite field}\label{ssec:curves over gerbes with prescribed field of definition, finite field case}

The goal of this section is to prove Theorem \ref{thm:finite cover of curve over a gerbe} over a finite field, thereby completing the proof of Theorem \ref{thm:main theorem for curves over gerbes} in general. The strategy of proof is broadly similar to the infinite field case, but requires combining the ingredients in a different way. In fact, the finite field case is in a sense easier than the infinite case, owing to the following result.

\begin{Proposition}\label{prop:invertible sheaf of degree 1}
    If $C$ is a curve over a finite field $K$, then $C$ admits an invertible sheaf of degree 1.
\end{Proposition}
\begin{proof}
    We consider the Leray spectral sequence for the sheaf $\mathbf{G}_m$ and the morphism $C\to\Spec K$, which yields an exact sequence
    \[
        0\to\Pic(C)\to\uPic_{C/K}(K)\to\Br(K).
    \]
    By Lang's theorem, $\uPic^1_{C/K}$ has a $K$-point. Furthermore, by Tsen's theorem, we have $\Br(K)=0$, so any choice of $K$-point of $\uPic^1_{C/K}$ lifts to an element of $\Pic(C)$ that classifies an invertible sheaf of degree 1 on $C$.
\end{proof}

There are two essential places in the proofs of Theorems \ref{thm:finite cover of curve over a gerbe, infinite field case} and \ref{thm:finite cover of curve over a gerbe, field of definition, infinite field case} where the assumption that $K$ was infinite is used: first, in the selection of a pencil $f:C\to\mathbf{P}^1$ such that $G\cong\sAut_{\ms G}(C/\mathbf{P}^1)$, and second in the consequence that $\mathbf{P}^1$ has infinitely many $K$-points, which allowed us to ensure that the branch loci of various morphisms were disjoint. In the following result we use some particular properties of curves over finite fields to select a pencil which will allow us to avoid both of these problems. We will use the following terminology.

\begin{Definition}
     A closed point $P$ of a curve $C$ over $K$ is \textit{defined over} a subfield $L\subset K$ if there exists a descent $(C_L,\alpha)$ of $C$ to $L$ and a closed point $P_L\in C_L$ such that the isomorphism $\alpha:C\iso C_L\otimes_LK$ maps $P$ isomorphically to the subscheme $P_L\otimes_LK\subset C_L\otimes_LK$.
\end{Definition}

\begin{Proposition}\label{prop:pencils with auts over finite fields}
     Let $C$ be a curve over a finite field $K$. For any positive integer $N$, there exists a finite separable morphism $f:C\to\mathbf{P}^1$ such that
     \begin{enumerate}
         \item\label{item:pencil 1} the degree of $f$ is a prime and is $\geqslant N$,
         \item\label{item:pencil 2} there exists a $K$-point of $\mathbf{P}^1$ whose preimage in $C$ contains two closed points at which $f$ ramifies to different degrees,
         \item\label{item:pencil 3} $f$ is ramified only over $\infty$, and
         \item\label{item:pencil 4} $f$ is not defined over any proper subfield of $K$.
     \end{enumerate}
\end{Proposition}
\begin{proof}
    By Proposition \ref{prop:invertible sheaf of degree 1} we may find an invertible sheaf $\ms L$ on $C$ of degree 1. By a result of Kedlaya \cite[Theorem 1]{MR2092132}, we may find an integer $d$ and sections $s_0,s_1\in\H^0(C,\ms L^d)$ which have no common zeros such that the resulting degree $d$ finite morphism $g:C\to\mathbf{P}^1$ is separable and ramifies only over a single $K$-point of $\mathbf{P}^1$, which we may choose to be not equal to either $0$ or $\infty$. We will now modify this pencil so that the desired properties hold. Our argument to ensure property~\eqref{item:pencil 3} holds is itself inspired by Kedlaya's proof.
    
    Let $\Sigma\subset C$ be the ramification divisor of $g$, let $Z_0=V(s_1)$ be the zero divisor of $g$, and let $Z_{\infty}=V(s_0)$ be the pole divisor of $g$. As $g$ is unramified over $0$ and $\infty$, the three divisors $\Sigma,Z_0,$ and $Z_{\infty}$ are pairwise disjoint. Choose distinct closed points $P,Q\in C$ which are not contained in $\Sigma,Z_0,$ or $Z_{\infty}$ such that $P$ is not defined over any proper subfield of $K$ (Lemma \ref{lem:field of definition lemma, finite field case}). By Riemann--Roch, for all sufficiently large integers $e$ we may find sections $t_0,t_1\in\H^0(C,\ms L^e)$ which have no common zeros such that $t_0$ vanishes on each point of $\Sigma+P+Q$ and does not vanish on any point of $Z_{\infty}$, the orders of vanishing of $t_0$ at $P$ and $Q$ are each $\geqslant 2$, and such that the order of vanishing of $t_0$ at $P$ is strictly greater than the order of vanishing of $t_0$ at any other point. By Dirichlet's theorem on primes in arithmetic progressions, we may assume that $d+pe$ is prime and that $d+pe\geqslant N$. Consider the sections $u_0,u_1\in\H^0(C,\ms L^{d+pe})$ defined by 
    \[
        u_0=s_0t_0^{p}\hspace{1cm}\mbox{and}\hspace{1cm} u_1=s_1t_0^p+s_0t_1^p.
    \]
    As $t_0,t_1$ and $t_0,s_0$ have no common zeros, $u_0$ and $u_1$ have no common zeros. Therefore the sections $u_0,u_1$ define a finite morphism say $f:C\to\mathbf{P}^1$. We claim that $f$ has the desired properties.

    As $g$ was separable, so is $f$. Furthermore, $f$ has degree $d+pe$, which by construction is prime and is $\geqslant N$, so~\eqref{item:pencil 1} holds. As $s_0$ does not vanish at $P$ or at $Q$, the orders of vanishing of $u_0$ at $P$ and at $Q$ are equal to $p$ times the orders of vanishing of $t_0$ at $P$ and at $Q$, which by assumption are different and each $\geqslant 2$. Therefore~\eqref{item:pencil 2} holds. To verify~\eqref{item:pencil 3}, let $R\in C$ be a point at which $f$ ramifies. Assume for the sake of contradiction that $u_0$ does not vanish at $R$. We may write the rational function $u_1/u_0$ as
    \[
        \dfrac{u_1}{u_0}=\dfrac{s_1}{s_0}+\left(\dfrac{t_1}{t_0}\right)^p.
    \]
    Using that the differential of a $p$th power is zero, we compute that
    \[
        \mathrm{d}\left(\dfrac{u_1}{u_0}\right)=\mathrm{d}\left(\dfrac{s_1}{s_0}\right).
    \]
    It follows that $\mathrm{d}(s_1/s_0)$ vanishes at $R$, and so $R\in\Sigma$. But $t_0$ vanishes at every point of $\Sigma$, so $t_0$ and hence $u_0$ vanish at $R$, which is a contradiction. We conclude that $u_0$ vanishes at $R$, so $f(R)=\infty$, and hence~\eqref{item:pencil 3} holds. Finally, by construction there is no other point of $C$ at which $f$ ramifies to the same degree as it does at $P$. Therefore if $f$ is defined over a subfield $L\subset K$ then $P$ is also defined over $L$. By our choice of $P$, we have $K=L$, so~\eqref{item:pencil 4} holds.
\end{proof}

The proof of the preceding proposition used the following lemma.

\begin{Lemma}\label{lem:field of definition lemma, finite field case}
    Let $C$ be a curve over a finite field $K$. There exist infinitely many closed points $P\in C$ each of which is not defined over any proper subfield of $K$ (Definition \ref{def:field of definition of some things}).
\end{Lemma}
\begin{proof}
    For a curve $D$ over a finite field $L$ with and a positive integer $n$, we let $N_n(D/L)$ denote the number of closed points of $D$ which have degree $n$ over $L$. Set $q=|L|$. We have the Hasse--Weil bound
    \[
        |N_n(D/L)-(q^n+1)|\leqslant 2g_Dq^{n/2}
    \]
    which implies the asymptotic expression
    \[
        N_n(D/L)=q^n+O(q^{n/2}).
    \]
    Say $K/L$ is a field extension with $K\neq L$ and $C_L$ is a curve over $L$ such that $C\cong C_L\otimes_LK$. The above asymptotics imply that $N_n(C/K)$ grows faster than any constant times $N_n(C_L/L)$, and in particular for all sufficiently large $n$ there is a point of $C$ of degree $n$ over $K$ which is not the base change of a point of $C_L$. There are only finitely possibilities for $L$ and for $C_L$, so we may even find a point of which is not the base change of a closed point of any descent of $C$ to a proper subfield of $K$.
\end{proof}

The following two results explain why the properties of Proposition \ref{prop:pencils with auts over finite fields} will be useful to us.

\begin{Lemma}\label{lem:transitive action?}
    Let $f:C\to D$ be a finite separable morphism of curves over a field $K$. Suppose that $f$ has prime degree, and that there exists a closed point of $D$ whose preimage under $f$ contains two closed points at which $f$ is ramified to different degrees. If $\widetilde{C}\to C$ is an unramified Galois cover, then the natural map $\Aut_K(\widetilde{C}/C)\to\Aut_K(\widetilde{C}/D)$ is an isomorphism. In particular, $\Aut_K(C/D)=1$.
\end{Lemma}
\begin{proof}
    The automorphism group scheme of a Galois cover acts transitively on each fiber. Thus, a finite morphism of curves satisfying the given condition on ramification cannot be Galois. In particular, $f:C\to D$ is not Galois. Furthermore, as $\widetilde{C}\to C$ is unramified, the composition $\widetilde{C}\to C\to D$ also satisfies the same condition on ramification, and hence $\widetilde{C}\to D$ is not Galois. 

    Now, set $E=\widetilde{C}/\Aut_K(\widetilde{C}/D)$. As $\widetilde{C}\to C$ is Galois, we have $C=\widetilde{C}/\Aut_K(\widetilde{C}/C)$, and therefore we have a factorization
    \[
        C\to E\to D
    \]
    of $f$. We assume that $f$ has prime degree, so either $C=E$ or $E=D$. The latter cannot be the case, as then the composition $\widetilde{C}\to C\to D$ would be Galois. We conclude that $C=E$, and therefore $\Aut_K(\widetilde{C}/C)=\Aut_K(\widetilde{C}/D)$.
\end{proof}

We now work over a gerbe. Let $\ms G$ be a Deligne--Mumford gerbe over a field $K$ with inertia $G_{\ms G}:=\ms I_{\ms G}$. Consider a diagram
\[
    \begin{tikzcd}
        C\arrow{d}[swap]{\pi}\arrow{r}{f}&D\\
        \ms G&
    \end{tikzcd}
\]
in which $f:C\to D$ is a finite separable morphism of curves over $K$ and $\pi:C\to\ms G$ is a morphism which is a relative curve. 

\begin{Lemma}\label{lem:pencils over finite fields}
    Suppose that $f$ has prime degree and that there exists a closed point of $D$ whose preimage under $f$ contains two points at which $f$ ramifies to different degrees. The inertial action on $C$ induces an isomorphism $G_{\ms G}\iso\sAut_{\ms G}(C/D)$, and also we have $\sAut_K(C/D)=1$.
\end{Lemma}
\begin{proof}
    By Lemma \ref{lem:transitive action?}, we have $\sAut_K(C/D)=1$. Let $L/K$ be a separable closure of $K$ and let $s\in\ms G(L)$ be a section. The pullback of $G_{\ms G}$ along $s$ is then the $L$-group scheme associated to a finite group, say $G$. Let $\widetilde{C}_L$ be the pullback of $C$ along $s$. As in Example \ref{ex:non split gerbe, once again}, the pullback of $f$ along $s$ factors as the composition
    \[
        \widetilde{C}_L\xrightarrow{|G|} C_L\xrightarrow{f_L} D_L
    \]
    where the left arrow is a torsor under $G$, ie. an unramified Galois cover. The pullback of the inertial action map $G_{\ms G}\to\sAut_{\ms G}(C/D)$ under $s$ is the map
    \[
         G\iso\sAut_{K}(\widetilde{C}_L/C_L)\to\sAut_K(\widetilde{C}_L/D_L).
    \]
    By Lemma \ref{lem:transitive action?}, this map is an isomorphism on $L$-points, hence an isomorphism.
\end{proof}

We now prove the main result of this section, which is the finite field case of Theorem \ref{thm:finite cover of curve over a gerbe}. We recall the notation: $K/K_0$ is a finitely generated field extension, $\ms G$ is a Deligne--Mumford gerbe over $K$ with inertia $G_{\ms G}:=\ms I_{\ms G}$, $C$ is a curve over $K$, and $\pi:C\to\ms G$ is a morphism which is a relative curve.

\begin{Theorem}\label{thm:finite cover of curve over a gerbe, field of definition, finite field case}
    Theorem \ref{thm:finite cover of curve over a gerbe} holds if $K$ is finite.
\end{Theorem}
\begin{proof}
    Our construction is similar to that employed in the proofs of Theorems \ref{thm:finite cover of curve over a gerbe, infinite field case} and \ref{thm:finite cover of curve over a gerbe, field of definition, infinite field case}, but requires assembling the pieces in a slightly different manner. We apply Proposition \ref{prop:pencils with auts over finite fields} to $C$ to obtain a finite morphism $f:C\to\mathbf{P}^1$ whose degree is prime and which satisfies conditions~\eqref{item:pencil 2},~\eqref{item:pencil 3}, and~\eqref{item:pencil 4} of Proposition \ref{prop:pencils with auts over finite fields}. By Lemma \ref{lem:pencils over finite fields}, we then have that the inertial action induces an isomorphism $G_{\ms G}\iso\sAut_{\ms G}(C/\mathbf{P}^1)$, and furthermore we have $\sAut_K(C/\mathbf{P}^1)=1$. Let $X$ be a curve over $K$ such that $\sAut_K(X)=1$ (Theorem \ref{thm:curves with no auts}). We apply Proposition \ref{prop:pencils with auts over finite fields} to $X$ to obtain a finite morphism $g:X\to\mathbf{P}^1$ whose degree is a prime which is strictly greater than $|G_{\ms G}|\deg(f)$ and which is ramified only over $\infty$ (note that, unlike in the proof of Theorem \ref{thm:finite cover of curve over a gerbe, infinite field case}, this morphism need not be \textit{totally} ramified over $\infty$). We may assume furthermore that $X$ and $g$ are defined over $K_0$. Let $\lambda$ be a positive integer, and let $\varphi:\mathbf{P}^1\to\mathbf{P}^1$ be a morphism which satisfies the conclusions of Lemma \ref{lem:incompressible lemma} applied to $f$ and $\lambda$. As before, we define $D$ and $E$ to be the fiber products in the diagram
    \[
        \begin{tikzcd}
            E\arrow{d}\arrow{r}{f''}&X\arrow{d}{g}\\
            D\arrow{d}\arrow{r}{f'}&\mathbf{P}^1\arrow{d}{\varphi}\\
            C\arrow{r}{f}\arrow{d}[swap]{\pi}&\mathbf{P}^1\\
            \ms G&
        \end{tikzcd}
    \]
    with Cartesian squares. The branch loci of $f$ and $\varphi$ are disjoint, and $\varphi$ is totally ramified over a $K$-point of $\mathbf{P}^1$. Therefore condition~\eqref{item:fiber product condition 1 gerbe} of Lemma \ref{lem:relative curve lemma 1} holds for the lower square, and so $D$ is a smooth proper geometrically integral curve over $K$, the morphism $D\to\ms G$ is a relative curve, and the lower square in the above diagram satisfies condition $(\ast\ast)$. As $g$ ramifies only over $\infty$, the branch loci of $f'$ and $g$ are disjoint, and by construction we have that $|G_{\ms G}|\deg(f)<\deg(g)$ and $\deg(g)$ is prime. Therefore condition~\eqref{item:fiber product condition 2 gerbe} of Lemma \ref{lem:relative curve lemma 1} holds for the upper square, so $E$ is a smooth proper geometrically integral curve over $K$, $E\to\ms G$ is a relative curve, and the upper square satisfies condition $(\ast\ast)$. We claim that if $\lambda$ is sufficiently large then the morphism $E\to C$ has the desired properties. We note that, as in the proof of Theorem \ref{thm:finite cover of curve over a gerbe, infinite field case}, $f'$ is $\lambda$-incompressible, so we have $g_E\geqslant g_D\geqslant\lambda$, and therefore we may assume that the genus of $E$ is arbitrarily large.
    
    The proofs of Claims \ref{claim:claim 1}, \ref{claim:claim 2}, and \ref{claim:claim 3} of the proof of Theorem \ref{thm:finite cover of curve over a gerbe, infinite field case} apply without change in our situation. (We remark that the proofs of these claims use the fact that the two squares in the above diagram each satisfy condition $(\ast\ast)$. While we have verified above that this does indeed hold in our case as well, this is for different reasons than in the infinite field case.) We conclude that the inertial action map $G_{\ms G}\iso\sAut_{\ms G}(E)$ is an isomorphism and that $\sAut_K(E)=1$. It remains to show that $E$ is not defined over any proper intermediate subfield. The proof of this is exactly the same as in the proof of Theorem \ref{thm:finite cover of curve over a gerbe, field of definition, infinite field case}: if $E$ is defined over a subfield $L\subset K$, then by Proposition \ref{prop:step 1} $f''$ is also defined over $L$, and by Proposition \ref{prop:step 2} $f$ is also defined over $L$. By our choice of $f$, we conclude that $K=L$.
\end{proof}

\begin{Remark}\label{rem:gerbes over finite fields are split}
    In fact, every Deligne--Mumford gerbe over a finite field is split. This follows from a theorem of Grothendieck on the vanishing of nonabelian $\H^2$ over fields of cohomological dimension $\leqslant 1$ \cite[Theorem 3.5]{SPRING}.
\end{Remark}

\section{Proofs of main results}\label{sec:checking murphys law}

In this section we give the proofs of the results stated in \S\ref{sec:introduction}. We first use our constructions in \S\ref{sec:curves over gerbes} to verify that Murphy's Law holds for the moduli stack of curves (Theorem \ref{thm:stacky murphy's law for curves}). Using geometric constructions starting from curves, we then prove that Murphy's Law holds for the stack of principally polarized abelian varieties (Theorem \ref{thm:stacky murphys law for pp abelian varieties}) and the stack of canonically polarized varieties (Theorem \ref{thm:stacky murphy's law for canonically polarized varieties}), thus completing the proof of Theorem \ref{thm:stacky murphy's law}. Finally, in \S\ref{ssec:proofs of further results}, we prove the results stated in \S\ref{ssec:further results}.

\subsection{Murphy's Law for the stack of curves}

In this section we will show that the stack of curves over a field satisfies Stacky Murphy's Law. Let $K_0$ be a field. Let 
\[
    \ms M_{\bullet}=\ms M_{0,3}\sqcup\ms M_{1,1}\sqcup\ms M_2\sqcup\ms M_3\sqcup\dots
\]
be the stack of smooth proper geometrically integral curves over $K_0$. Let $K/K_0$ be a finitely generated field extension and let $\ms G$ be a Deligne--Mumford gerbe over $K$ with inertia $G_{\ms G}:=\ms I_{\ms G}$.

\begin{Proposition}\label{prop:gerbe --> stack}
    Let $C$ be a curve over $K$ equipped with a $K$-morphism $\pi:C\to\ms G$ which is a relative curve of genus $g\geqslant 2$. Suppose that
    \begin{enumerate}
        \item\label{item:gerbe stack prop 1} the inertial action induces an isomorphism $G_{\ms G}\iso\sAut_{\ms G}(C)$, and
        \item\label{item:gerbe stack prop 2} $C$ is not defined over any intermediate extension $K/L/K_0$ such that $K\neq L$.
    \end{enumerate}
    The morphism $m_{\pi}:\ms G\to\ms M_{g}$ induced by $\pi$ induces an isomorphism between $\ms G$ and the residual gerbe of a point of $\ms M_{g}$.
\end{Proposition}
\begin{proof}
    Let $x\in|\ms M_{g}|$ be the image of the unique element of $|\ms G|$ under the map $|m_{\pi}|:|\ms G|\to|\ms M_g|$. Let $\ms G(x)$ be the residual gerbe and let $k(x)$ be the residue field of $\ms M_g$ at $x$. Let $\ms C_g\to\ms M_{g}$ be the universal family of curves of genus $g$ and let $\ms C(x)$ be the restriction of $\ms C_g$ to $\ms G(x)$. We have a diagram
    \begin{equation}\label{eq:two squares of gerbes and curves}
        \begin{tikzcd}
            C\arrow{d}[swap]{\pi}\arrow{r}&\ms C(x)\arrow{d}\arrow[hook]{r}&\ms C_g\arrow{d}\\
            \ms G\arrow{r}{m_{\pi}}\arrow{d}&\ms G(x)\arrow{d}\arrow[hook]{r}&\ms M_g\\
            \Spec K\arrow{r}&\Spec k(x)&
        \end{tikzcd}
    \end{equation}
    in which the top two squares are Cartesian. We have an identification $\ms I_{\ms M_g}=\sAut_{\ms M_g}(\ms C_g)$. The inclusion of the residual gerbe is stabilizer preserving (see Definition \ref{def:stabilizer preserving} and Example \ref{ex:inclusion of residual gerbe is stabilizer preserving}), hence induces an isomorphism $\ms I_{\ms G(x)}\iso\sAut_{\ms G(x)}(\ms C(x))$. With this identification the commutative square of inertia stacks~\eqref{eq:induced square on inertia} induced by $m_{\pi}$ becomes a commutative square
    \[
        \begin{tikzcd}
            G_{\ms G}\arrow{d}\arrow{r}&\sAut_{\ms G(x)}(\ms C(x))\arrow{d}\\
            \ms G\arrow{r}{m_{\pi}}&\ms G(x).
        \end{tikzcd}
    \]
    The induced map
    \[
        G_{\ms G}\to m_{\pi}^{-1}\sAut_{\ms G(x)}(\ms C(x))=\sAut_{\ms G}(C)
    \]
    of sheaves of groups on $\ms G$ is the inertial action map, which by assumption~\eqref{item:gerbe stack prop 1} is an isomorphism. Thus $m_{\pi}$ is stabilizer preserving, and therefore Lemma \ref{lem:morphism of gerbes} implies that the lower square of~\eqref{eq:two squares of gerbes and curves} is also Cartesian. It follows that $\ms C(x)$ is itself a curve over $\Spec k(x)$ (rather than merely a stacky curve) and that $C$ is defined over $k(x)$. From assumption~\eqref{item:gerbe stack prop 2} we deduce that the map $\Spec K\to\Spec k(x)$ is an isomorphism, and hence the map $m_{\pi}:\ms G\to\ms G(x)$ is an isomorphism.
\end{proof}

\begin{Theorem}\label{thm:stacky murphy's law for curves}
    The stack $\ms M_{\bullet}$ satisfies Stacky Murphy's Law.
\end{Theorem}
\begin{proof}
   Theorem \ref{thm:main theorem for curves over gerbes} and Proposition \ref{prop:gerbe --> stack} combined imply that if $K/K_0$ is a finitely generated field extension and $\ms G$ is a Deligne--Mumford gerbe over $K$, then $\ms G$ is isomorphic to the residual gerbe of a point of $\ms M_{g}$ for some $g$.
\end{proof}

\begin{Remark}
    The reader might be confused by the fact that we obtain nontrivial gerbes over $K$ using a curve $C$ whose field of definition is $K$. Note however that we are considering the point of $\ms M_{g}$ obtained from the family $\pi:C\to\ms G$, not from the $K$-curve $C\to\Spec K$. These need not be the same. Indeed, if $s\in\ms G(L)$ is a splitting of $\ms G$ over a field extension $L/K$, then setting $\widetilde{C}_L:=C\times_{\ms G}\Spec L$ and $C_L:=C\times_{\Spec K}\Spec L$ (as in Example \ref{ex:non split gerbe}), we have a finite morphism $\widetilde{C}_L\to C_L$ of curves over $L$ which is a torsor for the $L$-group scheme $ G_L:=s^{-1}G_{\ms G}$. Thus, the fiber of $C\to\ms G$ is a finite cover of the fiber of $C\to\Spec K$. Therefore, writing $g'$ for the genus of the curve $C$, we have that $g'\geqslant g$, and typically this inequality will be strict. In contrast to the point $x\in |\ms M_{g}|$ corresponding to $C\to\ms G$ considered above, the point $x'\in |\ms M_{g'}|$ corresponding to $C\to\Spec K$ has residue field $K$, and the residual gerbe at $x'$ is the split gerbe $\B\sAut_K(C)$ over $K$.
\end{Remark}

\subsection{Transfer of Murphy's Law along a morphism of stacks}

We record two results which relate the validity of Stacky Murphy's Law for stacks related by a suitable morphism.

\begin{Proposition}\label{prop:murphy and an immersion}
    Let $f:\ms M\hookrightarrow\ms N$ be a monomorphism of Deligne--Mumford stacks over a field $K_0$. If $\ms M$ satisfies Stacky Murphy's Law, then so does $\ms N$.
\end{Proposition}
\begin{proof}
    As $f$ is a monomorphism, Lemma \ref{lem:immersions preserve residual gerbes} implies that $f$ induces an isomorphism $\ms G(x)\iso\ms G(f(x))$ of residual gerbes for any point $x\in|\ms M|$.
\end{proof}

\begin{Proposition}
    Let $q:\ms M\to\ms N$ be a morphism of Deligne--Mumford stacks over a field $K_0$. Suppose that $q$ is representable and that there exists a Zariski open cover of $\ms N$ over which $q$ admits a section. If $\ms N$ satisfies Stacky Murphy's Law, then so does $\ms M$.
\end{Proposition}
\begin{proof}
    Let $\ms U\to\ms N$ be a Zariski open cover over which $q$ admits a section, say $\sigma:\ms U\to\ms M$. Say $\ms U$ is the disjoint union of stacks $\ms U_i$ where the maps $\ms U_i\to\ms N$ are open immersions. By Lemma \ref{lem:immersions preserve residual gerbes}, each residual gerbe of $\ms N$ is isomorphic to a residual gerbe of one of the $\ms U_i$. The maps $\ms U_i\to\ms M$ are immersions, so by the same result, $\ms M$ satisfies Stacky Murphy's Law.
\end{proof}

\subsection{Murphy's Law for the stack of principally polarized abelian varieties}

Let $K_0$ be a field. Let $\ms A_{\bullet}$ denote the stack of principally polarized abelian varieties over $K_0$. In this section we will prove that $\ms A_{\bullet}$ satisfies Stacky Murphy's Law.

Let $\ms M_{\bullet}$ be the stack of curves over $K_0$ and let $\ms H_{\bullet}\subset\ms M_{\bullet}$ denote the hyperelliptic locus (a closed substack). Let $\ms M_{\bullet}^{\mathrm{free}}\subset\ms M_{\bullet}$ denote the open substack parameterizing curves whose automorphism group acts freely (without fixed points). Equivalently, $\ms M^{\mathrm{free}}_{\bullet}$ is the complement in $\ms M_{\bullet}$ of the image of the stacky locus of the universal curve.

\begin{Theorem}\label{thm:stacky murphys law for pp abelian varieties}
     The stack $\ms M^{\mathrm{free}}_{\bullet}$, the stack $\ms M_{\bullet}\setminus\ms H_{\bullet}$ of non hyperelliptic curves, and the stack $\ms A_{\bullet}$ of principally polarized abelian varieties all satisfy Stacky Murphy's Law.
\end{Theorem}
\begin{proof}
    By Theorem \ref{thm:stacky murphy's law for curves}, every Deligne--Mumford gerbe over a finitely generated extension of $K_0$ appears as a residual gerbe of some point $x\in|\ms M_{\bullet}|$. In fact, the proof shows something stronger: the curves over gerbes produced in Theorem \ref{thm:main theorem for curves over gerbes} have free inertial action, so we may choose the point $x$ to lie in the open substack $\ms M^{\mathrm{free}}_{\bullet}\subset\ms M_{\bullet}$. By Lemma \ref{lem:immersions preserve residual gerbes} the stack $\ms M^{\mathrm{free}}_{\bullet}$ therefore also satisfies Stacky Murphy's Law. The automorphism group of a hyperelliptic curve never acts freely, as the hyperelliptic involution necessarily fixes points. Thus $\ms M^{\mathrm{free}}_{\bullet}$ is contained in the non-hyperelliptic locus, so the stack $\ms M_{\bullet}\setminus\ms H_{\bullet}$ also satisfies Stacky Murphy's Law. Consider the Torelli map $\tau:\ms M_{\bullet}\to \ms A_{\bullet}$. It follows from classical results of Oort and Steenbrink \cite[Theorem 2.6, 2.7]{MR0605341} that $\tau$ restricts to an immersion $\ms M_{\bullet}\setminus\ms H_{\bullet}\to \ms A_{\bullet}$. Applying Proposition \ref{prop:murphy and an immersion}, we conclude that $\ms A_{\bullet}$ satisfies Stacky Murphy's Law.
\end{proof}

\subsection{Murphy's Law for the stack of canonically polarized varieties}

In this section we will show by taking products of curves that the moduli stack of canonically polarized varieties of any fixed positive dimension satisfies Stacky Murphy's Law. We begin with the following result concerning isomorphisms of products of curves. Let $K$ be a field and let $d\geqslant 1$ be an integer.
 
\begin{Lemma}\label{lem:isomorphisms of products of curves}
    Let $C_1,\dots,C_d$ and $D_1,\dots,D_d$ be two collections of pairwise non-isomorphic curves of genus $\geqslant 2$ over $K$. If
    \[
        \alpha:C_1\times\dots\times C_d\iso D_1\times\dots\times D_d
    \]
    is an isomorphism over $K$, then there exists a permutation $\sigma$ of $\left\{1,\dots,d\right\}$ and isomorphisms $\beta_i:C_{\sigma(i)}\iso D_{i}$ such that $\alpha$ is equal to the composition
    \[
        C_1\times\dots\times C_d\xrightarrow{\widetilde{\sigma}}C_{\sigma(1)}\times\dots\times C_{\sigma(d)}\xrightarrow{\beta_1\times\dots\times\beta_d}D_1\times\dots\times D_d,
    \]
    where $\widetilde{\sigma}$ is the isomorphism which permutes the factors according to $\sigma$. Moreover, the permutation $\sigma$ and the isomorphisms $\beta_i$ are uniquely determined by $\alpha$.
\end{Lemma}
\begin{proof}
    For $1\leqslant i\leqslant d$ we let $\alpha_i:C_1\times\dots\times C_d\to D_i$ be the composition of $\alpha$ with the $i$th projection map $\pi_i:D_1\times\dots\times D_d\to D_i$. As $\alpha$ is an isomorphism, each $\alpha_i$ must be nonconstant on some factor. For each $i$ we choose such a factor, say $C_{\sigma(i)}$, and we let $\sigma:\left\{1,\dots,d\right\}\to\left\{1,\dots,d\right\}$ be the resulting function. We consider the diagrams
    \[
        \begin{tikzcd}
            C_1\times\dots\times C_d\arrow{r}{\alpha_i}\arrow{d}{\pi^{\sigma(i)}}&D_i\\
            C_1\times\dots\times\widehat{C}_{\sigma(i)}\times\dots\times C_d&
        \end{tikzcd}
    \]
    where the hat indicates that we omit the $C_{\sigma(i)}$-factor. The map $\alpha_i$ thus corresponds to a family of morphisms $C_{\sigma(i)}\to D_i$ parameterized by the lower product. We assume that the $D_i$ have genus $\geqslant 2$, so for each $i$ the Hom-scheme of surjective morphisms $C_{\sigma(i)}\to D_i$ is discrete. Therefore $\alpha_i$ factors as the projection onto $C_{\sigma(i)}$ followed by a surjective morphism $\beta_i:C_{\sigma(i)}\to D_i$. As $\alpha$ is an isomorphism, $C_{\sigma(i)}$ is the only factor on which $\alpha_i$ is nonconstant, and furthermore the function $\sigma:\left\{1,\dots,d\right\}\to\left\{1,\dots,d\right\}$ must be a bijection and each of the $\beta_i$ must be an isomorphism. This shows that $\alpha$ factors as claimed. Finally, as we assumed the $C_i$ were pairwise non-isomorphic, the permutation $\sigma$ and the isomorphisms $\beta_i$ are uniquely determined by $\alpha$.
\end{proof}

\begin{Corollary}\label{cor:automorphism group of product of curves}
    Let $C_1,\dots,C_d$ be pairwise non-isomorphic curves of genus $\geqslant 2$ over $K$. The map
    \[
        \Aut_{K}(C_1)\times\dots\times\Aut_{K}(C_d)\to\Aut_{K}(C_1\times\dots\times C_d)
    \]
    defined by $(\beta_1,\dots,\beta_d)\mapsto \beta_1\times\dots\times\beta_d$ is an isomorphism.
\end{Corollary}
\begin{proof}
    This map is certainly injective. To see that it is surjective, let $\alpha$ be an automorphism of $C_1\times\dots\times C_d$. Applying Lemma \ref{lem:isomorphisms of products of curves} we find a permutation $\sigma$ of the set $\left\{1,\dots,d\right\}$ and isomorphisms $\beta_i:C_{\sigma(i)}\iso C_i$ such that $\alpha=(\beta_1\times\dots\beta_d)\circ\widetilde{\sigma}$. We assumed that the $C_i$ were pairwise non-isomorphic, so in fact $\sigma(i)=i$ for all $i$. Therefore each of the $\beta_i$ is an automorphism of $C_i$ and we have $\alpha=\beta_1\times\dots\times\beta_d$, so $\alpha$ is indeed in the image of the map.
\end{proof}

Let $K_0$ be a field. Fix an integer $d\geqslant 1$ and let $\ms M_{\omega}^d$ be the stack over $K_0$ parameterizing smooth proper geometrically integral canonically polarized varieties of dimension $d$ with reduced automorphism group scheme. Let $\ms M_{\geqslant 2}$ be the stack of curves of genus $\geqslant 2$ over $K_0$. Let 
\[
    \Delta\subset\underbrace{\ms M_{\geqslant 2}\times\dots\times\ms M_{\geqslant 2}}_{d}
\]
be the big diagonal in the $d$-fold self-product of $\ms M_{\geqslant 2}$ and let $\ms U^d$ be its complement.

\begin{Proposition}\label{prop:U satisfies the law}
    The stack $\ms U^d$ satisfies Stacky Murphy's Law.
\end{Proposition}
\begin{proof}
    Let $K/K_0$ be a finitely generated field extension and let $\ms G$ be a Deligne--Mumford gerbe over $K$. By Theorem \ref{thm:stacky murphy's law for curves}, $\ms G$ is isomorphic to a residual gerbe of a point of $\ms M_{\geqslant 2}$, say via a map $x_1:\ms G\to\ms M_{\geqslant 2}$. The same is true for the trivial $K$-gerbe $\Spec K$. Moreover, by Theorem \ref{thm:main theorem for curves over gerbes}, the gerbe $\Spec K$ appears as a residual gerbe in $\ms M_g$ for arbitrarily large $g$. Thus, we may find genera $2\leqslant g_2<\dots<g_d$ and maps $x_i:\Spec K\to\ms M_{g_i}$  for $i=2,\ldots,d$ each of which realizes $\Spec K$ as a residual gerbe of $\ms M_{g_i}$. The morphism
    \[
        \ms G\cong\ms G\times_{K}\Spec K\times_{K}\dots\times_K\Spec K\xrightarrow{x_1\times x_2\times\dots\times x_d}\ms M_{\geqslant 2}\times_K\dots\times_K\ms M_{\geqslant 2}
    \]
    factors through $\ms U^d$, and realizes $\ms G$ as a residual gerbe of $\ms U^d$.
\end{proof}

Consider the morphism
\begin{equation}\label{eq:product morphism}
    \Pi:\underbrace{\ms M_{\geqslant 2}\times\dots\times\ms M_{\geqslant 2}}_{d}\to\ms M_{\omega}^d
\end{equation}
which sends a tuple $(C_1,\dots,C_d)$ to the product $C_1\times\dots\times C_d$.

\begin{Proposition}\label{prop:product of curves morphism}
    The morphism $\Pi$ restricts to an open immersion $\Pi_{\ms U^d}:\ms U^d\to\ms M_{\omega}^d$.
\end{Proposition}
\begin{proof}    
    We will show that $\Pi_{\ms U^d}$ is an \'{e}tale monomorphism. To show that $\Pi_{\ms U^d}$ is \'{e}tale it will suffice to verify that it is formally \'{e}tale. In the case $d=2$, this follows from \cite[Theorem 3.3]{MR3143705}, and the general case follows by induction. It remains to show that $\Pi_{\ms U^d}$ is a monomorphism. We will verify the conditions of Lemma \ref{lem:monomorphism lemma}. We have already shown that $\Pi_{\ms U^d}$ is \'{e}tale, and in particular unramified, so condition~\eqref{item:unramified} holds. Condition~\eqref{item:injective on pts} follows from Lemma \ref{lem:isomorphisms of products of curves}, and condition~\eqref{item:inertia is inert} follows from Corollary \ref{cor:automorphism group of product of curves}. This completes the proof.
\end{proof}

\begin{Theorem}\label{thm:stacky murphy's law for canonically polarized varieties}
    For each integer $d\geqslant 1$ the stack $\ms M_{\omega}^d$ satisfies Stacky Murphy's Law.
\end{Theorem}
\begin{proof}
    By Proposition \ref{prop:U satisfies the law} the stack $\ms U^d$ satisfies Stacky Murphy's Law. By Proposition \ref{prop:product of curves morphism}, the map $\Pi_{\ms U^d}:\ms U^d\to\ms M^d_{\omega}$ is an open immersion, so by Proposition \ref{prop:murphy and an immersion} the stack $\ms M_{\omega}^d$ also satisfies Stacky Murphy's Law.
\end{proof}

\subsection{Proofs of further results}\label{ssec:proofs of further results}

In this section we give the proofs of the results stated in \S\ref{ssec:further results}. This will complete the proofs of all results stated in \S\ref{sec:introduction}.

\begin{proof}[Proof of Theorem \ref{thm:Murphy's law for BG}]
    Suppose that $\ms M$ is a DM stack over $K_0$ satisfying Stacky Murphy's Law. Then for any finitely generated field extension $K/K_0$ and finite \'{e}tale group scheme $G$ over $K$, we may find a morphism $\iota:\B G\hookrightarrow\ms M$ which induces an isomorphism between $\B G$ and the residual gerbe of a point of $\ms M$. Let $s:\Spec K\to\B G$ be the canonical splitting and write $x_K:=\iota\circ s$, so that we have a diagram
    \[
        \begin{tikzcd}
            \Spec K\arrow[bend left=25]{rr}{x_K}\arrow[equals]{dr}\arrow{r}{s}&\B G\arrow{d}\arrow[hook]{r}{\iota}&\ms M\\
            &\Spec K.&
        \end{tikzcd}
    \]
    The inclusion of the residual gerbe is stabilizer preserving, hence $\iota$ induces an isomorphism $G_{\B G}\iso\iota^{-1}\ms I_{\ms M}$. The pullback under $s$ of this map gives an isomorphism $\sAut_K(s)\iso\sAut_K(x_K)$. Composing with the canonical identification $G=\sAut_K(s)$, we obtain an isomorphism $G\iso\sAut_K(x_K)$. Suppose that $K/L/K_0$ is an intermediate extension over which $x$ is defined. Then there exists a morphism $\Spec L\to\ms M$ representing $x$. By the universal property of the residual gerbe, this map factors through $\iota$, and hence we have that $K\subset L$, so $L=K$.
\end{proof}

\begin{proof}[Proof of Corollary \ref{cor:presecribed automorphisms of curves 1}]
    Let $K/K_0$ be a finitely generated field extension and let $G$ be a finite \'{e}tale group scheme over $K$. Choose a curve $C$ over $K$ and a $K$-morphism $\pi:C\to\B G$ which satisfy the conclusions of Theorem \ref{thm:main theorem for curves over gerbes}. Let $s:\Spec K\to\B G$ be the canonical splitting and define $\widetilde{C}$ by the Cartesian diagram
    \[
        \begin{tikzcd}
            \widetilde{C}\arrow{r}\arrow{d}&C\arrow{d}{\pi}\\
            \Spec K\arrow{r}{s}&\B G.
        \end{tikzcd}
    \]
    As $\pi$ is a relative curve, $\widetilde{C}$ is a curve over $K$. By assumption, we have that $C$ is not defined over any proper intermediate extension of $K/K_0$, and that $\sAut_K(C)=1$. Furthermore, the inertial action map $G_{\B G}\to\sAut_{\B G}(C)$ is an isomorphism. The pullback under $s$ of this map gives an isomorphism $G\iso\sAut_{K}(\widetilde{C})$, and we have $\widetilde{C}/G\iso C$. It remains to show that $\widetilde{C}$ is not defined over any proper intermediate extension of $K/K_0$. Let $x\in|\ms M|$ be the point corresponding to $\widetilde{C}$. Then the residual gerbe at $x$ is isomorphic to $\B G$ and the residue field at $x$ is isomorphic to $K$. Thus, if $\widetilde{C}$ is defined over an intermediate extension $K/L/K_0$, then by the universal property of the residual gerbe we obtain that $K\subset L$, so $K=L$.
\end{proof}

\begin{proof}[Proof of Theorem \ref{thm:residue fields of M_g}]
    Suppose that $\ms M$ is a DM stack over $K_0$ with finite diagonal which satisfies Stacky Murphy's Law, and let $\rho:\ms M\to M$ be the coarse moduli space of $\ms M$. Let $K/K_0$ be a finitely generated field extension. Regarding $\Spec K$ as the trivial gerbe over $K$, we may find a point $x\in|\ms M|$ such that the residual gerbe of $\ms M$ at $x$ is isomorphic to $\Spec K$. Then the residue field of $x$ is also isomorphic to $\Spec K$. As the stabilizer of $\ms M$ at $x$ is trivial, $\ms M$ is tame at $x$, so by Proposition \ref{prop:residue field with tame assumption} the residue field of the point $\rho(x)\in|M|$ is isomorphic to $K$. 
\end{proof}

It remains only to prove Theorem \ref{thm:ranks of faithful vector bundles}. We require some preparation first.

\begin{Lemma}\label{lem:representation dimension lemma}
    If $f:\ms N\to\ms M$ is a representable morphism of algebraic stacks, then $\rdim(\ms N)\leqslant\rdim(\ms M)$.
\end{Lemma}
\begin{proof}
    Let $\ms V$ be a faithful vector bundle on $\ms M$. We claim that the pullback $f^*\ms V=\ms V\times_{\ms M}\ms N$ is a faithful vector bundle on $\ms N$. Indeed, the inertial action map for the pullback factors as the composition
    \[
        \ms I_{\ms N}\to \ms I_{\ms M}\times_{\ms M}\ms N\xrightarrow{\alpha_{\ms V}\times\id_{\ms N}}\sAut_{\ms M}(\ms V)\times_{\ms M}\ms N=\sAut_{\ms N}(f^*\ms V)
    \]
    where the left map is the map on inertia induced by $f$ and $\alpha_{\ms V}:\ms I_{\ms M}\to\sAut_{\ms M}(\ms V)$ is the inertial action on $\ms V$. As $f$ is representable, the map $\ms I_{\ms N}\to\ms I_{\ms M}\times_{\ms M}\ms N$ is injective (see Example \ref{ex:stabilizer preserving implies representable}), and by assumption $\alpha_{\ms V}$ is injective. Thus, the composition is injective, so $f^*\ms V$ is faithful.
\end{proof}

\begin{Proposition}\label{prop:period index existence}
    Let $K_0$ be a field and $\ell$ be a prime number not equal to the characteristic of $K_0$. For any integers $a,b$ with $0\leqslant a\leqslant b$, there exists a finitely generated field extension $K/K_0$ and a class $\alpha\in\Br(K)$ with period $\ell^a$ and index $\ell^{b}$.
\end{Proposition}
\begin{proof}
    By enlarging $K_0$ we may assume that $K_0$ contains a primitive $\ell^b$th root of unity, say $\omega$. Consider the field extension $L:=K_0(X,Y)/K_0$ and the cyclic algebra $D=(X,Y)_{\omega}$ over $L$, which has both period and index equal to $\ell^b$ \cite[Lemma 5.5.3]{MR3727161}. Let $K$ be the function field of a Brauer--Severi variety over $L$ corresponding to the algebra $D^{\otimes \ell^{b-a}}$. We claim that the Brauer class $\alpha\in\Br(K)$ of the pullback $D\otimes_LK$ has the desired properties. Indeed, Amitsur's theorem \cite[Theorem 5.4.1]{MR3727161} states that the kernel of the restriction map $\Br(L)\to\Br(K)$ is the cyclic subgroup generated by the class of $D^{\otimes\ell^{b-a}}$, so $D\otimes_LK$ has period $\ell^a$. On the other hand, by Schofield--Van Den Bergh's index reduction results, the index of the pullback $D\otimes_LK$ remains equal to $\ell^b$ \cite[Theorem 2.1]{MR1061778}. 
\end{proof}

\begin{Proposition}\label{prop:gerbes with big index}
    Let $K_0$ be a field and let $G$ be a finite group whose center contains a nontrivial element with order coprime to the characteristic of $K_0$. For any integer $d$, there exists a finitely generated field extension $K/K_0$ and a gerbe $\ms G$ over $\Spec K$ which is locally isomorphic to $\B G$ and satisfies $\rdim(\ms G)\geqslant d$.
\end{Proposition}
\begin{proof}
        We first reduce to the case when $G=\mathbf{Z}/\ell$ for a prime $\ell$ which is not equal to the characteristic of $K_0$. Consider a central extension
        \[
            1\to N\to G\to H\to 1
        \]
        of finite groups. Let $K/K_0$ be a field extension and let $\ms G$ be an $N$-gerbe over $K$. Let $\ms G'$ be the induced $G$-gerbe. There is a canonical morphism $\ms G\to\ms G'$, which is an $H$-torsor, and in particular is representable. By Lemma \ref{lem:representation dimension lemma} we have $\rdim(\ms G)\leqslant\rdim(\ms G')$. It follows that if the claim holds for $N$ then it also holds for $G$. By assumption, $G$ contains a central subgroup of the form $\mathbf{Z}/\ell$ for some prime $\ell$ which is not equal to the characteristic of $K_0$, so we may therefore reduce to this case.

        We now prove the result when $G=\mathbf{Z}/\ell$ for some prime $\ell$ not equal to the characteristic of $K_0$. After possibly extending $K_0$, we may assume that $K_0$ contains a primitive $\ell$th root of unity, and therefore that $\mathbf{Z}/\ell\cong\mu_\ell$. By Proposition \ref{prop:period index existence}, for any integer $m\geqslant 1$ there exists a finitely generated field extension $K/K_0$ and a Brauer class $\alpha\in\Br(K)$ of period $\ell$ and index $\ell^m$. If $\ms G$ is a $\mu_\ell$-gerbe whose cohomology class $[\ms G]\in\H^2(\Spec K,\mu_{\ell})$ maps to $\alpha\in\Br(K)$, then $\rdim(\ms G)=\ell^m$ (see eg. \cite[Proposition 3.1.2.1]{MR2388554}). This completes the proof.
\end{proof}

\begin{Remark}
    Proposition \ref{prop:gerbes with big index} is sharp, in the sense that if $G$ is a finite group whose center has order a power of the exponential characteristic of $K_0$, then if $K/K_0$ is a field extension and $\ms G$ is a gerbe over $K$ which is locally isomorphic to $\B G$ there is a universal bound on $\rdim_K\ms G$ depending only on $G$.
\end{Remark}

\begin{proof}[Proof of Theorem \ref{thm:ranks of faithful vector bundles}]
   By Proposition \ref{prop:gerbes with big index}, for any integer $d$ we may find a finitely generated field extension $K/K_0$ and a Deligne--Mumford gerbe $\ms G$ over $K$ which is locally isomorphic to $\B G$ and has the property that $\rdim(\ms G)\geqslant d$. By assumption, $\ms M$ satisfies Stacky Murphy's Law, so $\ms G$ is isomorphic to the residual gerbe of a point of $\ms M$, which is contained in the substack $\ms M^G$. The inclusion $\ms M^G\subset\ms M$ is a monomorphism, so by Lemma \ref{lem:immersions preserve residual gerbes} $\ms G$ is isomorphic to the residual gerbe of a point of $\ms M^G$. The inclusion of the residual gerbe of a point is a monomorphism, and in particular is representable, so by Lemma \ref{lem:representation dimension lemma} we have $d\leqslant \rdim(\ms G)\leqslant\rdim(\ms M^G)$. We conclude that $\rdim(\ms M^G)=\infty$.
\end{proof}

\appendix

\section{Algebraic stacks and residual gerbes}\label{sec:gerbes and stacks}

In this appendix we give some background on algebraic stacks. These results are all well known, and are included only for lack of a suitable reference. The contents of this section are as follows. In \S\ref{ssec:points of algebraic stacks} we define the underlying set of points of an algebraic stack and summarize the properties of their associated residual gerbes and residue fields. We give the definitions of these objects and the proofs in \S\ref{ssec:construction of residual gerbes and residue fields}. We have opted here to take a path which differs from what we have seen in the literature; see Remark \S\ref{rem:comparison of approaches to res gerbes} for a comparison. We then explain the relationship between our notion of the residue field and the classical notion of the field of moduli \S\ref{ssec:field of moduli}, and describe the structure of residual gerbes and residue fields in the presence of a coarse moduli space \S\ref{ssec:coarse moduli spaces}. Finally, in \S\ref{ssec:stabilizer preserving morphisms and res gerbes} we record a few results describing the action of a morphism of stacks on residual gerbes and residue fields.

\subsection{Points of algebraic stacks}\label{ssec:points of algebraic stacks}

Let $\ms M$ be an algebraic stack. We will define the set of \emph{topological points} $|\ms M|$ of $\ms M$ and the residue field and residual gerbe associated to a point $x\in|\ms M|$.

\begin{Definition}\cite[04XE, 04XG]{stacks-project}
\label{def:topological points}
A \textit{topological point} of $\ms M$ is an equivalence class of pairs $(L,x_L)$, where $L$ is a field and $x_L\in\ms M(L)$, and where two such pairs $(L,x_L)$ and $(L',x_{L'})$ are declared equivalent if there exists a field $F$ and inclusions $L\subset F$ and $L'\subset F$ such that $x_L|_F\cong x_{L'}|_F$ as objects of $\ms M(F)$, or in other words if there exists a field $F$ and a 2-commutative diagram
\begin{equation}\label{eq:points diagram}
        \begin{tikzcd}[column sep=small, row sep=small]
            &\Spec L\arrow{dr}{x_L}&\\
            \Spec F\arrow{ur}\arrow{dr}&&\ms M\\
            &\Spec L'.\arrow{ur}[swap]{x_{L'}}&
        \end{tikzcd}
\end{equation}
We let $|\ms M|$ denote the set of topological points of $\ms M$. If we think confusion is unlikely, we may refer to an element of $|\ms M|$ simply as a \emph{point}.
\end{Definition}

A morphism $f:\ms M\to\ms N$ of algebraic stacks induces a map $|f|:|\ms M|\to |\ms N|$. The set $|\ms M|$ can be equipped with a canonical topology \cite[04XL]{stacks-project}, with respect to which the map on sets of topological points induced by a morphism of stacks is continuous. Furthermore, if $\ms M$ is a scheme, then $|\ms M|$ is the set of points of the underlying topological space of $\ms M$ in the usual sense, and moreover the canonical topology on $|\ms M|$ reduces to the usual Zariski topology.

\begin{Remark}
    The inclusion $\ms M_{\red}\subset\ms M$ of the reduced substack of $\ms M$ induces a bijection $|\ms M_{\red}|\iso |\ms M|$ on sets of topological points.
\end{Remark}

\begin{Example}
    If $\ms G$ is a gerbe over a field $K$, then $|\ms G|$ is a singleton and the induced map $|\ms G|\to |\Spec K|$ is a bijection.
\end{Example}

Fix a point $x\in |\ms M|$. 

\begin{Definition}
    We say that an $L$-point $x_L\in\ms M(L)$ \textit{represents} $x$ if $x_L$ is in the equivalence class $x$. We say that $x$ is \textit{defined over} $L$ or that $L$ is a \textit{field of definition} for $x$ if there exists an $L$-point of $\ms M$ representing $x$.
\end{Definition}

We are interested in the collection of possible fields $L$ and $L$-points $x_L$ representing $x$. If $\ms M$ is a scheme, then this collection is entirely controlled by the residue field of $x$. Indeed, the canonical inclusion $\Spec k(x)\hookrightarrow\ms M$ of the spectrum of the residue field of $x$ is terminal among all field--valued points representing $x$. Thus, if $L$ is a field, then there is a bijection between $L$-points $x_L\in\ms M(L)$ representing $x$ and inclusions $k(x)\subset L$. In particular, $x$ is defined over $L$ if and only if $L$ contains $k(x)$.

Suppose now that $\ms M$ is an algebraic stack. Then there need not exist such a terminal object associated to a point $x\in|\ms M|$. Instead, we will define two objects associated to $x$ which approximate this property in different ways: the \textit{residual gerbe} $\ms G(x)$ of $x$, a stack which receives a map from every representative of $x$, but is not itself the spectrum of a field, and the \textit{residue field} $k(x)$ of $x$, which is a field contained in every field of definition for $x$, but whose spectrum does not itself map to $\ms M$.
These objects fit into a diagram
\begin{equation}\label{eq:residue field and gerbe}
    \begin{tikzcd}
        \ms G(x)\arrow{d}[swap]{\rho}\arrow[hook]{r}{\iota}&\ms M\\
        \Spec k(x).&
    \end{tikzcd}    
\end{equation}

We summarize the key properties of these objects. To ensure good behavior, we assume that $\ms M$ has quasi-compact diagonal. We then have the following.

\begin{itemize}
    \item $\ms G(x)$ is an algebraic stack, and is a gerbe over $k(x)$. In particular, $\ms G(x)$ is reduced and $|\ms G(x)|$ is a singleton.
    \item The diagram~\eqref{eq:residue field and gerbe} behaves functorially with respect to morphisms of stacks.
    \item The morphism $\iota$ is a monomorphism, and hence is \emph{stabilizer preserving}, meaning that it induces an isomorphism $\ms I_{\ms G(x)}\iso\iota^{-1}\ms I_{\ms M}$ on inertia.
    \item The map $|\iota|:|\ms G(x)|\hookrightarrow|\ms M|$ induced by $\iota$ sends the unique element of $|\ms G(x)|$ to $x$.
    \item The residual gerbe has the following universal property: if $L$ is a field and $x_L\in\ms M(L)$ is an $L$-point representing $x$, then $x_L$ factors uniquely through $\ms G(x)$, yielding a diagram
\[
    \begin{tikzcd}
        \Spec L\arrow{r}\arrow{dr}\arrow[bend left=35]{rr}{x_L}&\ms G(x)\arrow{d}{\rho}\arrow[hook]{r}{\iota}&\ms M\\
        &\Spec k(x).&
    \end{tikzcd}
\]
\end{itemize}
In particular, the final property implies that if $L$ is any field then there is a bijection between $L$-points representing $x$ and pairs $(i,s)$, where $i:k(x)\hookrightarrow L$ is an inclusion of fields and $s\in(\ms G(x))(L)$ is an $L$-point of $\ms G(x)$ over $k(x)$.
Thus, a necessary and sufficient condition for $x$ to be defined over $L$ is that there exists an inclusion $k(x)\subset L$ such that the field extension $L/k(x)$ splits the residual gerbe $\ms G(x)$. In particular, $x$ is defined over its own residue field if and only if $\ms G(x)$ is split.

\subsection{Construction of the residual gerbe and residue field}\label{ssec:construction of residual gerbes and residue fields}

Let $\ms M$ be an algebraic stack and let $x\in|\ms M|$ be a point. We will now define the residue field and residual gerbe associated to $x$ and prove some of their basic properties. We first define the \emph{category of representatives} of $x$.

\begin{Definition}
    Let $\Rep(x)$ be the category whose objects are pairs $(L,x_L)$, where $L$ is a field and $x_L\in\ms M(L)$ is an $L$-point which represents $x$. A morphism in $\Rep(x)$ between two objects $(L',x_{L'})$ and $(L,x_{L})$ is a 2-commutative diagram
\begin{equation}\label{eq:top points equiv diagram}
    \begin{tikzcd}[row sep=small]
        \Spec L'\arrow{dd}\arrow{dr}{x_{L'}}&\\
        &\ms M\\
        \Spec L.\arrow{ur}[swap]{x_{L}}&
    \end{tikzcd}
\end{equation}
\end{Definition}

The category $\Rep(x)$ is fibered in groupoids over the opposite of the category of fields, with fiber $\Rep(x)(L)$ over a field $L$ given by the full subcategory of $\ms M(L)$ whose objects are those $L$-points which represent $x$.

\begin{Remark}\label{rem:residue field is terminal}
    If $\ms M$ is a scheme, then the inclusion $\iota:\Spec k(x)\hookrightarrow\ms M$ of the spectrum of the residue field of $x$ gives a canonical representative of $x$ through which every other representative factors uniquely. Thus, the pair $(k(x),\iota)\in\Rep(x)$ is the terminal object of $\Rep(x)$. Moreover, $\Rep(x)$ itself is the fibered category associated to the functor on the category of spectra of fields represented by $\Spec k(x)$.
\end{Remark}

We use the category $\Rep(x)$ to define the \emph{residue field} of the point $x\in|\ms M|$.

\begin{Definition}[The residue field of a point]
\label{def:residue field}
    We define the \textit{residue field} $k(x)$ of $x$ as follows. An element of $k(x)$ is a collection
    \[
        \lambda=\left\{\lambda_{(L,x_L)}\right\}_{(L,x_L)\in\Rep(x)}
    \]
    consisting of a choice for each object $(L,x_L)$ of $\Rep(x)$ of an element $\lambda_{(L,x_L)}\in L$, with the property that if $(L',x_{L'})\to (L,x_L)$ is a morphism in $\Rep(x)$ then we have $\lambda_{(L,x_L)}|_{L'}=\lambda_{(L',x_{L'})}$. We define addition and multiplication in $k(x)$ componentwise: given $\lambda,\mu\in k(x)$, we define $\lambda+\mu$ and $\lambda\mu$ by
    \[
        (\lambda+\mu)_{(L,x_L)}:=\lambda_{(L,x_L)}+\mu_{(L,x_L)}\hspace{.5cm}\text{and}\hspace{.5cm}(\lambda\mu)_{(L,x_L)}:=\lambda_{(L,x_L)}\mu_{(L,x_L)}.
    \]
\end{Definition}

It is immediate that $k(x)$ is a ring. We verify in the following lemma that $k(x)$ is in fact a field.

\begin{Lemma}
    The set $k(x)$ equipped with the above addition and multiplication is a field. 
\end{Lemma}
\begin{proof}
    Let $\left\{\lambda_{(L,x_L)}\right\}_{(L,x_L)\in\Rep(x)}$ be an element of $k(x)$. We claim that $\lambda_{(L,x_L)}=0$ for some $(L,x_L)$ if and only if $\lambda_{(L,x_L)}=0$ for all $(L,x_L)$. Indeed, suppose given two elements $$(L,x_L),(L',x_{L'})\in\Rep(x).$$ By assumption, $x_L$ and $x_{L'}$ both represent $x$, so we may find a field $F$ and a 2-commutative diagram as in~\eqref{eq:top points equiv diagram}. Both maps $L\subset F$ and $L'\subset F$ are injective, so $\lambda_{(L,x_L)}=0$ if and only if $\lambda_{(L',x_{L'})}=0$. This proves the claim.
    
    We now show that $k(x)$ is a field. Suppose given an element $\lambda\in k(x)$ which is not identically zero. By the above claim, we then have that $\lambda_{(L,x_L)}\neq 0$ for all $(L,x_L)\in\Rep(x)$, and so we may define the inverse $\lambda^{-1}$ of $\lambda$ by taking the inverse componentwise:
    \[
        (\lambda^{-1})_{(L,x_L)}:=(\lambda_{(L,x_L)})^{-1}.
    \]
\end{proof}

If $(L,x_L)$ is any object of $\Rep(x)$, then the map $\lambda\mapsto\lambda_{(L,x_L)}$ defines a canonical inclusion $k(x)\subset L$. These inclusions are compatible with morphisms in $\Rep(x)$, in the sense that given any morphism in $\Rep(x)$ between objects $(L',x_{L'})$ and $(L,x_L)$, the diagram
    \[
        \begin{tikzcd}
            \Spec L'\arrow{dr}\arrow[bend left=25]{rrd}{x_{L'}}\arrow[bend right=30]{ddr}&&\\
            &\Spec L\arrow{r}{x_L}\arrow{d}&\ms M\\
            &\Spec k(x)&
        \end{tikzcd}
    \]
is 2-commutative.

\begin{Lemma}
    If $\ms M$ is a scheme, then the residue field of $\ms M$ at $x$ as defined in Definition \ref{def:residue field} is naturally identified with the residue field at $x$ in the classical sense.
\end{Lemma}
\begin{proof}
    Let us temporarily write $k(x)$ for the classical residue field at $x$ and $k(x)'$ for the residue field at $x$ as defined above. There is a canonical inclusion $\iota:\Spec k(x)\hookrightarrow\ms M$, so as described above we have a canonical map $k(x)'\hookrightarrow k(x)$. We define an inverse. As noted in Remark \ref{rem:residue field is terminal}, the pair $(k(x),\iota)\in\Rep(x)$ is a terminal object, so each object $(L,x_L)\in\Rep(x)$ is equipped with a unique morphism to $(k(x),\iota)$. Let $i_{(L,x_L)}:k(x)\hookrightarrow L$ denote the corresponding map of fields. The association $f\mapsto\left\{i_{(L,x_L)}(f)\right\}$ defines a field homomorphism $k(x)\hookrightarrow k(x)'$ which is the desired inverse.
\end{proof}

We now define the residual gerbe associated to the point $x\in|\ms M|$. Our definition is essentially that it should be the minimal substack $\ms G(x)\subset\ms M$ through which every representative of $x$ factors.

\begin{Remark}\label{rem:comparison of approaches to res gerbes}
    Our approach to the residual gerbes of a stack differs somewhat from those in the literature. Laumon and Moret-Bailly \cite[\S 11]{MR1771927} define the residual gerbe of a point $x\in|\ms M|$ using a choice of geometric point representing $x$. With this approach, one is then faced with the problem of showing the independence of the construction from this choice (see Rydh \cite[Appendix B]{MR2774654} for a discussion of this issue). The Stacks Project \cite{stacks-project} makes the definition that \emph{the residual gerbe at $x$ exists} \cite[06MU]{stacks-project} if there exists a reduced locally noetherian algebraic stack $\ms Z$ such that $|\ms Z|$ is a singleton and a monomorphism $\ms Z\hookrightarrow\ms M$ with image $x$. It is then shown that if this is the case then there is a unique substack of $\ms M$ which is the image of any such map. We instead define the residual gerbe of a point directly as a certain explicit substack of $\ms M$. We remark that all of these notions agree where they may be reasonably compared, see Lemma \ref{lem:residual gerbe is unique}.
\end{Remark}

\begin{Definition}[The residual gerbe of a point]
\label{def:residual gerbe}
    The \textit{residual gerbe} of $\ms M$ at $x$ is the full subcategory $\ms G(x)\subset\ms M$ whose fiber over a scheme $T$ is the full subgroupoid $\ms G(x)(T)\subset\ms M(T)$ consisting of those $T$-points $t\in\ms M(T)$ such that there exists a 2-commutative diagram
    \[
        \begin{tikzcd}
            U\arrow{d}\arrow{r}&T\arrow{d}{t}\\
            \Spec L\arrow{r}{x_L}&\ms M
        \end{tikzcd}
    \]
    where $U$ is a scheme, $U\to T$ is an fppf cover, $L$ is a field, and $x_L\in\ms M(L)$ is an $L$-point of $\ms M$ representing $x$.
\end{Definition}

It follows from the definition that $\ms G(x)$ is a stack in the fppf topology. What is not clear from this definition is whether $\ms G(x)$ is an algebraic stack (much less a gerbe). We will shortly give conditions under which this is the case. We first record some weaker properties which are always held by $\ms G(x)$.

\begin{Lemma}\label{lem:representable diagonal}
    The diagonal of the stack $\ms G(x)$ is representable by algebraic spaces.
\end{Lemma}
\begin{proof}
    Let $T$ be a scheme and let $t,u:T\to\ms G(x)$ be morphisms. As the map $\ms G(x)\subset\ms M$ is a monomorphism, the natural morphism $T\times_{\ms G(x)}T\to T\times_{\ms M}T$ is an isomorphism. Because $\ms M$ is an algebraic stack, $T\times_{\ms M}T$ is an algebraic space.
\end{proof}

As a consequence of Lemma \ref{lem:representable diagonal}, if $T$ is a scheme, then any morphism $T\to\ms G(x)$ is representable by algebraic spaces.

\begin{Lemma}\label{lem:fpqc cover}
    If $(L,x_L)$ is an object of $\Rep(x)$, then the corresponding morphism $\Spec L\to\ms G(x)$ is represented by fpqc covers of algebraic spaces.
\end{Lemma}
\begin{proof}
    Let $T$ be a scheme and let $t:T\to\ms M$ be a morphism which factors through $\ms G(x)$. We need to verify that $T\times_{\ms G(x)}\Spec L=T\times_{\ms M}\Spec L\to T$ is an fpqc cover. Choose an fppf cover $U\to T$, an object $(L',x_{L'})$ of $\Rep(x)$, and a 2-commutative diagram
    \[
        \begin{tikzcd}
            U\arrow{d}\arrow{r}&T\arrow{d}{t}\\
            \Spec L'\arrow{r}{x_{L'}}&\ms M.
        \end{tikzcd}
    \]
    As $U\to T$ is an fppf cover, it will suffice to show that $U\times_{\ms M}\Spec L\to U$ is an fpqc cover. Consider the diagram
    \[
        \begin{tikzcd}
            U\times_{\ms M}\Spec L\arrow{d}\arrow{r}&\Spec L'\times_{\ms M}\Spec L\arrow{d}\arrow{r}&\Spec L\arrow{d}\\
            U\arrow{r}&\Spec L'\arrow{r}&\ms M
        \end{tikzcd}
    \]
    which has 2-cartesian squares. As $x_{L'}$ represents $x$, the fiber product $\Spec L'\times_{\ms M}\Spec L$ is nonempty. As $L'$ is a field, the middle vertical arrow is an fpqc cover, and it follows that $U\times_{\ms M}\Spec L\to U$ is an fpqc cover.
\end{proof}

By construction, we have a canonical monomorphism $\iota:\ms G(x)\hookrightarrow\ms M$. To complete the diagram~\eqref{eq:residue field and gerbe}, it remains to define the morphism $\rho:\ms G(x)\to\Spec k(x)$. We consider the \textit{structure sheaf} $\ms O_{\ms G(x)}$ of $\ms G(x)$, which is the sheaf of rings on the category of schemes over $\ms G(x)$ given by $T\mapsto\Gamma(T,\ms O_T)$. An element of the ring $\Gamma(\ms G(x),\ms O_{\ms G(x)})$ of global sections is an assignment to each scheme $T$ over $\ms G(x)$ of an element of $\Gamma(T,\ms O_T)$, these assignments being furthermore required to be compatible with pullbacks in the category of schemes over $\ms G(x)$. Restriction defines a map $\Gamma(\ms G(x),\ms O_{\ms G(x)})\to k(x)$ of rings.

\begin{Lemma}\label{lem:iso on global sections}
    The restriction map $\Gamma(\ms G(x),\ms O_{\ms G(x)})\to k(x)$ is an isomorphism.
\end{Lemma}
\begin{proof}
    Let $(L,x_L)$ be an object of $\Rep(x)$. By Lemma \ref{lem:fpqc cover}, the corresponding morphism $\Spec L\to\ms G(x)$ is representable by fpqc covers of algebraic spaces. Thus, if we set $Z=\Spec L\times_{\ms M}\Spec L$, let $\pi_1,\pi_2:Z\to\Spec L$ denote the two projections, and put $R=\Gamma(Z,\ms O_Z)$, then $\Gamma(\ms G(x),\ms O_{\ms G(x)})$ may be identified with the equalizer of the two homomorphisms
    \[
        \pi_1^{\#},\pi_2^{\#}:L\to R.
    \]
    On the other hand, we have a canonical inclusion $k(x)\subset L$, and with the above identification we have inclusions $\Gamma(\ms G(x),\ms O)\subset k(x)\subset L$. We will show that the left inclusion is an equality. We note that as $L$ is a field, the image of $L$ in $R$ has trivial intersection with the nilradical of $R$, and so $\Gamma(\ms G(x),\ms O)$ is also the equalizer of the compositions of $\pi_1^{\#}$ and $\pi_2^{\#}$ with the reduction map $R\to R_{\red}$. Let $P$ be the set of minimal primes of $R_{\red}$. There is a canonical embedding
    \[
        R_{\red}\hookrightarrow\prod_{\mathfrak{p}\in P}(R_{\red})_{\mathfrak{p}}
    \]
    of $R$ into a product of fields \cite[00EW]{stacks-project}, and hence $\Gamma(\ms G(x),\ms O)$ is the equalizer of the two homomorphisms
    \[
        \begin{tikzcd}
            L\arrow{r}{\pi_i^{\#}}&R\arrow{r}&R_{\red}\arrow{r}&\prod_{\mathfrak{p}\in P}(R_{\red})_{\mathfrak{p}}.
        \end{tikzcd}
    \]
    If $F$ is any field and $f:\Spec F\to Z$ is a morphism, then $k(x)$ is contained in the equalizer of the two homomorphisms $f^{\#}\circ\pi_1^{\#},f^{\#}\circ\pi_2^{\#}:L\to F$. In particular, the above presentation of $\Gamma(\ms G(x),\ms O_{\ms G(x)})$ shows that we have the reverse inclusion $k(x)\subset\Gamma(\ms G(x),\ms O_{\ms G(x)})$.
\end{proof}

\begin{Definition}
    We define the map $\rho:\ms G(x)\to\Spec k(x)$ to be the composition of the canonical affinization morphism $\ms G(x)\to\Spec\Gamma(\ms G(x),\ms O_{\ms G(x)})$ with the inverse of the isomorphism $\Spec k(x)\iso\Spec \Gamma(\ms G(x),\ms O_{\ms G(x)})$ induced by restriction.
\end{Definition}

\begin{Remark}
    The morphism $\rho:\ms G(x)\to\Spec k(x)$ is initial with respect to maps from $\ms G(x)$ to affine schemes. If $\ms G(x)$ is a gerbe over $k(x)$, then $\rho$ is initial with respect to maps from $\ms G(x)$ to algebraic spaces.
\end{Remark}

We now prove some properties of the residual gerbe and residue field. We first observe that the diagram~\eqref{eq:residue field and gerbe} behaves functorially with respect to morphisms of stacks, in the following sense. Let $f:\ms M\to\ms N$ be a morphism of algebraic stacks, let $x\in |\ms M|$ be a topological point, and let $y=f(x)\in |\ms N|$ be the image of $x$ in $|\ms N|$. There is then an induced commutative diagram
    \[
        \begin{tikzcd}[row sep=small, column sep=tiny]
            \ms G(x)\arrow{dd}\arrow{dr}\arrow[hook]{rr}&&\ms M\arrow{dr}{f}&\\
            &\ms G(y)\arrow[hook]{rr}\arrow{dd}&&\ms N\\
            \Spec k(x)\arrow{dr}&&\phantom{\Spec k(x)}&\\
            &\Spec k(y).&&\phantom{\Spec k(x)}
        \end{tikzcd}
    \]
By construction, the residual gerbe $\ms G(x)$ has the universal property of receiving a map from every representative of $x$. We strengthen this slightly in the following result.

\begin{Proposition}\label{prop:universal property of the residual gerbe}
    Let $\ms G$ be an algebraic stack which admits an fppf cover by the spectrum of a field (for instance, a gerbe over a field). If $f:\ms G\to\ms M$ is a morphism whose induced map $|\ms G|\to|\ms M|$ sends the unique element of $|\ms G|$ to $x$, then $f$ factors through $\ms G(x)$.
\end{Proposition}
\begin{proof}
    By assumption, there exists a field $L$ and an fppf cover $\Spec L\to\ms G$. Let $T$ be a scheme and $T\to\ms G$ be a morphism. Form the pullback $T'=T\times_{\ms G}\Spec L$, and consider the diagram
    \[
        \begin{tikzcd}
            T'\arrow{r}\arrow{d}&\Spec L\arrow{d}\arrow{dr}{x_L}\\
            T\arrow{r}&\ms G\arrow{r}{f}&\ms M.
        \end{tikzcd}
    \]
    The map $T'\to T$ is an fppf cover, and the composition $x_L:\Spec L\to\ms M$ represents $x$, so the composition $T\to\ms M$ factors through $\ms G(x)$.
\end{proof}

\begin{Lemma}\label{lem:residual gerbe is unique}
    In the situation of Proposition \ref{prop:universal property of the residual gerbe}, if the map $f:\ms G\to\ms M$ is a monomorphism, then $f$ factors through an isomorphism $\ms G\iso\ms G(x)$.
\end{Lemma}
\begin{proof}
    It will suffice to show that if $L$ is a field then any $L$-point $x_L\in\ms M(L)$ representing $x$ factors through $\ms G$. To see this, consider the pullback diagram
    \[
        \begin{tikzcd}
            T\arrow[hook]{r}\arrow{d}&\Spec L\arrow{d}{x_L}\\
            \ms G\arrow[hook]{r}&\ms M.
        \end{tikzcd}
    \]
    We know that there exists a map to $\ms G$ from the spectrum of a field. Such a map necessarily represents $x$, as does $x_L$. Choosing a common field extension compatible with the maps to $\ms M$ gives a map from a nonempty scheme to $T$. In particular, $T$ is nonempty. As $\ms G\hookrightarrow\ms M$ is a monomorphism, so is the map $T\hookrightarrow\Spec L$. By \cite[06MG]{stacks-project}, $T\hookrightarrow\Spec L$ is an isomorphism, so $x_L$ factors through $\ms G$.
\end{proof}

The following result of Rydh \cite[Theorem B.2]{MR2774654} (see also \cite[06RD]{stacks-project}) shows that if $\ms M$ has quasi-compact diagonal then the residual gerbe at any point is an algebraic stack and is a gerbe over $k(x)$.

\begin{Theorem}\label{thm:algebraicity of residual gerbes}
    If $\ms M$ is an algebraic stack with quasi-compact diagonal, then for every point $x\in|\ms M|$ the residual gerbe $\ms G(x)$ is an algebraic stack which is a gerbe over $k(x)$.
\end{Theorem}
\begin{proof}
    Rydh shows \cite[Theorem B.2]{MR2774654} that for any $x\in|\ms M|$ there exists an algebraic stack $\ms G$ which is a gerbe over a field and a monomorphism $\ms G\hookrightarrow\ms M$ with topological image $x$. By Lemma \ref{lem:residual gerbe is unique}, this map factors through an isomorphism $\ms G\iso\ms G(x)$, and thus $\ms G(x)$ is algebraic and is a gerbe over a field. This implies that the sheafification of $\ms G(x)$ is given by the canonical map $\ms G(x)\to\Spec\Gamma(\ms G(x),\ms O_{\ms G(x)})$, which by construction is isomorphic to the map $\ms G(x)\to\Spec k(x)$.
\end{proof}

\subsection{The field of moduli}\label{ssec:field of moduli}

We will show that the residue field of a point of an algebraic stack is naturally identified with the field of moduli of a representative. Suppose that $\ms M$ is an algebraic stack defined over a field $K_0$. Let $L$ be a separably closed field and let $L/K_0$ be a field extension. Let $x_{L}\in\ms M(L)$ be an $L$-point of $\ms M$. Let $x\in |\ms M|$ be the corresponding topological point. Suppose that $\ms G(x)$ is algebraic and is a gerbe over $k(x)$ (this is the case if $\ms M$ has quasi-compact diagonal, by Theorem \ref{thm:algebraicity of residual gerbes}). In this situation we obtain two naturally defined fields intermediate to the extension $L/K_0$. On the one hand, applying the universal property of the residual gerbe and of the sheafification morphism $\ms G(x)\to\Spec k(x)$, we obtain a 2-commutative diagram
    \[
        \begin{tikzcd}
            \Spec L\arrow{r}\arrow{dr}[swap]{i^*}\arrow[bend left=25]{rr}{x_{L}}&\ms G(x)\arrow{d}{\rho}\arrow[hook]{r}&\ms M\arrow{d}\\
            &\Spec k(x)\arrow{r}&\Spec K_0.
        \end{tikzcd}
    \]
In particular, we have field extensions $L/k(x)/K_0$. On the other hand, given an automorphism $\sigma\in\Gal(L/K_0)$, let us write $(x_{L})^{\sigma}=x_{L}\circ\sigma^*\in\ms M(L)$ for the pullback of $x_{L}$ by $\sigma$, so that we have a 2-commutative diagram
\[
    \begin{tikzcd}
        \Spec L\arrow{r}[swap]{\sigma^*}\arrow[bend left=25]{rr}{(x_{L})^{\sigma}}&\Spec L\arrow{r}[swap]{x_{L}}&\ms M.
    \end{tikzcd}
\]
 We define a subgroup $\Gamma\subset\Gal(L/K_0)$ by
\[
    \Gamma=\left\{\sigma\in\Gal(L/K_0)|(x_{L})^{\sigma}\cong x_{L}\text{ in }\ms M(L)\right\}.
\]

\begin{Definition}\label{def:field of moduli}
    The \textit{field of moduli} of $x_{L}$ (over the ground field $K_0$) is the subfield $L^{\Gamma}\subset L$ consisting of those elements of $L$ which are fixed by every element of $\Gamma$.
\end{Definition}

By construction, we have extensions $L/L^{\Gamma}/K_0$. The following result shows that under suitable conditions the subfields $k(x)$ and $L^{\Gamma}$ of $L$ are equal.

\begin{Proposition}\label{prop:residue field is equal to field of moduli}
    Suppose that $L/k(x)$ is separable and algebraic (this is automatic if, for instance, $L$ is a separable closure of the ground field $K_0$) and that the residual gerbe $\ms G(x)$ has smooth inertia. The residue field $k(x)$ and the field of moduli $L^{\Gamma}$ are equal as subfields of $L$.
\end{Proposition}
\begin{proof}
    By the Galois correspondence and our assumption that $L/k(x)$ is separable, it will suffice to show that if $\sigma\in\Gal(L/K_0)$ is a field automorphism then $(x_L)^{\sigma}\cong x_L$ if and only if $\sigma$ fixes $k(x)$ pointwise. Suppose that $\sigma\in\Gamma$ is an automorphism such that $(x_L)^{\sigma}\cong x_L$. Such an isomorphism gives a 2-isomorphism rendering the diagram
    \[
        \begin{tikzcd}
            &\ms G(x)\arrow{d}{\rho}\\
            \Spec L\arrow[bend left=20]{ur}{(x_L)^{\sigma}}\arrow[bend right=20]{ur}{x_L}\arrow{r}[swap]{i^*}&\Spec k(x)
        \end{tikzcd}
    \]
    2-commutative, where $i^*:\Spec L\to\Spec k(x)$ is the map induced by $x_L$ and $i:k(x)\hookrightarrow L$ is the corresponding field extension. Thus, the compositions $\rho\circ (x_L)^{\sigma}$ and $\rho\circ x_L$ are equal. The former is equal to $i\circ\sigma$, while the latter is equal $i$. It follows that $i\circ\sigma=i$, so $k(x)$ is fixed by $\sigma$.

    Conversely, suppose that $\sigma\in\Gal(L/K_0)$ is a field automorphism which fixes $k(x)$ pointwise. We then have $i\circ\sigma=i$. Thus, $(x_L)^{\sigma}$ and $x_L$ correspond to two sections of the gerbe $\rho:\ms G(x)\to\Spec k(x)$ over the $k(x)$-scheme $i:\Spec L\to\Spec k(x)$. 
    We claim that there exists an isomorphism $(x_L)^{\sigma}\cong x_L$ in the groupoid $\ms G(x)(L)$. Consider the isomorphism sheaf $\sIsom_L((x_L)^{\sigma},x_L)$. As $\ms G(x)$ is a gerbe, this sheaf is nonempty. It also admits a simply transitive action by the smooth group scheme $\sAut_L(x_L)$, and thus is a torsor under $\sAut_L(x_L)$. A torsor under a smooth group scheme is itself smooth, and so is \'{e}tale locally trivial. But $L$ is separably closed, so any \'{e}tale morphism to $\Spec L$ admits a section. Thus, the torsor $\sIsom_L((x_L)^{\sigma},x_L)$ has an $L$-point, so there exists an isomorphism $(x_L)^{\sigma}\cong x_L$, as claimed.
\end{proof}

\begin{Corollary}
    With the assumptions of Proposition \ref{prop:residue field is equal to field of moduli}, the object $x_L\in\ms M(L)$ is defined over an intermediate extension $L/K/L^{\Gamma}$ if and only if the residual gerbe $\ms G(x)$ splits over $K$. In particular, $x_L$ is defined over its field of moduli if and only if the residual gerbe $\ms G(x)$ splits.
\end{Corollary}
\begin{proof}
    If $x_L$ is defined over an intermediate extension $L/K/L^{\Gamma}$, then the by the universal property of the residual gerbe the corresponding $K$-point gives a splitting of $\ms G(x)$ over $K$. Conversely, suppose that $x_K:\Spec K\to\ms G(x)$ is a splitting over $K$. Write $x_L'$ for the pullback of $x_K$ to $L$. As in the proof of Proposition \ref{prop:residue field is equal to field of moduli}, any two sections of $\ms G(x)$ over $L$ are isomorphic, so $x_L'$ is isomorphic to $x_L$.
\end{proof}

\subsection{Coarse moduli spaces}\label{ssec:coarse moduli spaces}

In this section we will consider the residual gerbes and residue fields of an algebraic stack in the presence of a coarse space morphism. In particular, we will show that the residual gerbes of points can be recovered as the reductions of the fibers of the coarse moduli space morphism. Let $\ms M$ be an algebraic stack.

\begin{Definition}
    A \textit{coarse moduli space} for $\ms M$ consists of an algebraic space $M$ and a morphism $\rho:\ms M\to M$ which is universal for maps from $\ms M$ to algebraic spaces and which induces a bijection $|\ms M|\iso |M|$ (equivalently, which induces a bijection $\pi_0(\ms M(L))\iso M(L)$ for every algebraically closed field $L$).
\end{Definition}

\begin{Remark}\label{rem:coarse moduli space remark}
    By a fundamental result of Keel--Mori \cite{MR1432041} (since generalized by Conrad \cite{conradstacks} and Rydh \cite{MR3084720}), if $\ms M$ is an algebraic stack which has finite diagonal (eg. a separated Deligne--Mumford stack) then $\ms M$ admits a coarse moduli space $\rho:\ms M\to M$. This has the additional properties that the map $\rho$ is proper and the pullback map $\rho^{\#}:\ms O_M\to\rho_*\ms O_{\ms M}$ is an isomorphism. Furthermore, $\rho$ is compatible with flat base change, in the sense that if $U\to M$ is a flat morphism of algebraic spaces then the base change $\ms M\times_MU\to U$ is a coarse space morphism.
\end{Remark}

\begin{Example}
    If $\ms G$ is a gerbe over a field $K$, then the structural morphism $\ms G\to\Spec K$ is a coarse moduli space for $\ms G$.
\end{Example}

Let $\ms M$ be a Deligne--Mumford stack with finite diagonal which is of finite type over a locally noetherian scheme $S$ and let $\rho:\ms M\to M$ be its coarse moduli space. Let $x\in |\ms M|$ be a point, let $\ms G(x)$ and $k(x)$ be the residual gerbe and residue field of $\ms M$ at $x$, and let $x'=\rho(x)\in|M|$ be the image of $x$ in $M$. The canonical morphism $\ms G(x)\to\Spec k(x)$ is a coarse space morphism, so by the universal property there is an induced map $\Spec k(x)\to M$ and a commutative diagram
\[
    \begin{tikzcd}
        \ms G(x)\arrow{d}\arrow[hook]{r}&\ms M\arrow{d}{\rho}\\
        \Spec k(x)\arrow{r}&M.
    \end{tikzcd}
\]
The lower horizontal arrow has image equal to the point $x'$, and hence factors through the inclusion $\Spec k(x')\hookrightarrow M$ of the residue field of the point $x'\in|M|$. Thus, we obtain a field extension $k(x)/k(x')$.

\begin{Proposition}\label{prop:residue field with tame assumption}
    The extension $k(x)/k(x')$ is purely inseparable, and if $\ms M$ is tame at $x$, then $\rho$ induces an isomorphism $k(x')\iso k(x)$. In particular, if $\ms M$ is tame (eg. if $\ms M$ is defined over a field of characteristic 0) then the bijection $|\rho|:|\ms M|\to|M|$ preserves residue fields.
\end{Proposition}
\begin{proof}
    This follows from \cite[Theorem 11.3.6]{MR3495343}.
\end{proof}

Let $\ms M_{x'}=\Spec k(x')\times_M\ms M$ be the fiber of $\rho$ over the canonical inclusion $\Spec k(x')\subset M$ of the spectrum of the residue field of the point $x'\in |M|$. By the universal property of the 2-fiber product, we obtain a morphism $\ms G(x)\to\ms M_{x'}$. As $\ms G(x)$ is reduced, this map factors through the reduction, yielding a map $\tau:\ms G(x)\to(\ms M_{x'})_{\red}$ and a commutative diagram
\[
    \begin{tikzcd}
        \ms G(x)\arrow{r}{\tau}\arrow{d}&(\ms M_{x'})_{\red}\arrow{d}\arrow[hook]{r}&\ms M\arrow{d}{\rho}\\
        \Spec k(x)\arrow{r}&\Spec k(x')\arrow[hook]{r}&M.
    \end{tikzcd}
\]
\begin{Lemma}\label{lem:residual gerbes and the coarse space}
    The map $\tau:\ms G(x)\to(\ms M_{x'})_{\red}$ is an isomorphism.
\end{Lemma}
\begin{proof}
We observe that the stack $(\ms M_{x'})_{\red}$ is reduced and that $|(\ms M_{x'})_{\red}|=\left\{x\right\}$ is a singleton. The morphism $\rho$ is finite, so we can find a field $L$ and a finite morphism $\Spec L\to(\ms M_{x'})_{\red}$. This morphism is automatically flat \cite[06MM,06MN]{stacks-project}, and hence is an fppf cover. The topological image of the inclusion $(\ms M_{x'})_{\red}\subset\ms M$ is equal to $x$, so by the universal property of the residual gerbe (Proposition \ref{prop:universal property of the residual gerbe}) the inclusion factors through a morphism $(\ms M_{x'})_{\red}\to\ms G(x)$, which is an inverse to $\tau$.
\end{proof}

\begin{Remark}
    Let us consider the situation in which $\ms M$ is tame and is given as a moduli stack parameterizing some algebro-geometric objects. There is then a universal family over $\ms M$, but this family need not descend to a family over the coarse space $M$. The residual gerbes of $\ms M$ provide an obstruction to this being the case. Indeed, let $x\in|\ms M|$ be a point and write $x'\in|M|$ for the image of $x$ under the bijection $|\rho|:|\ms M|\iso|M|$. We obtain a diagram
    \begin{equation}\label{eq:coarse moduli space diagram}
    \begin{tikzcd}
        &\ms G(x)\arrow[hook]{r}{\iota}\arrow{d}\arrow{dl}&\ms M\arrow{d}{\rho}\\
        \Spec k(x)\arrow{r}{\sim}&\Spec k(x')\arrow[dashed]{ur}\arrow[hook]{r}&M
    \end{tikzcd}
    \end{equation} 
    of solid arrows. Suppose that there were a family over the point $\Spec k(x')\subset M$ whose pullback along $\rho$ was isomorphic to the restriction of the universal family. Giving such a family is the same thing as giving a dashed arrow $\Spec k(x')\to\ms M$ rendering the above diagram commutative. The morphism $\iota$ identifies $\ms G(x)$ with the reduction of the fiber product $\Spec k(x')\times_M\ms M$, so any such arrow factors uniquely through a section of the map $\ms G(x)\to\Spec k(x')$. We conclude that the universal family on $\ms M$ descends to a family over the point $\Spec k(x')\subset M$ if and only if the residual gerbe $\ms G(x)$ is split. 
\end{Remark}

\subsection{Stabilizer-preserving morphisms and residual gerbes}\label{ssec:stabilizer preserving morphisms and res gerbes}

In this section we consider the action of a morphism of stacks on residual gerbes. Let $f:\ms M\to\ms N$ be a morphism of algebraic stacks. There is an induced commutative square
\begin{equation}\label{eq:induced square on inertia}
    \begin{tikzcd}
        \ms I_{\ms M}\arrow{r}{\ms I_f}\arrow{d}&\ms I_{\ms N}\arrow{d}\\
        \ms M\arrow{r}{f}&\ms N.
    \end{tikzcd}
\end{equation}

\begin{Definition}\label{def:stabilizer preserving}
\cite[0DU6]{stacks-project}
    A morphism $f:\ms M\to\ms N$ of algebraic stacks is \textit{stabilizer preserving} if the square~\eqref{eq:induced square on inertia} is 2-Cartesian, or equivalently if the induced morphism $\ms I_{\ms M}\to\ms I_{\ms N}\times_{\ms N}\ms M$ is an isomorphism.
\end{Definition}

\begin{Example}\label{ex:stabilizer preserving implies representable}
    A morphism $f:\ms M\to\ms N$ of algebraic stacks is representable if and only if the induced morphism $\ms I_{\ms M}\to\ms I_{\ms N}\times_{\ms N}\ms M$ is injective \cite[04YY]{stacks-project}. In particular, a stabilizer preserving morphism is representable.
\end{Example}

\begin{Example}
    If $f:\ms M\to\ms N$ is a monomorphism \cite[04ZW]{stacks-project} then $f$ is stabilizer preserving \cite[0CBB]{stacks-project}. In particular, an immersion of stacks is stabilizer preserving.
\end{Example}

\begin{Example}\label{ex:inclusion of residual gerbe is stabilizer preserving}
    If $\ms M$ is an algebraic stack and $x\in |\ms M|$ is a point whose residual gerbe $\ms G(x)$ is algebraic, then the canonical inclusion $\ms G(x)\hookrightarrow\ms M$ is a monomorphism of algebraic stacks, and hence is stabilizer preserving.
\end{Example}

Let $\ms G$ be a gerbe over a field $K$ and let $\ms G'$ be a gerbe over a field $L$. Let $G_{\ms G}\to\ms G$ and $G_{\ms G'}\to\ms G'$ be the respective inertia stacks of $\ms G$ and $\ms G'$. Consider a morphism $\tau:\ms G'\to\ms G$ of stacks. This induces a 2-commutative diagram
\begin{equation}\label{eq:square of gerbes}
    \begin{tikzcd}
        G_{\ms G'}\arrow{d}\arrow{r}{\zeta}&G_{\ms G}\arrow{d}\\
        \ms G'\arrow{d}\arrow{r}{\tau}&\ms G\arrow{d}\\
        \Spec L\arrow{r}{\epsilon}&\Spec K
    \end{tikzcd}
\end{equation}
of algebraic stacks.

\begin{Lemma}\label{lem:morphism of gerbes}
    We have the following.
    \begin{enumerate}
        \item\label{item:gerbe maps 1} $\tau$ is representable if and only if the map $G_{\ms G'}\to G_{\ms G}\times_{\ms G}\ms G'$ of group spaces over $\ms G'$ induced by $\zeta$ is injective.
        \item\label{item:gerbe maps 2} $\tau$ is stabilizer-preserving if and only if the lower square of the above diagram~\eqref{eq:square of gerbes} is 2-Cartesian (equivalently, the induced map $\ms G'\to\ms G\otimes_KL$ of gerbes over $L$ is an isomorphism).
        \item\label{item:gerbe maps 3} $\tau$ is a monomorphism if and only if the maps $\tau,\epsilon,$ and $\zeta$ are all isomorphisms.
    \end{enumerate}
\end{Lemma}
\begin{proof}
    A morphism $f:\ms M\to\ms N$ of stacks is representable if and only if the induced map $\ms I_{\ms M}\to\ms I_{\ms N}\times_{\ms N}\ms M$ is injective \cite[04YY]{stacks-project}. This gives $(1)$.
    
    For $(2)$, suppose that $\tau$ is stabilizer-preserving. By definition, this means that the upper square of the diagram~\eqref{eq:square of gerbes} is 2-Cartesian. The map $\ms G'\to\ms G\otimes_KL$ is then a map of gerbes which induces an isomorphism on inertia, and is therefore automatically an isomorphism. Conversely, the property of being stabilizer-preserving is preserved under base change. The map $\epsilon$ is stabilizer preserving, so if the lower square of the diagram is 2-Cartesian, then $\tau$ is also stabilizer preserving.

    For $(3)$, suppose that $\tau$ is a monomorphism. Then $\tau$ is in particular stabilizer-preserving, so by $(2)$, both squares of the diagram are 2-Cartesian. The morphism $\ms G\to\Spec K$ admits a section over a finite extension of $K$. This shows that $\epsilon$ is a monomorphism fppf locally on the target, and therefore is itself a monomorphism. By \cite[03DP]{stacks-project}, $\epsilon$ is necessarily an isomorphism. As both squares are 2-Cartesian, this implies that $\tau$ and $\zeta$ are also isomorphisms.
\end{proof}

\begin{Lemma}\label{lem:immersions preserve residual gerbes}
    Let $f:\ms M\hookrightarrow\ms N$ be a monomorphism of algebraic stacks. If $x\in|\ms M|$ is a point with image $y=f(x)\in|\ms N|$, then $f$ induces isomorphisms $\ms G(x)\iso\ms G(y)$ and $\Spec k(x)\iso\Spec k(y)$.
\end{Lemma}
\begin{proof}   
    The morphism $f$ induces a diagram
    \[
        \begin{tikzcd}
            \ms G(x)\arrow{d}\arrow[hook]{r}&\ms M\arrow[hook]{d}{f}\\
            \ms G(y)\arrow[hook]{r}&\ms N.
        \end{tikzcd}
    \]
    As $f$ and the inclusion $\ms G(x)\subset\ms M$ are monomorphisms, the map $\ms G(x)\to\ms G(y)$ is a monomorphism. The result follows from Lemma \ref{lem:morphism of gerbes}~\eqref{item:gerbe maps 3}.
\end{proof}

To verify that certain morphisms are monomorphisms, we will use the following result.

\begin{Lemma}\label{lem:monomorphism lemma}
    Let $f:\ms M\to\ms N$ be a morphism of algebraic stacks. Suppose that 
    \begin{enumerate}
        \item\label{item:unramified} $f$ is unramified (in the sense of \cite[0CIT]{stacks-project}),
        \item\label{item:injective on pts} if $K$ is a field, then the map
            \[
                \pi_0(\ms M(K))\to\pi_0(\ms N(K))
            \]
            is injective, and
        \item\label{item:inertia is inert} given a field $K$ and an object $x\in\ms M(K)$, the induced map
        \[
            \Aut_{\ms M(K)}(x)\to\Aut_{\ms N(K)}(f(x))
        \]
        is an isomorphism.
    \end{enumerate}
    Then $f$ is a monomorphism.
\end{Lemma}
\begin{proof}
    By \cite[04ZZ]{stacks-project}, it will suffice to show that the diagonal $\Delta_f:\ms M\to\ms M\times_{\ms N}\ms M$ is an equivalence. By \cite[0CJ0]{stacks-project}, $f$ is unramified if and only if $f$ is locally of finite type and $\Delta_f:\ms M\to\ms M\times_{\ms N}\ms M$ is \'{e}tale. Condition~\eqref{item:injective on pts} and the surjectivity in condition~\eqref{item:inertia is inert} combined imply that for any field $K$, the diagonal map $\pi_0(\ms M(K))\to\pi_0((\ms M\times_{\ms N}\ms M)(K))$ is surjective. Therefore $\Delta_f$ is an \'{e}tale cover. Consider the 2-Cartesian diagram
    \[
        \begin{tikzcd}
            \ms I_{\ms M/\ms N}\arrow{d}\arrow{r}{\delta}&\ms M\arrow{d}{\Delta_f}\\
            \ms M\arrow{r}{\Delta_f}&\ms M\times_{\ms N}\ms M.
        \end{tikzcd}
    \]
    The morphism $\delta$ is thus also an \'{e}tale cover. Let $K$ be a scheme. The fiber of the relative inertia stack $\ms I_{\ms M/\ms N}$ over a point $x\in\ms M(K)$ is (by definition) the set of pairs $(x,\alpha)$ where $x\in\ms M(K)$ and $\alpha$ is an automorphism of $x$ which is in the kernel of the map $\Aut_{\ms M(K)}(x)\to\Aut_{\ms N}(f(x))$. By~\eqref{item:inertia is inert}, this kernel is trivial, and therefore the map $\delta$ is universally injective. As $\delta$ is representable, it is therefore an equivalence. We conclude that $\Delta_f$ is an equivalence, and therefore $f$ is a monomorphism. 
\end{proof}

\bibliographystyle{plain}
\bibliography{biblio}{}

\end{document}